\documentclass[11pt,leqno]{article}
\normalfont
\renewcommand{\baselinestretch}{1.118}

\usepackage{amsfonts}

\newfont{\eulercursive}{eurm10 at 11pt}

\newcommand{\QED}{\raisebox{0.5mm}{\fbox{\rule{0mm}{1.5mm}\ }}}

\newcounter{myfn}[page]

\newcommand{\DynkinDiagramFigure}{Figure 1.1}
\newcommand{\BtwoMarsDemo}{Figure 1.2} 
\newcommand{\MainResults}{Propositions 2.3 and 3.1}

\newcommand{\StrongConvergenceTheorem}{Theorem 1.1}
\newcommand{\StrongConvergenceCorollary}{Lemma 1.2}
\newcommand{\ComparisonTheorem}{Theorem 1.3}
\newcommand{\ErikssonTheorems}{Theorems 1.1 and 1.3}
\newcommand{\ComparisonCorollary}{Lemma 1.4}
\newcommand{\ComparisonResults}{Lemmas 1.2 and 1.4}
\newcommand{\NotMarsFriendlyLemma}{Lemma 1.5}

\newcommand{\FourFamiliesFirstLemma}{Lemma 2.1}
\newcommand{\FourFamiliesSecondLemma}{Lemma 2.2}
\newcommand{\ConvergenceProposition}{Proposition 2.3}
\newcommand{\RemarkOnGameLength}{Remark 2.4}

\newcommand{\NotMarsFriendlyCatalog}{Proposition 3.1}
\newcommand{\NotMarsFriendlyFigure}{Figure 3.1}

\newcommand{\IntroNum}{1}
\newcommand{\FirstProofNum}{2}
\newcommand{\SecondProofNum}{3}

\newcommand{\myA}{\mbox{\sffamily A}}
\newcommand{\myB}{\mbox{\sffamily B}}
\newcommand{\myC}{\mbox{\sffamily C}}
\newcommand{\myD}{\mbox{\sffamily D}}
\newcommand{\myE}{\mbox{\sffamily E}}
\newcommand{\myF}{\mbox{\sffamily F}}
\newcommand{\myG}{\mbox{\sffamily G}}

\newcommand{\BTwoGraphForFigure}{\setlength{\unitlength}{1in}
\begin{picture}(1,0.4)
\put(0,0){\begin{picture}(1,0)
            \put(0,0.1){\circle*{0.05}}
            \put(-0.11,-0.05){\large $\gamma_{1}$}
            \put(1,0.1){\circle*{0.05}}
            \put(1.05,-0.05){\large $\gamma_{2}$}
            \put(0,0.1){\line(1,0){1}}
            \put(0.2,0.1){\vector(1,0){0.1}}
            \put(0.8,0.1){\vector(-1,0){0.1}}
            \put(0.7,0.1){\vector(-1,0){0.1}}
            \end{picture}}
\end{picture}}

\newcommand{\GTwoMarsOK}{\setlength{\unitlength}{0.75in}
\begin{picture}(1,0.55)
\put(0,0){\begin{picture}(1,0)
            \put(0,0.1){\circle*{0.075}}
            \put(1,0.1){\circle*{0.075}}
            \put(0,0.1){\line(1,0){1}}
            \put(0.2,0.1){\vector(1,0){0.1}}
            \put(0.8,0.1){\vector(-1,0){0.1}}
            \put(0.7,0.1){\vector(-1,0){0.1}}
            \put(0.6,0.1){\vector(-1,0){0.1}}
            \end{picture}}
\end{picture}}

\newcommand{\FFourMarsOK}{\setlength{\unitlength}{0.75in}
\begin{picture}(1,0.85)
\put(0,0){\begin{picture}(1,0)
            \put(0,0.1){\circle*{0.075}}
            \put(1,0.1){\circle*{0.075}}
            \put(2,0.1){\circle*{0.075}}
            \put(3,0.1){\circle*{0.075}}
            \put(0,0.1){\line(1,0){3}}
            \put(1.2,0.1){\vector(1,0){0.1}}
            \put(1.3,0.1){\vector(1,0){0.1}}
            \put(1.8,0.1){\vector(-1,0){0.1}}
            \end{picture}}
\end{picture}}

\newcommand{\CnMarsOK}{\setlength{\unitlength}{0.75in}
\begin{picture}(6,0.55)
\put(0,0){\begin{picture}(1,0)
            \put(0,0.1){\circle*{0.075}}
            \put(1,0.1){\circle*{0.075}}
            \put(2,0.1){\circle*{0.075}}
            \put(4,0.1){\circle*{0.075}}
            \put(5,0.1){\circle*{0.075}}
            \put(6,0.1){\circle*{0.075}}
            \put(0,0.1){\line(1,0){2}}
            \multiput(2,0.1)(0.4,0){5}{\line(1,0){0.2}}
            \put(4,0.1){\line(1,0){2}}
            \put(5.8,0.1){\vector(-1,0){0.1}}
            \put(5.7,0.1){\vector(-1,0){0.1}}
            \put(5.2,0.1){\vector(1,0){0.1}}
           \end{picture}}
\end{picture}}

\newcommand{\BnMarsOK}{\setlength{\unitlength}{0.75in}
\begin{picture}(6,0.55)
\put(0,0){\begin{picture}(1,0)
            \put(0,0.1){\circle*{0.075}}
            \put(1,0.1){\circle*{0.075}}
            \put(2,0.1){\circle*{0.075}}
            \put(4,0.1){\circle*{0.075}}
            \put(5,0.1){\circle*{0.075}}
            \put(6,0.1){\circle*{0.075}}
            \put(0,0.1){\line(1,0){2}}
            \multiput(2,0.1)(0.4,0){5}{\line(1,0){0.2}}
            \put(4,0.1){\line(1,0){2}}
            \put(5.8,0.1){\vector(-1,0){0.1}}
            \put(5.2,0.1){\vector(1,0){0.1}}
            \put(5.3,0.1){\vector(1,0){0.1}}
           \end{picture}}
\end{picture}}

\newcommand{\DnMarsOK}{\setlength{\unitlength}{0.75in}
\begin{picture}(6,0.75)
\put(0,-0.25){\begin{picture}(1,0)
            \put(0,0.35){\circle*{0.075}}
            \put(1,0.35){\circle*{0.075}}
            \put(2,0.35){\circle*{0.075}}
            \put(4,0.35){\circle*{0.075}}
            \put(5,0.35){\circle*{0.075}}
            \put(6,0.1){\circle*{0.075}}
            \put(6,0.6){\circle*{0.075}}
            \put(0,0.35){\line(1,0){2}}
            \multiput(2,0.35)(0.4,0){5}{\line(1,0){0.2}}
            \put(4,0.35){\line(1,0){1}}
            \put(5,0.35){\line(4,1){1}}
            \put(5,0.35){\line(4,-1){1}}
           \end{picture}}
\end{picture}}

\newcommand{\AnMarsOK}{\setlength{\unitlength}{0.75in}
\begin{picture}(6,0.85)
\put(0,0){\begin{picture}(1,0)
            \put(0,0.1){\circle*{0.075}}
            \put(1,0.1){\circle*{0.075}}
            \put(2,0.1){\circle*{0.075}}
            \put(4,0.1){\circle*{0.075}}
            \put(5,0.1){\circle*{0.075}}
            \put(6,0.1){\circle*{0.075}}
            \put(0,0.1){\line(1,0){2}}
            \multiput(2,0.1)(0.4,0){5}{\line(1,0){0.2}}
            \put(4,0.1){\line(1,0){2}}
           \end{picture}}
\end{picture}}

\newcommand{\EEightMarsOK}{\setlength{\unitlength}{0.75in}
\begin{picture}(6,0.75)
\put(0,-0.25){\begin{picture}(1,0)
            \put(0,0.1){\circle*{0.075}}
            \put(1,0.1){\circle*{0.075}}
            \put(2,0.1){\circle*{0.075}}
            \put(2,0.6){\circle*{0.075}}
            \put(3,0.1){\circle*{0.075}}
            \put(4,0.1){\circle*{0.075}}
            \put(5,0.1){\circle*{0.075}}
            \put(6,0.1){\circle*{0.075}}
            \put(0,0.1){\line(1,0){6}}
            \put(2,0.1){\line(0,1){0.5}}
           \end{picture}}
\end{picture}}

\newcommand{\ESevenMarsOK}{\setlength{\unitlength}{0.75in}
\begin{picture}(6,0.75)
\put(0,-0.25){\begin{picture}(1,0)
            \put(0,0.1){\circle*{0.075}}
            \put(1,0.1){\circle*{0.075}}
            \put(2,0.1){\circle*{0.075}}
            \put(2,0.6){\circle*{0.075}}
            \put(3,0.1){\circle*{0.075}}
            \put(4,0.1){\circle*{0.075}}
            \put(5,0.1){\circle*{0.075}}
            \put(0,0.1){\line(1,0){5}}
            \put(2,0.1){\line(0,1){0.5}}
           \end{picture}}
\end{picture}}

\newcommand{\ESixMarsOK}{\setlength{\unitlength}{0.75in}
\begin{picture}(6,0.75)
\put(0,-0.25){\begin{picture}(1,0)
            \put(0,0.1){\circle*{0.075}}
            \put(1,0.1){\circle*{0.075}}
            \put(2,0.1){\circle*{0.075}}
            \put(2,0.6){\circle*{0.075}}
            \put(3,0.1){\circle*{0.075}}
            \put(4,0.1){\circle*{0.075}}
            \put(0,0.1){\line(1,0){4}}
            \put(2,0.1){\line(0,1){0.5}}
           \end{picture}}
\end{picture}}

\newcommand{\TwoCitiesGraphWithLabels}{
\setlength{\unitlength}{1in}
\begin{picture}(1,0.25)
\put(0,0){\begin{picture}(1,0)
            \put(0,0.1){\circle*{0.05}}
            \put(-0.15,-0.05){\large $\gamma_{1}$}
            \put(1,0.1){\circle*{0.05}}
            \put(1.05,-0.05){\large $\gamma_{2}$}
            \put(0,0.1){\line(1,0){1}}
            \put(0.2,0.1){\vector(1,0){0.1}}
            \put(0.8,0.1){\vector(-1,0){0.1}}
            \put(0.225,0){\footnotesize $p$}
            \put(0.71,0){\footnotesize $q$}
            \end{picture}}
\end{picture}}

\newcommand{\ATwoGraphWithLabels}{\setlength{\unitlength}{0.75in}
\begin{picture}(1,0.3)
\put(0.4,0.225){\large $\mbox{\sffamily A}_{2}$}
\put(0,0){\begin{picture}(1,0)
            \put(0,0.1){\circle*{0.05}}
            \put(-0.15,-0.05){\large $\gamma_{1}$}
            \put(1,0.1){\circle*{0.05}}
            \put(1.05,-0.05){\large $\gamma_{2}$}
            \put(0,0.1){\line(1,0){1}}
            \put(0.2,0.1){\vector(1,0){0.1}}
            \put(0.8,0.1){\vector(-1,0){0.1}}
            \end{picture}}
\end{picture}}

\newcommand{\ATwoGraphNoEdgeLabels}{\setlength{\unitlength}{0.75in}
\begin{picture}(1.65,0.25)
\put(0.25,0){\begin{picture}(1,0)
            \put(0,0.1){\circle*{0.05}}
            \put(-0.20,-0.05){\large $\gamma_{1}$}
            \put(1,0.1){\circle*{0.05}}
            \put(1.05,-0.05){\large $\gamma_{2}$}
            \put(0,0.1){\line(1,0){1}}
            \end{picture}}
\end{picture}}

\newcommand{\BTwoGraphWithLabels}{\setlength{\unitlength}{0.75in}
\begin{picture}(1,0.35)
\put(0.4,0.225){\large $\mbox{\sffamily B}_{2}$}
\put(0,0){\begin{picture}(1,0)
            \put(0,0.1){\circle*{0.05}}
            \put(-0.15,-0.05){\large $\gamma_{1}$}
            \put(1,0.1){\circle*{0.05}}
            \put(1.05,-0.05){\large $\gamma_{2}$}
            \put(0,0.1){\line(1,0){1}}
            \put(0.2,0.1){\vector(1,0){0.1}}
            \put(0.8,0.1){\vector(-1,0){0.1}}
            \put(0.7,0.1){\vector(-1,0){0.1}}
            \end{picture}}
\end{picture}}

\newcommand{\GTwoGraphWithLabels}{\setlength{\unitlength}{0.75in}
\begin{picture}(1,0.35)
\put(0.4,0.225){\large $\mbox{\sffamily G}_{2}$}
\put(0,0){\begin{picture}(1,0)
            \put(0,0.1){\circle*{0.05}}
            \put(-0.15,-0.05){\large $\gamma_{1}$}
            \put(1,0.1){\circle*{0.05}}
            \put(1.05,-0.05){\large $\gamma_{2}$}
            \put(0,0.1){\line(1,0){1}}
            \put(0.2,0.1){\vector(1,0){0.1}}
            \put(0.8,0.1){\vector(-1,0){0.1}}
            \put(0.7,0.1){\vector(-1,0){0.1}}
            \put(0.6,0.1){\vector(-1,0){0.1}}
            \end{picture}}
\end{picture}}

\newcommand{\SmallCycles}{\setlength{\unitlength}{0.75in}
\begin{picture}(6.5,3)
\put(-0.5,1.5){\begin{picture}(1,0)
            \put(0,0.6){\circle*{0.075}}
            \put(0.5,0.1){\circle*{0.075}}
            \put(0.5,1.1){\circle*{0.075}}
            \put(0,0.6){\line(1,1){0.5}}
            \put(0,0.6){\line(1,-1){0.5}}
            \put(0.5,0.1){\line(0,1){1}}
            \put(0.5,0.3){\vector(0,1){0.1}}
            \put(0.5,0.9){\vector(0,-1){0.1}}
            \put(0,0.6){\vector(1,1){0.2}}
            \put(0.5,1.1){\vector(-1,-1){0.2}}
            \put(0,0.6){\vector(1,-1){0.2}}
            \put(0.5,0.1){\vector(-1,1){0.2}}
            \put(-0.05,0.8){\footnotesize $q_{1}$}
            \put(0.175,1){\footnotesize $p_{1}$}
            \put(-0.05,0.35){\footnotesize $q_{2}$}
            \put(0.175,0.15){\footnotesize $p_{2}$} 
            \put(0.7,0.6){\footnotesize $p_{1}q_{1} \geq 1, 
            p_{2}q_{2} \geq 1$}
           \end{picture}}
\put(2.3,1.5){\begin{picture}(1,0)
            \put(0,0.6){\circle*{0.075}}
            \put(0.5,0.1){\circle*{0.075}}
            \put(0.5,1.1){\circle*{0.075}}
            \put(0,0.6){\line(1,1){0.5}}
            \put(0,0.6){\line(1,-1){0.5}}
            \put(0.5,0.1){\line(0,1){1}}
            \put(0.5,0.3){\vector(0,1){0.1}}
            \put(0.5,0.9){\vector(0,-1){0.1}}
            \put(0.5,0.8){\vector(0,-1){0.1}}
            \put(0,0.6){\vector(1,1){0.2}}
            \put(0.5,1.1){\vector(-1,-1){0.2}}
            \put(0,0.6){\vector(1,-1){0.2}}
            \put(0.5,0.1){\vector(-1,1){0.2}}
            \put(-0.05,0.8){\footnotesize $q_{1}$}
            \put(0.175,1){\footnotesize $p_{1}$}
            \put(-0.05,0.35){\footnotesize $q_{2}$}
            \put(0.175,0.15){\footnotesize $p_{2}$} 
            \put(0.7,0.6){\footnotesize $p_{1}q_{1} \geq 2, 
            p_{2}q_{2} \geq 2$}
           \end{picture}}
\put(5.15,1.5){\begin{picture}(1,0)
            \put(0,0.6){\circle*{0.075}}
            \put(0.5,0.1){\circle*{0.075}}
            \put(0.5,1.1){\circle*{0.075}}
            \put(0,0.6){\line(1,1){0.5}}
            \put(0,0.6){\line(1,-1){0.5}}
            \put(0.5,0.1){\line(0,1){1}}
            \put(0.5,0.3){\vector(0,1){0.1}}
            \put(0.5,0.9){\vector(0,-1){0.1}}
            \put(0.5,0.8){\vector(0,-1){0.1}}
            \put(0.5,0.7){\vector(0,-1){0.1}}
            \put(0,0.6){\vector(1,1){0.2}}
            \put(0.5,1.1){\vector(-1,-1){0.2}}
            \put(0,0.6){\vector(1,-1){0.2}}
            \put(0.5,0.1){\vector(-1,1){0.2}}
            \put(-0.05,0.8){\footnotesize $q_{1}$}
            \put(0.175,1){\footnotesize $p_{1}$}
            \put(-0.05,0.35){\footnotesize $q_{2}$}
            \put(0.175,0.15){\footnotesize $p_{2}$} 
            \put(0.7,0.6){\footnotesize $p_{1}q_{1} \geq 3$, 
            $p_{2}q_{2} \geq 3$}
           \end{picture}}
\put(0.5,0){\begin{picture}(1,0)
            \put(0,0.6){\circle*{0.075}}
            \put(0.5,0.1){\circle*{0.075}}
            \put(0.5,1.1){\circle*{0.075}}
            \put(1,0.6){\circle*{0.075}}
            \put(0,0.6){\line(1,1){0.5}}
            \put(0,0.6){\line(1,-1){0.5}}
            \put(1,0.6){\line(-1,1){0.5}}
            \put(1,0.6){\line(-1,-1){0.5}}
            \put(0.95,0.65){\vector(-1,1){0.1}}
            \put(0.55,1.05){\vector(1,-1){0.1}}
            \put(0.65,0.95){\vector(1,-1){0.1}}
           \end{picture}}
\put(2,0){\begin{picture}(1,0)
            \put(0,0.6){\circle*{0.075}}
            \put(0.5,0.1){\circle*{0.075}}
            \put(0.5,1.1){\circle*{0.075}}
            \put(1,0.6){\circle*{0.075}}
            \put(0,0.6){\line(1,1){0.5}}
            \put(0,0.6){\line(1,-1){0.5}}
            \put(1,0.6){\line(-1,1){0.5}}
            \put(1,0.6){\line(-1,-1){0.5}}
            \put(0.95,0.65){\vector(-1,1){0.1}}
            \put(0.55,1.05){\vector(1,-1){0.1}}
            \put(0.65,0.95){\vector(1,-1){0.1}}
            \put(0.45,0.15){\vector(-1,1){0.1}}
            \put(0.05,0.55){\vector(1,-1){0.1}}
            \put(0.15,0.45){\vector(1,-1){0.1}}
           \end{picture}}
\put(3.5,0){\begin{picture}(1,0)
            \put(0,0.6){\circle*{0.075}}
            \put(0.5,0.1){\circle*{0.075}}
            \put(0.5,1.1){\circle*{0.075}}
            \put(1,0.6){\circle*{0.075}}
            \put(0,0.6){\line(1,1){0.5}}
            \put(0,0.6){\line(1,-1){0.5}}
            \put(1,0.6){\line(-1,1){0.5}}
            \put(1,0.6){\line(-1,-1){0.5}}
            \put(0.95,0.65){\vector(-1,1){0.1}}
            \put(0.55,1.05){\vector(1,-1){0.1}}
            \put(0.65,0.95){\vector(1,-1){0.1}}
            \put(0.45,0.15){\vector(-1,1){0.1}}
            \put(0.05,0.55){\vector(1,-1){0.1}}
            \put(0.35,0.25){\vector(-1,1){0.1}}
           \end{picture}}
\put(5,0){\begin{picture}(1,0)
            \put(0,0.6){\circle*{0.075}}
            \put(0,0.1){\circle*{0.075}}
            \put(0.5,0.1){\circle*{0.075}}
            \put(0.5,1.1){\circle*{0.075}}
            \put(1,0.6){\circle*{0.075}}
            \put(0,0.6){\line(1,1){0.5}}
            \put(0,0.6){\line(0,-1){0.5}}
            \put(0,0.1){\line(1,0){0.5}}
            \put(1,0.6){\line(-1,1){0.5}}
            \put(1,0.6){\line(-1,-1){0.5}}
            \put(0.95,0.65){\vector(-1,1){0.1}}
            \put(0.55,1.05){\vector(1,-1){0.1}}
            \put(0.65,0.95){\vector(1,-1){0.1}}
           \end{picture}}
\end{picture}}

\begin{document}

\newpage
\setcounter{page}{1} 
\renewcommand{\baselinestretch}{1}

\vspace*{-0.7in}
\hfill {\footnotesize October 29, 2008}

\begin{center}
{\large \bf Convergent and divergent numbers games for certain 
collections\\ of edge-weighted graphs}

Robert G.\ Donnelly

\vspace*{-0.05in} 
Department of Mathematics and Statistics, Murray State
University, Murray, KY 42071

\end{center}

\begin{abstract}
The numbers game is a one-player game played on  
a finite simple 
graph with certain ``amplitudes'' assigned to its edges 
and with an initial assignment of real 
numbers to its nodes.  The moves of the game successively transform 
the numbers at the nodes using the amplitudes in a certain way.  
Here, the edge amplitudes will be negative integers. 
Combinatorial methods are used to investigate the convergence and 
divergence of numbers games played on certain such graphs. 
The results obtained here provide support for results 
in a companion paper. 
\begin{center}

\ 


\end{center}
\end{abstract}
\def\abstractname{Contents}
\begin{abstract}
\begin{center}
\parbox{3.82in}{
\IntroNum. Introduction, definitions, and preliminary results\\
\FirstProofNum.  Convergent numbers games on Dynkin diagrams of 
finite type\\
\SecondProofNum.  Divergent games for some families of graphs}
\end{center}
\end{abstract}

\vspace{2ex} 

\noindent
{\Large \bf \IntroNum.\ \ Introduction, definitions, and preliminary results}

\vspace{1ex} 
The numbers game is a one-player game played on a finite simple graph 
with weights (which we call ``amplitudes'') on its edges 
and with an initial assignment of 
real numbers 
to its nodes.  
Here, each of the two edge amplitudes (one for each direction) 
will be negative integers. 
The move a player can make 
is to ``fire'' one of the nodes with a positive number.  This move 
transforms the number at the fired node 
by changing its sign, and it also 
transforms the number at each adjacent node in a certain way 
using an amplitude   
along the incident edge.  
The player fires the nodes in some sequence of 
the player's choosing, continuing until no node has a positive 
number.  

The numbers game as formulated by Mozes \cite{Mozes}  
has also been studied by Proctor 
\cite{PrEur}, \cite{PrDComplete}, Bj\"{o}rner 
\cite{Bjorner}, 
Eriksson \cite{ErikssonLinear}, \cite{ErikssonThesis}, 
\cite{ErikssonJerusalem}, \cite{ErikssonReachable}, 
\cite{ErikssonDiscrete}, \cite{ErikssonEur}, \cite{DE}, 
Wildberger \cite{WildbergerAdv}, \cite{WildbergerEur}, 
\cite{WildbergerPreprint}, and Donnelly \cite{DonEnumbers}.  
Wildberger studies a dual version which 
he calls the ``mutation game.'' 
See Alon {\em et al} \cite{AKP} for a 
brief and readable treatment of the numbers game on ``unweighted'' 
cyclic graphs. Much of the numbers game discussion in \S 4.3 of 
the book \cite{BB} by Bj\"{o}rner and Brenti can be 
found in \cite{ErikssonThesis} and \cite{ErikssonDiscrete}.  
See these references for discussions of how the numbers game is a 
combinatorial encoding of information for geometric representations of 
Weyl groups (and more generally Coxeter groups) and has uses for 
computing orbits, finding reduced 
decompositions of Weyl group elements,  
solving the word problem, and obtaining combinatorial 
models for Weyl groups. 
Proctor developed this process in \cite{PrEur} 
to compute Weyl group 
orbits of weights with respect to the fundamental weight basis.  Here 
we use his perspective of firing 
nodes with positive, as opposed to negative, numbers. 
Mozes studied numbers games on graphs for 
which the matrix $M$ of integer amplitudes is ``symmetrizable'' 
(i.e.\ there is a nonsingular  
diagonal matrix $D$ such that $D^{-1}M$ is symmetric); 
in \cite{Mozes} he obtained ``strong convergence'' results and a 
geometric characterization of the initial positions for which the 
game terminates.  There will be no symmetrizable assumption here. 

Such graphs-with-amplitudes will henceforth be called ``GCM graphs'' 
for reasons explained below.  
Given any such graph, 
an initial ``position'' is an assignment of numbers to the nodes. 
The position is 
is ``nonzero'' if at least one of the numbers is nonzero.  
A numbers game played from some initial position is ``convergent'' if 
it terminates after a finite number of node 
firings; otherwise we say the game is ``divergent.''  

Here we investigate convergence and divergence of numbers games 
played on certain GCM graphs. 
The purpose is to provide supporting details for the proof of a 
result of \cite{DE}. 
In particular, we aim 
to give straightforward combinatorial proofs of \MainResults. 
These results are used in the proof of the first main result of 
\cite{DE}: A connected GCM graph has a convergent numbers 
game played from a nonzero initial position with nonnegative 
numbers if and only if the graph is one of the ``Dynkin diagrams'' 
of \DynkinDiagramFigure, in which case all numbers games played from 
a given initial position will converge to the same terminal position 
in the same number of steps. 
\ConvergenceProposition\ asserts that 
for the GCM graphs in \DynkinDiagramFigure, all numbers games are convergent. 
Applying results of Eriksson, we will then see that 
two numbers games played from the same initial position on one of these graphs 
converge to the same terminal position in the same number of steps. 
\NotMarsFriendlyCatalog\ asserts that for the GCM graphs of 
\NotMarsFriendlyFigure, a 
numbers game is divergent from any nonzero initial position with nonnegative 
numbers. 
Two key results needed for proofs of both propositions 
are Eriksson's Strong Convergence 
and Comparison Theorems (see \ErikssonTheorems\ below).  
Proofs of both propositions also involve 
case analysis arguments that are fairly routine, can be checked by 
hand, and are often easily expedited using a computer algebra system to 
automate some of the computations. 
Complete details are provided here.  

Fix a positive integer $n$ and a totally ordered set $I_{n}$ with $n$ 
elements (usually $I_{n} := \{1<\ldots<n\}$).  
A {\em generalized Cartan matrix} (or {\em GCM}) is an $n \times n$ 
matrix $M = (M_{ij})_{i,j \in I_{n}}$  
with integer entries satisfying the requirements that each 
main diagonal matrix entry is 2, that all other matrix entries are 
nonpositive, and that if a matrix entry $M_{ij}$ is nonzero then its 
transpose entry $M_{ji}$ is also nonzero.  
Generalized Cartan matrices 
are the starting point for the study of 
Kac--Moody 
algebras: beginning with a GCM, one can 
write down a list of the defining relations for a Kac--Moody 
algebra as well as the associated Weyl group 
(see \cite{Kac} or \cite{Kumar}).  
To an $n \times n$ generalized Cartan matrix 
$M = (M_{ij})_{i,j \in I_{n}}$ we associate a finite 
graph $\Gamma$ (which has undirected edges, 
no loops, and no multiple edges) 
as follows:     
The nodes $(\gamma_{i})_{i \in I_{n}}$ of $\Gamma$ are indexed 
by the set $I_{n}$, 
and   an edge is placed between nodes $\gamma_{i}$ and $\gamma_{j}$ 
if and only if $i \not= j$ 
and the matrix entries $M_{ij}$ and $M_{ji}$ are nonzero.  We call the pair 
$(\Gamma,M)$ a {\em GCM graph}. 
We consider two 
GCM graphs $(\Gamma, M = (M_{ij})_{i,j \in I_{n}})$ 
and $(\Gamma', M' = (M'_{pq})_{p,q \in I'_{n}})$ 
to be the same if under some bijection $\sigma: I_{n} \rightarrow 
I'_{n}$ we have nodes $\gamma_{i}$ and $\gamma_{j}$ in 
$\Gamma$ adjacent if and only if 
$\gamma'_{\sigma(i)}$ and $\gamma'_{\sigma(j)}$ are adjacent in $\Gamma'$ with 
$M_{ij} = M'_{\sigma(i),\sigma(j)}$. 
With $p = -M_{12}$ and $q = -M_{21}$, 
we depict a generic connected two-node GCM graph as follows:

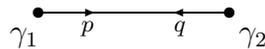
\vspace*{-0.075in}
\noindent
\begin{center}
\TwoCitiesGraphWithLabels
\end{center}

\vspace*{-0.075in}
\noindent 
We use special names and notation to refer to 
two-node GCM 
graphs which have $p = 1$ and $q = 1$, $2$, or $3$ respectively:
\noindent
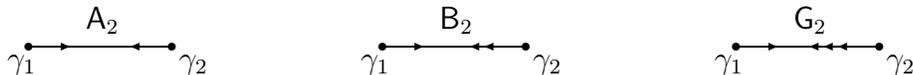
\begin{center}
\ATwoGraphWithLabels
\hspace*{1in}
\BTwoGraphWithLabels
\hspace*{1in}
\GTwoGraphWithLabels
\end{center}
When $p=1$ and $q=1$ it is convenient to use the graph 
\ATwoGraphNoEdgeLabels\ to 
represent the GCM graph $\myA_{2}$.  
A GCM graph $(\Gamma,M)$ is a {\em Dynkin diagram of finite type} if each connected 
component of $(\Gamma,M)$ 
is one of the graphs of \DynkinDiagramFigure.  We number our nodes as 
in \S 11.4 of \cite{Hum}. In these cases the GCMs  are ``Cartan'' 
matrices. 

\begin{figure}[t]
\begin{center}
\DynkinDiagramFigure: Connected Dynkin diagrams of finite type. 
\end{center}

\vspace*{-0.55in}
\begin{tabular}{cl}
$\myA_{n}$ ($n \geq 1$) & \AnMarsOK\\

$\myB_{n}$ ($n \geq 2$) & \BnMarsOK\\

$\myC_{n}$ ($n \geq 3$) & \CnMarsOK\\

$\myD_{n}$ ($n \geq 4$) & \DnMarsOK\\

$\myE_{6}$ & \ESixMarsOK\\

$\myE_{7}$ & \ESevenMarsOK\\

$\myE_{8}$ & \EEightMarsOK\\

$\myF_{4}$ & \FFourMarsOK\\

$\myG_{2}$ & \GTwoMarsOK
\end{tabular}

\vspace*{-0.25in}
\end{figure}

A {\em position} $\lambda = (\lambda_i)_{i \in 
I_{n}}$ is an assignment of real numbers to the nodes of the GCM graph 
$(\Gamma,M)$. 
The position $\lambda$ is 
{\em dominant} (respectively, {\em strongly dominant}) if 
$\lambda_{i} \geq 0$ 
(resp. $\lambda_i > 0$) for all $i \in I_{n}$; 
$\lambda$ is {\em nonzero} if at least one $\lambda_i \not= 0$. 
For $i \in I_{n}$, the {\em 
fundamental position} $\omega_i$ is the assignment of the number  
$1$ at node $\gamma_{i}$ and the number $0$ at all other nodes.  
Given a position $\lambda$ on a GCM graph $(\Gamma,M)$, to 
{\em fire} a node $\gamma_{i}$ is to change the number at each node 
$\gamma_{j}$ of $\Gamma$ by the transformation  
\[\lambda_j	 \longmapsto \lambda_j - 
M_{ij}\lambda_i,\] provided the number at node 
$\gamma_{i}$ is positive; otherwise node $\gamma_{i}$ is not allowed 
to be fired. 
Since the generalized Cartan  
matrix $M$ assigns a pair of {\em amplitudes} ($M_{ij}$ and 
$M_{ji}$) to each edge of the 
graph $\Gamma$, we sometimes refer to GCMs as 
{\em amplitude matrices}.  
The {\em numbers game} 
is the one-player game on a GCM graph $(\Gamma,M)$ in which the player 
(1) Assigns an initial position  
to the nodes of $\Gamma$; (2) Chooses a node with a positive 
number and fires the node to obtain a new position; and (3) 
Repeats step (2) for the new position if there is at least one node 
with a positive number.  
\begin{figure}[thb]
\begin{center}
\BtwoMarsDemo: The numbers game for the GCM graph $\myB_{2}$. 

\vspace*{0.25in}
\setlength{\unitlength}{0.9in}
\begin{picture}(1.5,5.15) 
\put(0.25,5.1){\BTwoGraphForFigure}
\put(0.25,5.3){$a$}
\put(1.4,5.3){$b$}
\put(0.3,4.9){\vector(-4,-3){0.7}}
\put(1.4,4.9){\vector(4,-3){0.7}}
\put(1.75,3.8){\BTwoGraphForFigure}
\put(1.5,4.0){$a+2b$}
\put(2.8,4.0){$-b$}
\put(2.4,3.6){\vector(0,-1){0.6}}
\put(1.75,2.5){\BTwoGraphForFigure}
\put(1.45,2.7){$-a-2b$}
\put(2.75,2.7){$a+b$}
\put(2.4,2.3){\vector(0,-1){0.6}}
\put(1.75,1.2){\BTwoGraphForFigure}
\put(1.75,1.4){$a$}
\put(2.65,1.4){$-a-b$}
\put(2.1,1.0){\vector(-4,-3){0.7}}
\put(-1.25,3.8){\BTwoGraphForFigure}
\put(-1.35,4.0){$-a$}
\put(-0.25,4.0){$a+b$}
\put(-0.7,3.6){\vector(0,-1){0.6}}
\put(-1.25,2.5){\BTwoGraphForFigure}
\put(-1.5,2.7){$a+2b$}
\put(-0.35,2.7){$-a-b$}
\put(-0.7,2.3){\vector(0,-1){0.6}}
\put(-1.25,1.2){\BTwoGraphForFigure}
\put(-1.55,1.4){$-a-2b$}
\put(-0.1,1.4){$b$}
\put(-0.4,1.0){\vector(4,-3){0.7}}
\put(0.25,-0.1){\BTwoGraphForFigure}
\put(0.15,0.1){$-a$}
\put(1.3,0.1){$-b$}
\end{picture}
\end{center}
\end{figure}
Consider now the GCM graph $\myB_{2}$.  As we can see in \BtwoMarsDemo, 
the numbers game terminates in a finite number of steps for any 
initial position and any legal sequence of node firings, 
if it is understood that the player  
will continue to fire as long as there is at least one 
node with a positive number.  In general, 
given a position $\lambda$, a {\em game sequence 
for} $\lambda$ is the (possibly empty, possibly infinite) sequence 
$(\gamma_{i_{1}}, \gamma_{i_{2}},\ldots)$, where 
$\gamma_{i_{j}}$ 
is the $j$th node that is fired in some 
numbers game with initial position $\lambda$.  
More generally, a {\em firing sequence} from some position $\lambda$ is an 
initial portion of some game sequence played from $\lambda$; the 
phrase {\em legal firing sequence} is used to emphasize that all node 
firings in the sequence are known or assumed to be possible. 
Note that a game sequence 
$(\gamma_{i_{1}}, 
\gamma_{i_{2}},\ldots,\gamma_{i_{l}})$ 
is of finite length $l$ 
(possibly with $l = 0$) if 
the number is nonpositive at each node after the $l$th firing; in 
this case we say the game sequence is {\em convergent} and the 
resulting position is the {\em terminal position} for the game 
sequence.  
We say a connected GCM graph $(\Gamma,M)$ is {\em admissible} if 
there exists a nonzero dominant initial position with a convergent 
game sequence. 

The following preliminary results are 
needed for the proofs of \MainResults. 
These results also appear in \cite{DE} and \cite{DonEnumbers} for use 
in proofs of key theorems of those papers.  Proofs or references for 
these results are also given here.  
Following \cite{ErikssonThesis} and \cite{ErikssonEur}, we say 
the numbers game on a GCM graph $(\Gamma,M)$ is {\em strongly 
convergent} if given any initial position, any two game sequences 
either both diverge or both converge to the same terminal position in the 
same number of steps. The next result  
follows from  
Theorem 3.1 of \cite{ErikssonEur} 
(or see Theorem 3.6 of \cite{ErikssonThesis}).  

\noindent
{\bf \StrongConvergenceTheorem\ (Eriksson's Strong Convergence Theorem)}\ \ 
{\sl The numbers game on a connected GCM graph 
is strongly 
convergent.}  

The 
following weaker result also applies when the 
GCM graph is not connected:  

\noindent
{\bf \StrongConvergenceCorollary} \ \ {\sl For any GCM graph, 
if a game sequence for an initial position  
$\lambda$ diverges, then all game  
sequences for $\lambda$ diverge.}
 
The next result is an immediate consequence of Theorem 4.3 of 
\cite{ErikssonThesis} or Theorem 4.5 of 
\cite{ErikssonDiscrete}.  
Eriksson's proof of this result in \cite{ErikssonThesis} uses only 
combinatorial and linear algebraic methods. 
 
\noindent 
{\bf \ComparisonTheorem\ (Eriksson's Comparison Theorem)}\ \ 
{\sl Given a GCM graph, suppose that a game sequence 
for an initial position $\lambda = (\lambda_{i})_{i \in 
I_{n}}$ converges.  Suppose that a position $\lambda' := (\lambda'_{i})_{i 
\in I_{n}}$ has the property that $\lambda'_{i} \leq 
\lambda_{i}$ for all $i \in I_{n}$.  Then some game sequence 
for the initial position $\lambda'$ also converges.}

Let $r$ be a positive real number.  
Observe that if 
$(\gamma_{i_{1}},\ldots,\gamma_{i_{l}})$ 
is a convergent 
game sequence for an initial position $\lambda = 
(\lambda_{i})_{i \in I_{n}}$, then 
$(\gamma_{i_{1}},\ldots,\gamma_{i_{l}})$ 
is a convergent 
game sequence for the initial position $r\lambda := 
(r\lambda_{i})_{i \in I_{n}}$. 
This observation and \ComparisonTheorem\ imply the following 
result: 

\noindent 
{\bf \ComparisonCorollary}\ \ {\sl 
Let $\lambda = (\lambda_{i})_{i \in 
I_{n}}$ be a dominant initial position such that $\lambda_{j} > 0$ 
for some $j \in I_{n}$. Suppose that a game sequence for $\lambda$ 
converges.  
Then some game   
sequence for the fundamental position $\omega_{j}$ also 
converges.}

The following is an immediate consequence of \ComparisonResults: 

\noindent 
{\bf \NotMarsFriendlyLemma}\ \ {\sl A GCM graph is not admissible if 
for each fundamental position there is a divergent game  
sequence.}

\vspace{1ex} 

\noindent
{\Large \bf \FirstProofNum.\ \ Convergent numbers games on Dynkin diagrams 
of finite type}

\vspace{1ex} 
Eriksson's Strong Convergence and Comparison Theorems are key steps in 
our proof of \ConvergenceProposition.  The remaining step, which 
accounts for most of the length of this section, is to 
provide convergent game sequences for numbers games played from 
strongly dominant positions on 
connected Dynkin diagrams of finite type. 
Finding convergent game sequences for numbers games played on Dynkin 
diagrams of finite type may seem like a difficult task at first, but 
in view of \ConvergenceProposition, there is no way to go wrong: any 
two numbers games played from the same initial position 
will terminate in the same finite number of steps.  

A general theory connecting the numbers game and 
Coxeter/Weyl group actions was 
developed by Eriksson in \cite{ErikssonThesis} and 
\cite{ErikssonDiscrete}.  
From this theory 
it follows that for a numbers game played from a strongly dominant 
position on a Dynkin diagram of finite type, 
any game sequence corresponds to a reduced expression for the longest 
element of the corresponding Weyl group 
and conversely any reduced expression for the longest Weyl group element 
corresponds to a game sequence 
(see Propositions 4.1 and 4.2 of \cite{ErikssonDiscrete} or 
Theorem 4.3.1 part ({\em iv}) of \cite{BB}).  
Moreover, the length of the 
game sequence 
is the length of any such reduced expression and is also equal to the 
number of positive roots in the associated root system.  For 
further discussion of this phenomenon, see \cite{DonEnumbers}. 

For the four infinite families of connected 
Dynkin diagrams of finite type, 
the next results are proved by induction on $n$, the number of 
nodes.  This is effected by observing natural ``GCM subgraph'' 
inclusions 
$\myA_{n-1} \hookrightarrow \myA_{n}$ ($n \geq 2$), 
$\myB_{n-1} \hookrightarrow \myB_{n}$ ($n \geq 3$), 
$\myC_{n-1} \hookrightarrow \myC_{n}$ ($n \geq 4$), and 
$\myD_{n-1} \hookrightarrow \myD_{n}$ ($n \geq 5$). 
If $I'_{m}$ is a subset of the node set $I_{n}$ of a GCM graph 
$(\Gamma,M)$, 
then let $\Gamma'$ be the subgraph of $\Gamma$ with node set $I'_{m}$ and 
the induced set of edges, and let $M'$ be the corresponding 
submatrix 
of the amplitude matrix $M$; we call $(\Gamma',M')$ a {\em GCM 
subgraph} of $(\Gamma,M)$.

\noindent 
{\bf \FourFamiliesFirstLemma}\\
\noindent {\bf A.}\ \ {\sl For  $n \geq 2$  
and for any strongly dominant position $(a_{1},\ldots,a_{n})$ 
on} $\myA_{n}${\sl , one can obtain the 
position $(a_{1}+\cdots+a_{n},-a_{n},\ldots,-a_{3},-a_{2})$ 
by the sequence 
$(\mathbf{s}_{n}$, $\mathbf{s}_{n-1}$, $\ldots$, 
$\mathbf{s}_{2})$ of 
legal node firings where  
$\mathbf{s}_{i}$ is the subsequence 
$(\gamma_{i}$, $\gamma_{i+1}$, $\ldots$, $\gamma_{n-1}$, $\gamma_{n})$ 
for $2 \leq i \leq n$.}

\noindent {\bf B.}\ \ {\sl 
A similar statement holds for} $\myB_{n}$ {\sl with $n \geq 2$: The initial 
strongly dominant position is $(a_{1},\ldots,a_{n})$, and the position 
$(a_{1}+2a_{2}+\cdots+2a_{n-1}+a_{n},-a_{2},-a_{3},\ldots,-a_{n})$ is 
obtained by the sequence 
$(\mathbf{s}_{n}$, $\mathbf{s}_{n-1}$, $\ldots$, 
$\mathbf{s}_{2})$ of 
legal node firings where  
$\mathbf{s}_{i}$ is the subsequence 
$(\gamma_{i}$, $\gamma_{i+1}$, $\ldots$, $\gamma_{n-1}$, $\gamma_{n}$, 
$\gamma_{n-1}$, $\ldots$, $\gamma_{i+1}$, 
$\gamma_{i})$ for $2 \leq i \leq n-1$ and 
$\mathbf{s}_{n} = (\gamma_{n})$.}

\noindent {\bf C.}\ \ {\sl  
A similar statement holds for} $\myC_{n}$ {\sl with $n \geq 3$: The initial 
strongly dominant position is $(a_{1},\ldots,a_{n})$, and the position 
$(a_{1}+2a_{2}+\cdots+2a_{n-1}+2a_{n},-a_{2},-a_{3},\ldots,-a_{n})$ is 
obtained by the sequence 
$(\mathbf{s}_{n}$, $\mathbf{s}_{n-1}$, $\ldots$, 
$\mathbf{s}_{2})$ of 
legal node firings where  
$\mathbf{s}_{i}$ is the subsequence 
$(\gamma_{i}$, $\gamma_{i+1}$, $\ldots$, $\gamma_{n-1}$, $\gamma_{n}$, 
$\gamma_{n-1}$, $\ldots$, $\gamma_{i+1}$, 
$\gamma_{i})$ for $2 \leq i \leq n-1$ 
and $\mathbf{s}_{n} = (\gamma_{n})$.}

\noindent {\bf D.}\ \ {\sl 
A similar statement holds for} $\myD_{n}$ {\sl with $n \geq 4$: The initial 
strongly dominant position is $(a_{1},\ldots,a_{n})$.  
Let $b_{n-1} := a_{n-1}$ and $b_{n} := 
a_{n}$ when $n$ is odd 
and where $b_{n-1} := a_{n}$ and $b_{n} := a_{n-1}$ 
when $n$ is even.  
The position 
$(a_{1}+2a_{2}+\cdots+2a_{n-2}+a_{n-1}+a_{n},-a_{2},-a_{3},\ldots,
-a_{n-2},-b_{n-1},-b_{n})$ 
is obtained by the sequence 
$(\mathbf{s}_{n-1}$, $\mathbf{s}_{n-2}$, $\ldots$, 
$\mathbf{s}_{2})$ of 
legal node firings where  
$\mathbf{s}_{i}$ is the subsequence 
$(\gamma_{i}$, $\gamma_{i+1}$, $\ldots$, $\gamma_{n-2}$, $\gamma_{n-1}$, $\gamma_{n}$, 
$\gamma_{n-2}$, $\ldots$, $\gamma_{i+1}$, $\gamma_{i})$ for $2 \leq i \leq n-2$ and 
$\mathbf{s}_{n-1} = (\gamma_{n-1}, \gamma_{n})$.}

{\em Proof.} In case $\myA$, the result clearly holds for the two-node 
graph.  As our 
induction hypothesis, assume the lemma statement holds for all type 
$\myA$ Dynkin diagrams with fewer than $n$ nodes.  Given $\myA_{n}$, the 
GCM subgraph determined by the $n-1$ rightmost nodes is an $\myA_{n-1}$ 
Dynkin diagram.  Applying the induction hypothesis, the legal firing 
sequence 
$(\mathbf{s}_{n}$, $\ldots$, $\mathbf{s}_{3})$ from the strongly 
dominant position $(a_{1}, a_{2}, \ldots, a_{n})$ results in the 
position $(a_{1}, a_{2}+\cdots+a_{n}, -a_{n}, \ldots, -a_{3})$. To 
this position we now apply the sequence $\mathbf{s}_{2} = 
(\gamma_{2}$, $\ldots$, $\gamma_{n-1}$, $\gamma_{n})$.  
Once $\gamma_{2}$ is fired, then 
for $3 \leq i \leq n$ it is easily seen that 
just before $\gamma_{i}$ is fired 
in the sequence $\mathbf{s}_{2}$ 
the position is $(a_{1}+a_{2}+\cdots+a_{n}$, $-a_{n}$, $-a_{n-1}$, 
$\ldots$, $-a_{n+5-i}$, $-a_{n+4-i}$, $-a_{2}-a_{3}-\cdots-a_{n+3-i}$, 
$a_{2}+\cdots+a_{n+2-i}$, $-a_{n+2-i}$, $-a_{n+1-i}$, $\ldots$, 
$-a_{4}$, $-a_{3})$.   
Then each firing in the sequence $\mathbf{s}_{2}$ 
is legal, and the resulting position is 
 $(a_{1}+a_{2}+\cdots+a_{n}$, $-a_{n}$, $-a_{n-1}$, $\ldots$, 
$-a_{3}$, $-a_{2})$.

In case $\myB$, the result clearly holds for the two-node graph. 
As our induction hypothesis, assume the lemma statement holds for all type 
$\myB$ Dynkin diagrams with fewer than $n$ nodes.  Given $\myB_{n}$, the 
GCM subgraph determined by the $n-1$ rightmost nodes is a $\myB_{n-1}$ 
Dynkin diagram.  Applying the induction hypothesis, the legal firing 
sequence 
$(\mathbf{s}_{n}$, $\ldots$, $\mathbf{s}_{3})$ from the strongly 
dominant position $(a_{1}, a_{2}, \ldots, a_{n})$ results in the 
position $(a_{1}, a_{2}+2a_{3}+\cdots+2a_{n-1}+a_{n}, -a_{3}, \ldots, 
-a_{n})$. To 
this position we now apply the sequence $\mathbf{s}_{2} = 
(\gamma_{2}$, $\ldots$, $\gamma_{n-1}$, $\gamma_{n}$, 
$\gamma_{n-1}$, $\ldots$ $\gamma_{2})$.  
Once $\gamma_{2}$ is fired, then 
for $3 \leq i \leq n-1$ it is easily seen that 
just before $\gamma_{i}$ is fired 
for the first time 
in the sequence $\mathbf{s}_{2}$ 
the position is 
$(a_{1}+a_{2}+2a_{3}+\cdots+2a_{n-1}+a_{n}$, 
$-a_{3}$, $-a_{4}$,  
$\ldots$, $-a_{i-1}$, 
$-a_{2}-a_{3}-\cdots-a_{i-1}-2a_{i}-\cdots-2a_{n-1}-a_{n}$, 
$a_{2}+a_{3}+\cdots+a_{i-1}+a_{i}+2a_{i+1}+\cdots+2a_{n-1}+a_{n}$, 
$-a_{i+1}$, $-a_{i+2}$, $\ldots$, 
$-a_{n})$. 
One now sees that 
just before $\gamma_{n}$ is fired    
the position is 
$(a_{1}+a_{2}+2a_{3}+\cdots+2a_{n-1}+a_{n}$, 
$-a_{3}$, $-a_{4}$,  
$\ldots$, $-a_{n-1}$, 
$-a_{2}-a_{3}-\cdots-a_{n-1}-a_{n}$, 
$2a_{2}+2a_{3}+\cdots+2a_{n-1}+a_{n})$. 
Now for $3 \leq i \leq n-1$, it is easily seen that 
just before $\gamma_{i}$ is fired 
for the second time 
in the sequence $\mathbf{s}_{2}$    
the position is 
$(a_{1}+a_{2}+2a_{3}+\cdots+2a_{n-1}+a_{n}$, 
$-a_{3}$, $-a_{4}$,  
$\ldots$, $-a_{i}$, 
$a_{2}+a_{3}+\cdots+a_{i-1}+a_{i}$, 
$-a_{2}-a_{3}+\cdots-a_{i}-a_{i+1}$, 
$-a_{i+2}$, $-a_{i+3}$, $\ldots$, 
$-a_{n})$. 
Finally, fire $\gamma_{2}$ from the position 
$(a_{1}+a_{2}+2a_{3}+\cdots+2a_{n-1}+a_{n}$, 
$a_{2}$, $-a_{2}-a_{3}$, $-a_{4}$,  
$\ldots$, $-a_{n})$. 
Then each firing in the sequence $\mathbf{s}_{2}$ 
is legal, and the resulting position is 
$(a_{1}+2a_{2}+\cdots+2a_{n-1}+a_{n}$, $-a_{2}$, $-a_{3}$, $\ldots$, 
$-a_{n})$. 

In case $\myC$, it is easy to confirm that the specified sequence of 
four legal node firings 
played from a strongly dominant position on the three-node 
graph yields the stated resulting position.  As our 
induction hypothesis, assume the lemma statement holds for all type 
$\myC$ Dynkin diagrams with fewer than $n$ nodes.  Given $\myC_{n}$, the 
GCM subgraph determined by the $n-1$ rightmost nodes is a $\myC_{n-1}$ 
Dynkin diagram.  Applying the induction hypothesis, the legal firing 
sequence 
$(\mathbf{s}_{n}$, $\ldots$, $\mathbf{s}_{3})$ from the strongly 
dominant position $(a_{1}, a_{2}, \ldots, a_{n})$ results in the 
position $(a_{1}, a_{2}+2a_{3}+\cdots+2a_{n}, -a_{3}, \ldots, 
-a_{n})$. To 
this position we now apply the sequence $\mathbf{s}_{2} = 
(\gamma_{2}$, $\ldots$, $\gamma_{n-1}$, $\gamma_{n}$, 
$\gamma_{n-1}$, $\ldots$ $\gamma_{2})$.  
Once $\gamma_{2}$ is fired, then 
for $3 \leq i \leq n-1$ it is easily seen that 
just before $\gamma_{i}$ is fired 
for the first time 
in the sequence $\mathbf{s}_{2}$    
the position is 
$(a_{1}+a_{2}+2a_{3}+\cdots+2a_{n}$, 
$-a_{3}$, $-a_{4}$,  
$\ldots$, $-a_{i-1}$, 
$-a_{2}-a_{3}-\cdots-a_{i-1}-2a_{i}-\cdots-2a_{n}$, 
$a_{2}+a_{3}+\cdots+a_{i-1}+a_{i}+2a_{i+1}+\cdots+2a_{n}$, 
$-a_{i+1}$, $-a_{i+2}$, $\ldots$, 
$-a_{n})$.  
One now sees that 
just before $\gamma_{n}$ is fired    
the position is 
$(a_{1}+a_{2}+2a_{3}+\cdots+2a_{n}$, 
$-a_{3}$, $-a_{4}$,  
$\ldots$, $-a_{n-1}$, 
$-a_{2}-a_{3}-\cdots-a_{n-1}-2a_{n}$, 
$a_{2}+a_{3}+\cdots+a_{n})$. 
Now for $3 \leq i \leq n-1$, it is easily seen that 
just before $\gamma_{i}$ is fired 
for the second time 
in the sequence $\mathbf{s}_{2}$     
the position is 
$(a_{1}+a_{2}+2a_{3}+\cdots+2a_{n}$, 
$-a_{3}$, $-a_{4}$,  
$\ldots$, $-a_{i}$, 
$a_{2}+a_{3}+\cdots+a_{i-1}+a_{i}$, 
$-a_{2}-a_{3}+\cdots-a_{i}-a_{i+1}$, 
$-a_{i+2}$, $-a_{i+3}$, $\ldots$, 
$-a_{n})$.  
Finally, fire $\gamma_{2}$ from the position 
$(a_{1}+a_{2}+2a_{3}+\cdots+2a_{n}$, 
$a_{2}$, $-a_{2}-a_{3}$, $-a_{4}$,  
$\ldots$, $-a_{n})$. 
Then each firing in the sequence $\mathbf{s}_{2}$ 
is legal, and the resulting position is 
$(a_{1}+2a_{2}+\cdots+2a_{n}$, $-a_{2}$, $-a_{3}$, $\ldots$, 
$-a_{n})$. 

In case $\myD$, it is easy to confirm that the specified sequence of 
six legal node firings 
played from a strongly dominant position on the four-node 
graph yields the stated resulting position.  
As our induction hypothesis, assume the lemma statement holds 
for all type 
$\myD$ Dynkin diagrams with fewer than $n$ nodes.  
Given $\myD_{n}$, the 
GCM subgraph determined by the $n-1$ rightmost nodes is a $\myD_{n-1}$ 
Dynkin diagram.  
Assume for the moment that $n$ is even, so $n-1$ is odd. 
Applying the induction hypothesis, the legal firing 
sequence 
$(\mathbf{s}_{n-1}$, $\ldots$, $\mathbf{s}_{3})$ from the strongly 
dominant position $(a_{1}, a_{2}, \ldots, a_{n})$ results in the 
position $(a_{1}, a_{2}+2a_{3}+\cdots+2a_{n-2}+a_{n-1}+a_{n}$, 
$-a_{3}$, $\ldots$, 
$-a_{n-2}$, $-a_{n-1}$, $-a_{n})$. 
To this position we now apply the sequence $\mathbf{s}_{2} = 
(\gamma_{2}$, $\ldots$, $\gamma_{n-2}$, $\gamma_{n-1}$, 
$\gamma_{n}$, $\gamma_{n-2}$, $\ldots$ $\gamma_{2})$.  
Once $\gamma_{2}$ is fired, then 
for $3 \leq i \leq n-3$ 
it is easily seen that 
just before $\gamma_{i}$ is fired 
for the first time 
in the sequence $\mathbf{s}_{2}$     
the position is 
$(a_{1}+a_{2}+2a_{3}+\cdots+2a_{n-2}+a_{n-1}+a_{n}$, 
$-a_{3}$, $-a_{4}$,  
$\ldots$, $-a_{i-1}$, 
$-a_{2}-a_{3}-\cdots-a_{i-1}-2a_{i}-\cdots-2a_{n-2}-a_{n-1}-a_{n}$, 
$a_{2}+a_{3}+\cdots+a_{i-1}+a_{i}+2a_{i+1}+\cdots+2a_{n-2}+a_{n-1}+a_{n}$, 
$-a_{i+1}$, $-a_{i+2}$, $\ldots$, 
$-a_{n-2}$, $-a_{n-1}$, $-a_{n})$.  
One now sees that 
just before $\gamma_{n-2}$ is fired for the first time in the sequence  
the position is 
$(a_{1}+a_{2}+2a_{3}+\cdots+2a_{n-2}+a_{n-1}+a_{n}$, 
$-a_{3}$, $-a_{4}$,  
$\ldots$, $-a_{n-3}$, 
$-a_{2}-a_{3}-\cdots-a_{n-3}-2a_{n-2}-a_{n-1}-a_{n}$, 
$a_{2}+a_{3}+\cdots+a_{n-2}+a_{n-1}+a_{n}$, 
$-a_{n-1}$, $-a_{n})$.  
Then just before $\gamma_{n-1}$ is fired the position is 
$(a_{1}+a_{2}+2a_{3}+\cdots+2a_{n-2}+a_{n-1}+a_{n}$, 
$-a_{3}$, $-a_{4}$,  
$\ldots$, $-a_{n-3}$, 
$-a_{n-2}$, 
$-a_{2}-a_{3}-\cdots-a_{n-2}-a_{n-1}-a_{n}$, 
$a_{2}+a_{3}+\cdots+a_{n-2}+a_{n}$, 
$a_{2}+a_{3}+\cdots+a_{n-2}+a_{n-1})$,    
and just before $\gamma_{n}$ is fired the position is 
$(a_{1}+a_{2}+2a_{3}+\cdots+2a_{n-2}+a_{n-1}+a_{n}$, 
$-a_{3}$, $-a_{4}$,  
$\ldots$, $-a_{n-3}$, 
$-a_{n-2}$, 
$-a_{n-1}$, 
$-a_{2}-a_{3}-\cdots-a_{n-2}-a_{n}$, 
$a_{2}+a_{3}+\cdots+a_{n-2}+a_{n-1})$.   
So just before $\gamma_{n-2}$ is fired for the second time in the 
sequence $\mathbf{s}_{2}$   
the position is 
$(a_{1}+a_{2}+2a_{3}+\cdots+2a_{n-2}+a_{n-1}+a_{n}$, 
$-a_{3}$, $-a_{4}$,  
$\ldots$, $-a_{n-3}$, 
$-a_{n-2}$, 
$a_{2}+a_{3}+\cdots+a_{n-2}$, 
$-a_{2}-a_{3}-\cdots-a_{n-2}-a_{n}$, 
$-a_{2}-a_{3}-\cdots-a_{n-2}-a_{n-1})$.  
Now for $3 \leq i \leq n$, it is easily seen that 
just before $\gamma_{i}$ is fired 
for the second time 
in the sequence $\mathbf{s}_{2}$     
the position is 
$(a_{1}+a_{2}+2a_{3}+\cdots+2a_{n}$, 
$-a_{3}$, $-a_{4}$,  
$\ldots$, $-a_{i}$, 
$a_{2}+a_{3}+\cdots+a_{i-1}+a_{i}$, 
$-a_{2}-a_{3}+\cdots-a_{i}-a_{i+1}$, 
$-a_{i+2}$, $-a_{i+3}$, $\ldots$, $-a_{n-2}$, 
$-a_{n}$, $-a_{n-1})$.   
Finally, fire $\gamma_{2}$ from the position 
$(a_{1}+a_{2}+2a_{3}+\cdots+2a_{n-2}+a_{n-1}+a_{n}$, 
$a_{2}$, $-a_{2}-a_{3}$, $-a_{4}$,  
$\ldots$, $-a_{n-2}$, $-a_{n}$, $-a_{n-1})$. 
Then each firing in the sequence $\mathbf{s}_{2}$ 
is legal, and the resulting position is 
$(a_{1}+2a_{2}+\cdots+2a_{n-2}+a_{n-1}+a_{n}$, $-a_{2}$, $-a_{3}$, $\ldots$, 
$a_{n-2}$, $-a_{n}$, $-a_{n-1})$.  When $n$ is odd, the argument is 
entirely similar.\hfill\QED 

From this we immediately obtain the following: 

\noindent 
{\bf \FourFamiliesSecondLemma}\\
\noindent {\bf A.}\ \ {\sl For any positive integer $n$  
and for any strongly dominant position $(a_{1},\ldots,a_{n})$ 
on} $\myA_{n}${\sl , one can obtain the 
position $(-a_{n},\ldots,-a_{2},-a_{1})$ 
by the sequence 
$(\mathbf{s}_{n}$, $\mathbf{s}_{n-1}$, $\ldots$, 
$\mathbf{s}_{1})$ of 
legal node firings where  
$\mathbf{s}_{i}$ is the subsequence 
$(\gamma_{i}$, $\gamma_{i+1}$, $\ldots$, $\gamma_{n-1}$, $\gamma_{n})$ 
for $1 \leq i \leq n$.}

\noindent {\bf B.}\ \ {\sl 
A similar statement holds for} $\myB_{n}$ {\sl with $n \geq 2$: The initial 
strongly dominant position is $(a_{1},\ldots,a_{n})$, and the position 
$(-a_{1},-a_{2},\ldots,-a_{n})$ is 
obtained by the sequence 
$(\mathbf{s}_{n}$, $\mathbf{s}_{n-1}$, $\ldots$, 
$\mathbf{s}_{1})$ of 
legal node firings where  
$\mathbf{s}_{i}$ is the subsequence 
$(\gamma_{i}$, $\gamma_{i+1}$, $\ldots$, $\gamma_{n-1}$, $\gamma_{n}$, 
$\gamma_{n-1}$, $\ldots$, $\gamma_{i+1}$, 
$\gamma_{i})$ for $1 \leq i \leq n-1$ and 
$\mathbf{s}_{n} = (\gamma_{n})$.}

\noindent {\bf C.}\ \ {\sl  
A similar statement holds for} $\myC_{n}$ {\sl with $n \geq 3$: The initial 
strongly dominant position is $(a_{1},\ldots,a_{n})$, and the position 
$(-a_{1},-a_{2},\ldots,-a_{n})$ is 
obtained by the sequence 
$(\mathbf{s}_{n}$, $\mathbf{s}_{n-1}$, $\ldots$, 
$\mathbf{s}_{1})$ of 
legal node firings where  
$\mathbf{s}_{i}$ is the subsequence 
$(\gamma_{i}$, $\gamma_{i+1}$, $\ldots$, $\gamma_{n-1}$, $\gamma_{n}$, 
$\gamma_{n-1}$, $\ldots$, $\gamma_{i+1}$, 
$\gamma_{i})$ for $1 \leq i \leq n-1$ 
and $\mathbf{s}_{n} = (\gamma_{n})$.}

\noindent {\bf D.}\ \ {\sl 
A similar statement holds for} $\myD_{n}$ {\sl with $n \geq 4$: The initial 
strongly dominant position is $(a_{1},\ldots,a_{n})$.  
Let $b_{n-1} := a_{n-1}$ and $b_{n} := 
a_{n}$ when $n$ is even 
and where $b_{n-1} := a_{n}$ and $b_{n} := a_{n-1}$ 
when $n$ is odd.  
The position 
$(a_{1}+2a_{2}+\cdots+2a_{n-2}+a_{n-1}+a_{n},-a_{2},-a_{3},\ldots,
-a_{n-2},-b_{n-1},-b_{n})$ 
is obtained by the sequence 
$(\mathbf{s}_{n-1}$, $\mathbf{s}_{n-2}$, $\ldots$, 
$\mathbf{s}_{1})$ of 
legal node firings where  
$\mathbf{s}_{i}$ is the subsequence 
$(\gamma_{i}$, $\gamma_{i+1}$, $\ldots$, $\gamma_{n-2}$, $\gamma_{n-1}$, $\gamma_{n}$, 
$\gamma_{n-2}$, $\ldots$, $\gamma_{i+1}$, $\gamma_{i})$ for $1 \leq i \leq n-2$ and 
$\mathbf{s}_{n-1} = (\gamma_{n-1}, \gamma_{n})$.}

{\em Proof.}  View $\myA_{n}$ (respectively $\myB_{n}$, $\myC_{n}$, 
$\myD_{n}$) 
as a GCM subgraph of $\myA_{n+1}$ (respectively $\myB_{n+1}$, $\myC_{n+1}$, 
$\myD_{n+1}$) by adding a 
node ``to the left'' of $\gamma_{1}$.  Conclude by applying 
\FourFamiliesFirstLemma.\hfill\QED

\noindent 
{\bf \ConvergenceProposition}\ \ {\sl A connected Dynkin diagram 
$(\Gamma,M)$ of 
finite type is admissible.  Moreover, for any initial position on 
$(\Gamma,M)$, all game sequences converge to the same terminal 
position in the same finite number of steps.} 

{\em Proof.} 
The Strong Convergence Theorem  
shows that if a game sequence for some initial position $\lambda$ on 
$(\Gamma,M)$ converges, 
then all game sequences from $\lambda$ converge to the same 
terminal position in the 
same finite number of steps.  
Then in light of 
The Comparison Theorem, 
it  suffices to show that for 
any strongly dominant initial position on  
$(\Gamma,M)$, 
there is a convergent game sequence.  

For the four infinite families ($\myA_{n}$, $\myB_{n}$, $\myC_{n}$, and 
$\myD_{n}$), this is handled by \FourFamiliesSecondLemma. 
For the exceptional graphs 
($\myE_{6}$, $\myE_{7}$, $\myE_{8}$, $\myF_{4}$, and $\myG_{2}$) this can be checked by 
hand (requiring 36, 63, 120, 24, and 6 firings respectively).  
For $\myG_{2}$, begin with strongly dominant position $\lambda = (a,b)$ 
with positive number $a$ on node $\gamma_{1}$ and positive 
number $b$ on 
$\gamma_{2}$.  
The following sequence of six node firings is easily seen 
to be legal: $(\gamma_{1}$, $\gamma_{2}$, $\gamma_{1}$, 
$\gamma_{2}$, $\gamma_{1}$, $\gamma_{2})$.  
The resulting position is $(-a,-b)$. 
For $\myF_{4}$, begin with strongly dominant position $(a,b,c,d)$.  
Play the numbers game from this initial position to see that the 
following sequence of 24 node firings is legal: 
$(\gamma_{1}$, 
$\gamma_{2}$, 
$\gamma_{3}$, 
$\gamma_{4}$, 
$\gamma_{3}$, 
$\gamma_{2}$, 
$\gamma_{1}$, 
$\gamma_{2}$, 
$\gamma_{3}$, 
$\gamma_{4}$, 
$\gamma_{2}$, 
$\gamma_{3}$, 
$\gamma_{2}$, 
$\gamma_{1}$, 
$\gamma_{4}$, 
$\gamma_{3}$, 
$\gamma_{2}$, 
$\gamma_{3}$, 
$\gamma_{4}$, 
$\gamma_{2}$, 
$\gamma_{1}$, 
$\gamma_{3}$, 
$\gamma_{2}$, 
$\gamma_{3})$.  
The resulting position is $(-a,-b,-c,-d)$. 
For $\myE_{6}$, begin with strongly dominant position $(a,b,c,d,e,f)$.  
Play the numbers game from this initial position to see that the 
following sequence of 36 node firings is legal: 
$(\gamma_{1}$, 
$\gamma_{2}$, 
$\gamma_{3}$, 
$\gamma_{4}$, 
$\gamma_{3}$, 
$\gamma_{2}$, 
$\gamma_{1}$, 
$\gamma_{4}$, 
$\gamma_{3}$, 
$\gamma_{4}$, 
$\gamma_{5}$, 
$\gamma_{4}$, 
$\gamma_{2}$, 
$\gamma_{3}$, 
$\gamma_{4}$, 
$\gamma_{1}$, 
$\gamma_{3}$, 
$\gamma_{5}$, 
$\gamma_{6}$, 
$\gamma_{4}$, 
$\gamma_{5}$, 
$\gamma_{4}$, 
$\gamma_{2}$, 
$\gamma_{3}$, 
$\gamma_{1}$, 
$\gamma_{4}$, 
$\gamma_{3}$, 
$\gamma_{1}$, 
$\gamma_{5}$, 
$\gamma_{4}$, 
$\gamma_{2}$, 
$\gamma_{6}$, 
$\gamma_{5}$, 
$\gamma_{4}$, 
$\gamma_{3}$, 
$\gamma_{1})$.  
The resulting position is $(-f,-b,-e,-d,-c,-a)$. 
For $\myE_{7}$, begin with strongly dominant position $(a,b,c,d,e,f,g)$.  
Play the numbers game from this initial position to see that the 
following sequence of 63 node firings is legal (the first 36 of these 
node firings are exactly the previous game sequence played on the 
$\myE_{6}$ subgraph): 
$(\gamma_{1}$, 
$\gamma_{2}$, 
$\gamma_{3}$, 
$\gamma_{4}$, 
$\gamma_{3}$, 
$\gamma_{2}$, 
$\gamma_{1}$, 
$\gamma_{4}$, 
$\gamma_{3}$, 
$\gamma_{4}$, 
$\gamma_{5}$, 
$\gamma_{4}$, 
$\gamma_{2}$, 
$\gamma_{3}$, 
$\gamma_{4}$, 
$\gamma_{1}$, 
$\gamma_{3}$, 
$\gamma_{5}$, 
$\gamma_{6}$, 
$\gamma_{4}$, 
$\gamma_{5}$, 
$\gamma_{4}$, 
$\gamma_{2}$, 
$\gamma_{3}$, 
$\gamma_{1}$, 
$\gamma_{4}$, 
$\gamma_{3}$, 
$\gamma_{1}$, 
$\gamma_{5}$, 
$\gamma_{4}$, 
$\gamma_{2}$, 
$\gamma_{6}$, 
$\gamma_{5}$, 
$\gamma_{4}$, 
$\gamma_{3}$, 
$\gamma_{1}$, 
$\gamma_{7}$, 
$\gamma_{6}$, 
$\gamma_{5}$, 
$\gamma_{4}$, 
$\gamma_{2}$, 
$\gamma_{3}$, 
$\gamma_{1}$, 
$\gamma_{4}$, 
$\gamma_{3}$, 
$\gamma_{5}$, 
$\gamma_{4}$, 
$\gamma_{2}$, 
$\gamma_{6}$, 
$\gamma_{5}$, 
$\gamma_{4}$, 
$\gamma_{3}$, 
$\gamma_{1}$, 
$\gamma_{7}$, 
$\gamma_{6}$, 
$\gamma_{5}$, 
$\gamma_{4}$, 
$\gamma_{2}$, 
$\gamma_{3}$, 
$\gamma_{4}$, 
$\gamma_{5}$, 
$\gamma_{6}$, 
$\gamma_{7})$.  
The resulting position is $(-a,-b,-c,-d,-e,-f,-g)$.
For $\myE_{8}$, begin with strongly dominant position $(a,b,c,d,e,f,g,h)$.  
Play the numbers game from this initial position to see that the 
following sequence of 120 node firings is legal (the first 63 of these 
node firings are exactly the previous game sequence played on the 
$\myE_{7}$ subgraph): 
$(\gamma_{1}$, 
$\gamma_{2}$, 
$\gamma_{3}$, 
$\gamma_{4}$, 
$\gamma_{3}$, 
$\gamma_{2}$, 
$\gamma_{1}$, 
$\gamma_{4}$, 
$\gamma_{3}$, 
$\gamma_{4}$, 
$\gamma_{5}$, 
$\gamma_{4}$, 
$\gamma_{2}$, 
$\gamma_{3}$, 
$\gamma_{4}$, 
$\gamma_{1}$, 
$\gamma_{3}$, 
$\gamma_{5}$, 
$\gamma_{6}$, 
$\gamma_{4}$, 
$\gamma_{5}$, 
$\gamma_{4}$, 
$\gamma_{2}$, 
$\gamma_{3}$, 
$\gamma_{1}$, 
$\gamma_{4}$, 
$\gamma_{3}$, 
$\gamma_{1}$, 
$\gamma_{5}$, 
$\gamma_{4}$, 
$\gamma_{2}$, 
$\gamma_{6}$, 
$\gamma_{5}$, 
$\gamma_{4}$, 
$\gamma_{3}$, 
$\gamma_{1}$, 
$\gamma_{7}$, 
$\gamma_{6}$, 
$\gamma_{5}$, 
$\gamma_{4}$, 
$\gamma_{2}$, 
$\gamma_{3}$, 
$\gamma_{1}$, 
$\gamma_{4}$, 
$\gamma_{3}$, 
$\gamma_{5}$, 
$\gamma_{4}$, 
$\gamma_{2}$, 
$\gamma_{6}$, 
$\gamma_{5}$, 
$\gamma_{4}$, 
$\gamma_{3}$, 
$\gamma_{1}$, 
$\gamma_{7}$, 
$\gamma_{6}$, 
$\gamma_{5}$, 
$\gamma_{4}$, 
$\gamma_{2}$, 
$\gamma_{3}$, 
$\gamma_{4}$, 
$\gamma_{5}$, 
$\gamma_{6}$, 
$\gamma_{7}$,  
$\gamma_{8}$,  
$\gamma_{7}$,  
$\gamma_{6}$,  
$\gamma_{5}$,  
$\gamma_{4}$,  
$\gamma_{2}$,  
$\gamma_{3}$,  
$\gamma_{1}$,  
$\gamma_{4}$,  
$\gamma_{3}$,  
$\gamma_{5}$,  
$\gamma_{4}$,  
$\gamma_{2}$,  
$\gamma_{6}$,  
$\gamma_{5}$,  
$\gamma_{4}$,  
$\gamma_{3}$,  
$\gamma_{1}$,  
$\gamma_{7}$,  
$\gamma_{6}$,  
$\gamma_{5}$,  
$\gamma_{4}$,  
$\gamma_{2}$,  
$\gamma_{3}$,  
$\gamma_{4}$,  
$\gamma_{5}$,  
$\gamma_{6}$,  
$\gamma_{7}$,  
$\gamma_{8}$,  
$\gamma_{7}$,  
$\gamma_{6}$,  
$\gamma_{5}$,  
$\gamma_{4}$,  
$\gamma_{2}$,  
$\gamma_{3}$,  
$\gamma_{1}$,  
$\gamma_{4}$,  
$\gamma_{3}$,  
$\gamma_{5}$,  
$\gamma_{4}$,  
$\gamma_{2}$,  
$\gamma_{6}$,  
$\gamma_{7}$,  
$\gamma_{5}$,  
$\gamma_{6}$,  
$\gamma_{4}$,  
$\gamma_{3}$,  
$\gamma_{1}$,  
$\gamma_{5}$,  
$\gamma_{4}$,  
$\gamma_{2}$,  
$\gamma_{3}$,  
$\gamma_{4}$,  
$\gamma_{5}$,  
$\gamma_{6}$,  
$\gamma_{7}$,  
$\gamma_{8})$.   
The resulting position is $(-a,-b,-c,-d,-e,-f,-g,-h)$.
\hfill\QED

\noindent 
{\bf \RemarkOnGameLength}\ \ As can be seen from 
\FourFamiliesSecondLemma\ and \ConvergenceProposition, 
the length of any game sequence played from 
any strongly dominant initial position on $\myA_{n}$ (respectively $\myB_{n}$, 
$\myC_{n}$, $\myD_{n}$) is $\frac{n(n+1)}{2}$ (respectively $n^{2}$, 
$n^{2}$, $n(n-1)$).  Similarly, from the 
statement and proof of \ConvergenceProposition\ it follows that the 
length of any game sequence played from any strongly dominant initial 
position on $\myE_{6}$ (respectively $\myE_{7}$, 
$\myE_{8}$, $\myF_{4}$, $\myG_{2}$) is $36$ (respectively $63$, 
$120$, $24$, $6$).\hfill\QED 

\vspace{2ex} 

\noindent
{\Large \bf \SecondProofNum.\ \ Divergent games for some families of graphs}

\vspace{1ex} 
\noindent 
{\bf \NotMarsFriendlyCatalog}\ \ {\sl The connected GCM graphs of 
\NotMarsFriendlyFigure\ are not admissible.} 

{\em Proof.} 
By \NotMarsFriendlyLemma\ it 
suffices to show that for each graph in 
\NotMarsFriendlyFigure\ and for each fundamental position, there is a 
divergent game sequence. 
In each case we 
exhibit a divergent 
game sequence which is a simple pattern of node firings.  
Remarkably, 
in all cases trial and error quickly lead us to 
these patterns. 
Our goal in this 
proof 
is not to develop any general theory for 
finding divergent game sequences for these cases, 
but rather to show that such game 
sequences can be found and presented in an elementary (though 
sometimes tedious) manner. 
The $\widetilde{\myA}$, $\widetilde{\myB}$, $\widetilde{\myC}$, 
$\widetilde{\myD}$, $\widetilde{\myE}$, and $\widetilde{\myF}$ cases are 
handled using a common line of reasoning: 
A sequence of legal node firings is applied to a 
position whose numbers are linear expressions in an index 
variable $k$.  It is then observed 
that the numbers for the resulting 
position are linear expressions of the same form with respect to 
the variable $k+1$ and that the firing sequence can be repeated.  
The $\widetilde{\myG}$ cases and the families of 
small cycles are 
handled using a variation of this kind of argument: A sequence of 
legal node firings is applied to a generic position satisfying certain 
inequalities, and it is shown that the resulting position also 
satisfies these inequalities so that the firing sequence can be 
repeated. 
Each paragraph in what follows 
demonstrates inadmissibility for some graph in our list. 
Our case-analysis argument is lengthy in part because we have tried 
to make each paragraph reasonably self-contained. 


\begin{figure}[ht]
\begin{center}
\NotMarsFriendlyFigure: Some connected GCM graphs that are not 
admissible.\\  (Figure continues on the next page.) 
\end{center}

\hspace*{0.75in}
\fbox{The $\widetilde{\myA}$ family} 
\hspace*{0.5in}
\parbox[c]{3.1in}{
\setlength{\unitlength}{0.5in}
\begin{picture}(6.1,0.7)
            \put(0,0.1){\circle*{0.075}}
            \put(1,0.1){\circle*{0.075}}
            \put(2,0.1){\circle*{0.075}}
            \put(3,0.6){\circle*{0.075}}
            \put(4,0.1){\circle*{0.075}}
            \put(5,0.1){\circle*{0.075}}
            \put(6,0.1){\circle*{0.075}}
            \put(0,0.1){\line(1,0){2}}
            \put(0,0.1){\line(6,1){3}}
            \multiput(2,0.1)(0.4,0){5}{\line(1,0){0.2}}
            \put(4,0.1){\line(1,0){2}}
            \put(6,0.1){\line(-6,1){3}}
           \end{picture}}

\vspace*{0.25in}
\hspace*{0.75in}
\fbox{The $\widetilde{\myB}$ family} 
\hspace*{0.52in}
\parbox[c]{3.1in}{
\setlength{\unitlength}{0.5in}
\begin{picture}(4.1,0.7)
            \put(0,0.6){\circle*{0.075}}
            \put(0,0.1){\circle*{0.075}}
            \put(1,0.35){\circle*{0.075}}
            \put(2,0.35){\circle*{0.075}}
            \put(4,0.35){\circle*{0.075}}
            \put(5,0.35){\circle*{0.075}}
            \put(6,0.35){\circle*{0.075}}
            \put(0,0.1){\line(4,1){1}}
            \put(0,0.6){\line(4,-1){1}}
            \put(1,0.35){\line(1,0){1}}
            \multiput(2,0.35)(0.4,0){5}{\line(1,0){0.2}}
            \put(4,0.35){\line(1,0){2}}
            \put(5.8,0.35){\vector(-1,0){0.1}}
            \put(5.2,0.35){\vector(1,0){0.1}}
            \put(5.3,0.35){\vector(1,0){0.1}}
           \end{picture}}

\vspace*{0.1in}
\hspace*{0.75in}
\hspace*{1.542in}
\parbox[c]{3.1in}{
\setlength{\unitlength}{0.5in}
\begin{picture}(4.1,0.2)
            \put(0,0.1){\circle*{0.075}}
            \put(1,0.1){\circle*{0.075}}
            \put(2,0.1){\circle*{0.075}}
            \put(4,0.1){\circle*{0.075}}
            \put(5,0.1){\circle*{0.075}}
            \put(6,0.1){\circle*{0.075}}
            \put(0,0.1){\line(1,0){2}}
            \multiput(2,0.1)(0.4,0){5}{\line(1,0){0.2}}
            \put(4,0.1){\line(1,0){2}}
            \put(0.2,0.1){\vector(1,0){0.1}}
            \put(0.8,0.1){\vector(-1,0){0.1}}
            \put(0.7,0.1){\vector(-1,0){0.1}}
            \put(5.8,0.1){\vector(-1,0){0.1}}
            \put(5.2,0.1){\vector(1,0){0.1}}
            \put(5.3,0.1){\vector(1,0){0.1}}
           \end{picture}} 

\vspace*{0.25in}
\hspace*{0.75in}
\fbox{The $\widetilde{\myC}$ family} 
\hspace*{0.525in}
\parbox[c]{3.1in}{
\setlength{\unitlength}{0.5in}
\begin{picture}(4.1,0.2)
            \put(0,0.1){\circle*{0.075}}
            \put(1,0.1){\circle*{0.075}}
            \put(2,0.1){\circle*{0.075}}
            \put(4,0.1){\circle*{0.075}}
            \put(5,0.1){\circle*{0.075}}
            \put(6,0.1){\circle*{0.075}}
            \put(0,0.1){\line(1,0){2}}
            \multiput(2,0.1)(0.4,0){5}{\line(1,0){0.2}}
            \put(4,0.1){\line(1,0){2}}
            \put(0.2,0.1){\vector(1,0){0.1}}
            \put(0.3,0.1){\vector(1,0){0.1}}
            \put(0.8,0.1){\vector(-1,0){0.1}}
            \put(5.8,0.1){\vector(-1,0){0.1}}
            \put(5.7,0.1){\vector(-1,0){0.1}}
            \put(5.2,0.1){\vector(1,0){0.1}}
           \end{picture}}  

\vspace*{0.1in}
\hspace*{0.75in}
\hspace*{1.542in}
\parbox[c]{3.1in}{
\setlength{\unitlength}{0.5in}
\begin{picture}(4.1,0.2)
            \put(0,0.1){\circle*{0.075}}
            \put(1,0.1){\circle*{0.075}}
            \put(2,0.1){\circle*{0.075}}
            \put(4,0.1){\circle*{0.075}}
            \put(5,0.1){\circle*{0.075}}
            \put(6,0.1){\circle*{0.075}}
            \put(0,0.1){\line(1,0){2}}
            \multiput(2,0.1)(0.4,0){5}{\line(1,0){0.2}}
            \put(4,0.1){\line(1,0){2}}
            \put(0.2,0.1){\vector(1,0){0.1}}
            \put(0.8,0.1){\vector(-1,0){0.1}}
            \put(0.7,0.1){\vector(-1,0){0.1}}
            \put(5.8,0.1){\vector(-1,0){0.1}}
            \put(5.7,0.1){\vector(-1,0){0.1}}
            \put(5.2,0.1){\vector(1,0){0.1}}
           \end{picture}}  

\vspace*{0.1in}
\hspace*{0.75in}
\hspace*{1.542in}
\parbox[c]{3.1in}{
\setlength{\unitlength}{0.5in}
\begin{picture}(4.1,0.7)
            \put(0,0.6){\circle*{0.075}}
            \put(0,0.1){\circle*{0.075}}
            \put(1,0.35){\circle*{0.075}}
            \put(2,0.35){\circle*{0.075}}
            \put(4,0.35){\circle*{0.075}}
            \put(5,0.35){\circle*{0.075}}
            \put(6,0.35){\circle*{0.075}}
            \put(0,0.1){\line(4,1){1}}
            \put(0,0.6){\line(4,-1){1}}
            \put(1,0.35){\line(1,0){1}}
            \multiput(2,0.35)(0.4,0){5}{\line(1,0){0.2}}
            \put(4,0.35){\line(1,0){2}}
            \put(5.2,0.35){\vector(1,0){0.1}}
            \put(5.8,0.35){\vector(-1,0){0.1}}
            \put(5.7,0.35){\vector(-1,0){0.1}}
           \end{picture}}

\vspace*{0.25in}
\hspace*{0.75in}
\fbox{The $\widetilde{\myD}$ family} 
\hspace*{0.5in}
\parbox[c]{3.1in}{
\setlength{\unitlength}{0.5in}
\begin{picture}(4.1,0.7)
            \put(0,0.6){\circle*{0.075}}
            \put(0,0.1){\circle*{0.075}}
            \put(1,0.35){\circle*{0.075}}
            \put(2,0.35){\circle*{0.075}}
            \put(4,0.35){\circle*{0.075}}
            \put(5,0.35){\circle*{0.075}}
            \put(6,0.1){\circle*{0.075}}
            \put(6,0.6){\circle*{0.075}}
            \put(0,0.1){\line(4,1){1}}
            \put(0,0.6){\line(4,-1){1}}
            \put(1,0.35){\line(1,0){1}}
            \multiput(2,0.35)(0.4,0){5}{\line(1,0){0.2}}
            \put(4,0.35){\line(1,0){1}}
            \put(5,0.35){\line(4,1){1}}
            \put(5,0.35){\line(4,-1){1}}
           \end{picture}} 

\vspace*{0.25in}
\fbox{The $\widetilde{\myE}$ family}\\ 
\hspace*{0.5in}
\parbox[c]{2.1in}{
\setlength{\unitlength}{0.5in}
\begin{picture}(4.1,1.2)
            \put(0,0.1){\circle*{0.075}}
            \put(1,0.1){\circle*{0.075}}
            \put(2,0.1){\circle*{0.075}}
            \put(2,0.6){\circle*{0.075}}
            \put(2,1.1){\circle*{0.075}}
            \put(3,0.1){\circle*{0.075}}
            \put(4,0.1){\circle*{0.075}}
            \put(0,0.1){\line(1,0){4}}
            \put(2,0.1){\line(0,1){1}}
           \end{picture}}
\hspace*{0.25in} 
\parbox[c]{3.1in}{
\setlength{\unitlength}{0.5in}
\begin{picture}(4.1,0.8)
            \put(0,0.1){\circle*{0.075}}
            \put(1,0.1){\circle*{0.075}}
            \put(2,0.1){\circle*{0.075}}
            \put(3,0.1){\circle*{0.075}}
            \put(3,0.6){\circle*{0.075}}
            \put(4,0.1){\circle*{0.075}}
            \put(5,0.1){\circle*{0.075}}
            \put(6,0.1){\circle*{0.075}}
            \put(0,0.1){\line(1,0){6}}
            \put(3,0.1){\line(0,1){0.5}}
           \end{picture}}

\vspace*{0.1in}
\hspace*{2.0in} 
\parbox[c]{3.1in}{
\setlength{\unitlength}{0.5in}
\begin{picture}(4.1,0.8)
            \put(0,0.1){\circle*{0.075}}
            \put(1,0.1){\circle*{0.075}}
            \put(2,0.1){\circle*{0.075}}
            \put(2,0.6){\circle*{0.075}}
            \put(3,0.1){\circle*{0.075}}
            \put(4,0.1){\circle*{0.075}}
            \put(5,0.1){\circle*{0.075}}
            \put(6,0.1){\circle*{0.075}}
            \put(7,0.1){\circle*{0.075}}
            \put(0,0.1){\line(1,0){7}}
            \put(2,0.1){\line(0,1){0.5}}
           \end{picture}}

\vspace*{0.25in}
\fbox{The $\widetilde{\myF}$ family}  
\hspace*{0.25in}
\parbox[c]{2.1in}{
\setlength{\unitlength}{0.5in}
\begin{picture}(4.1,0.2)
            \put(0,0.1){\circle*{0.075}}
            \put(1,0.1){\circle*{0.075}}
            \put(2,0.1){\circle*{0.075}}
            \put(3,0.1){\circle*{0.075}}
            \put(4,0.1){\circle*{0.075}}
            \put(0,0.1){\line(1,0){4}}
            \put(1.2,0.1){\vector(1,0){0.1}}
            \put(1.3,0.1){\vector(1,0){0.1}}
            \put(1.8,0.1){\vector(-1,0){0.1}}
            \end{picture}}
\hspace*{0.5in}
\parbox[c]{2.1in}{
\setlength{\unitlength}{0.5in}
\begin{picture}(4.1,0.2)
            \put(0,0.1){\circle*{0.075}}
            \put(1,0.1){\circle*{0.075}}
            \put(2,0.1){\circle*{0.075}}
            \put(3,0.1){\circle*{0.075}}
            \put(4,0.1){\circle*{0.075}}
            \put(0,0.1){\line(1,0){4}}
            \put(1.2,0.1){\vector(1,0){0.1}}
            \put(1.8,0.1){\vector(-1,0){0.1}}
            \put(1.7,0.1){\vector(-1,0){0.1}}
            \end{picture}}
\end{figure} 

\noindent 
\fbox{The $\widetilde{\myA}$ family}\ \ 
The infinite $\widetilde{\myA}$ 
family of GCM graphs of \NotMarsFriendlyFigure\  
is the family of cycles with 
amplitude products of unity on all edges.  
Such cycles were in fact the graphs that motivated Mozes' study of the 
numbers game in \cite{Mozes}.  The argument we give here 
demonstrating inadmissibility for each graph in this family is a 
special case of the proof of Lemma 3.1 of \cite{DonEnumbers}. 
For an $n$-node graph 
%
%
in the $\widetilde{\myA}$ family (we take $n \geq 3$), 
number the top node $\gamma_{1}$ 
and the remaining nodes $\gamma_{2}$,$\ldots$,$\gamma_{n}$ in 
succession in the clockwise order around the cycle.  
For each fundamental position, we exhibit a divergent game sequence 
as a short sequence of legal 
node firings which can be repeated indefinitely. 
By symmetry, it suffices to do so for the 
fundamental position $\omega_{1} = (1,0,\ldots,0)$.  
This is the 
$k=0$ version of the position $(2k+1,-k,0,\ldots,0,-k)$.  From any such 
position with $k \geq 0$, 
the following sequence of node firings is easily seen to be legal: 
$(\gamma_{1},$ $\gamma_{2},$ $\ldots$ , $\gamma_{n-1}$, $\gamma_{n}$, 
$\gamma_{n-1}$, $\ldots$ , $\gamma_{3}$, 
$\gamma_{2})$.  
This sequence results in the position 
$(2(k+1)+1,-(k+1),0,\ldots,0,-(k+1))$. 
This gives the desired divergent 
game sequence.  We conclude that any such GCM graph is inadmissible.


\begin{figure}[ht]
\begin{center}
\NotMarsFriendlyFigure\ (continued): Some connected GCM graphs that are not 
admissible. 

\vspace*{0.25in}
\fbox{The $\widetilde{\myG}$ family}  
\hspace*{0.25in}
\parbox[c]{1.1in}{
\setlength{\unitlength}{0.5in}
\begin{picture}(2.1,0.2)
            \put(0,0.1){\circle*{0.075}}
            \put(1,0.1){\circle*{0.075}}
            \put(2,0.1){\circle*{0.075}}
            \put(0,0.1){\line(1,0){2}}
            \put(0.2,0.1){\vector(1,0){0.1}}
            \put(0.8,0.1){\vector(-1,0){0.1}}
            \put(0.7,0.1){\vector(-1,0){0.1}}
            \put(0.6,0.1){\vector(-1,0){0.1}}
            \end{picture}}
\hspace*{0.5in}
\parbox[c]{1.1in}{
\setlength{\unitlength}{0.5in}
\begin{picture}(2.1,0.2)
            \put(0,0.1){\circle*{0.075}}
            \put(1,0.1){\circle*{0.075}}
            \put(2,0.1){\circle*{0.075}}
            \put(0,0.1){\line(1,0){2}}
            \put(0.2,0.1){\vector(1,0){0.1}}
            \put(0.8,0.1){\vector(-1,0){0.1}}
            \put(0.7,0.1){\vector(-1,0){0.1}}
            \put(0.6,0.1){\vector(-1,0){0.1}}
            \put(1.2,0.1){\vector(1,0){0.1}}
            \put(1.8,0.1){\vector(-1,0){0.1}}
            \put(1.7,0.1){\vector(-1,0){0.1}}
            \end{picture}}
\hspace*{0.5in}
\parbox[c]{1.1in}{
\setlength{\unitlength}{0.5in}
\begin{picture}(2.1,0.2)
            \put(0,0.1){\circle*{0.075}}
            \put(1,0.1){\circle*{0.075}}
            \put(2,0.1){\circle*{0.075}}
            \put(0,0.1){\line(1,0){2}}
            \put(0.2,0.1){\vector(1,0){0.1}}
            \put(0.8,0.1){\vector(-1,0){0.1}}
            \put(0.7,0.1){\vector(-1,0){0.1}}
            \put(0.6,0.1){\vector(-1,0){0.1}}
            \put(1.2,0.1){\vector(1,0){0.1}}
            \put(1.8,0.1){\vector(-1,0){0.1}}
            \put(1.7,0.1){\vector(-1,0){0.1}}
            \put(1.6,0.1){\vector(-1,0){0.1}}
            \end{picture}}

\vspace*{0.1in}
\hspace*{1.28in}
\parbox[c]{1.1in}{
\setlength{\unitlength}{0.5in}
\begin{picture}(2.1,0.2)
            \put(0,0.1){\circle*{0.075}}
            \put(1,0.1){\circle*{0.075}}
            \put(2,0.1){\circle*{0.075}}
            \put(0,0.1){\line(1,0){2}}
            \put(0.2,0.1){\vector(1,0){0.1}}
            \put(0.3,0.1){\vector(1,0){0.1}}
            \put(0.4,0.1){\vector(1,0){0.1}}
            \put(0.8,0.1){\vector(-1,0){0.1}}
            \end{picture}}
\hspace*{0.5in}
\parbox[c]{1.1in}{
\setlength{\unitlength}{0.5in}
\begin{picture}(2.1,0.2)
            \put(0,0.1){\circle*{0.075}}
            \put(1,0.1){\circle*{0.075}}
            \put(2,0.1){\circle*{0.075}}
            \put(0,0.1){\line(1,0){2}}
            \put(0.2,0.1){\vector(1,0){0.1}}
            \put(0.8,0.1){\vector(-1,0){0.1}}
            \put(0.7,0.1){\vector(-1,0){0.1}}
            \put(0.6,0.1){\vector(-1,0){0.1}}
            \put(1.8,0.1){\vector(-1,0){0.1}}
            \put(1.2,0.1){\vector(1,0){0.1}}
            \put(1.3,0.1){\vector(1,0){0.1}}
            \end{picture}} 
\hspace*{0.5in}
\parbox[c]{1.1in}{
\setlength{\unitlength}{0.5in}
\begin{picture}(2.1,0.2)
            \put(0,0.1){\circle*{0.075}}
            \put(1,0.1){\circle*{0.075}}
            \put(2,0.1){\circle*{0.075}}
            \put(0,0.1){\line(1,0){2}}
            \put(0.2,0.1){\vector(1,0){0.1}}
            \put(0.8,0.1){\vector(-1,0){0.1}}
            \put(0.7,0.1){\vector(-1,0){0.1}}
            \put(0.6,0.1){\vector(-1,0){0.1}}
            \put(1.8,0.1){\vector(-1,0){0.1}}
            \put(1.2,0.1){\vector(1,0){0.1}}
            \put(1.3,0.1){\vector(1,0){0.1}}
            \put(1.4,0.1){\vector(1,0){0.1}}
            \end{picture}}

\vspace*{0.25in}
\fbox{Families of small cycles}\\
\SmallCycles  
\end{center}
\end{figure}


\noindent 
\fbox{The $\widetilde{\myB}$ family}\ \ 
First, we show why  
\hspace*{0.025in}
\parbox[c]{1.1in}{
\setlength{\unitlength}{0.5in}
\begin{picture}(2.1,0.7)
            \put(0,0.6){\circle*{0.075}}
            \put(0,0.1){\circle*{0.075}}
            \put(1,0.35){\circle*{0.075}}
            \put(2,0.35){\circle*{0.075}}
            \put(0,0.1){\line(4,1){1}}
            \put(0,0.6){\line(4,-1){1}}
            \put(1,0.35){\line(1,0){1}}
            \put(1.8,0.35){\vector(-1,0){0.1}}
            \put(1.2,0.35){\vector(1,0){0.1}}
            \put(1.3,0.35){\vector(1,0){0.1}}
           \end{picture}}
%
is not admissible.  
Label the leftmost nodes as $\gamma_{1}$ and $\gamma_{2}$, the 
middle node as $\gamma_{3}$, and the rightmost node as $\gamma_{4}$.  
For each fundamental position, we exhibit a divergent game sequence 
as a short sequence of legal 
node firings which can be repeated indefinitely. 
From the fundamental position $\omega_{1} = (1,0,0,0)$, 
play the (legal) sequence $(\gamma_{1}$, $\gamma_{3}$, $\gamma_{2}$, 
$\gamma_{4}$, $\gamma_{3})$ to obtain the position 
$(2,1,-2,2)$.  This is the 
$k=0$ version of the position $(k+2,k+1,-2(k+1),2(k+1))$.  From any such 
position with $k \geq 0$, 
the following sequence of node firings is easily seen to be legal: 
$(\gamma_{1}$, $\gamma_{2}$, $\gamma_{4}$, $\gamma_{3}$, 
$\gamma_{1}$, $\gamma_{2}$, $\gamma_{4}$, $\gamma_{3})$.  
This sequence results in the position 
$((k+1)+2,(k+1)+1,-2[(k+1)+1],2[(k+1)+1])$. By symmetry, we also obtain 
a divergent game sequence from the fundamental position 
$\omega_{2}$. 
The fundamental position $\omega_{3} = (0,0,1,0)$ is the 
$k=0$ version of the position $(-2k,-2k,2k+1,0)$.  From any such 
position with $k \geq 0$, 
the following sequence of node firings is easily seen to be legal: 
$(\gamma_{3}$, $\gamma_{4}$, $\gamma_{3}$, $\gamma_{2}$, $\gamma_{1})$.  
This sequence results in the position 
$(-2(k+1),-2(k+1),2(k+1)+1,0)$. 
The fundamental position $\omega_{4} = (0,0,0,1)$ is the 
$k=0$ version of the position $(0,0,-k,2k+1)$.  From any such 
position with $k \geq 0$, 
the following sequence of node firings is easily seen to be legal: 
$(\gamma_{4}$, $\gamma_{3}$, $\gamma_{2}$, $\gamma_{1}$, $\gamma_{3})$.  
This sequence results in the position 
$(0,0,-(k+1),2(k+1)+1)$. 

Next, we show why  
\hspace*{0.025in}
\parbox[c]{3.1in}{
\setlength{\unitlength}{0.5in}
\begin{picture}(4.1,0.7)
            \put(0,0.6){\circle*{0.075}}
            \put(0,0.1){\circle*{0.075}}
            \put(1,0.35){\circle*{0.075}}
            \put(2,0.35){\circle*{0.075}}
            \put(4,0.35){\circle*{0.075}}
            \put(5,0.35){\circle*{0.075}}
            \put(6,0.35){\circle*{0.075}}
            \put(0,0.1){\line(4,1){1}}
            \put(0,0.6){\line(4,-1){1}}
            \put(1,0.35){\line(1,0){1}}
            \multiput(2,0.35)(0.4,0){5}{\line(1,0){0.2}}
            \put(4,0.35){\line(1,0){2}}
            \put(5.8,0.35){\vector(-1,0){0.1}}
            \put(5.2,0.35){\vector(1,0){0.1}}
            \put(5.3,0.35){\vector(1,0){0.1}}
           \end{picture}}
%
is not admissible when the graph has $n \geq 5$ nodes.  
Label the leftmost nodes as $\gamma_{1}$ and $\gamma_{2}$, and label 
the remaining nodes in succession from left to right as $\gamma_{3}$, 
$\ldots$ , $\gamma_{n-1}$, $\gamma_{n}$.  
For each fundamental position, we exhibit a divergent game sequence 
as a short sequence of legal 
node firings which can be repeated indefinitely. 
The fundamental position $\omega_{1} = (1,0,\ldots,0)$ is the 
$k=0$ version of the position $(2k+1,-2k,0,\ldots,0)$.  From any such 
position with $k \geq 0$, 
the following sequence of node firings is easily seen to be legal: 
$(\gamma_{1}$, $\gamma_{3}$, $\gamma_{4}$, $\ldots$ , $\gamma_{n-1}$, 
$\gamma_{n}$, 
$\gamma_{n-1}$, $\ldots$ , $\gamma_{3}$, 
$\gamma_{2}$, 
$\gamma_{1}$, $\gamma_{3}$, $\gamma_{4}$, $\ldots$ , $\gamma_{n-1}$, 
$\gamma_{n}$, 
$\gamma_{n-1}$, $\ldots$ , $\gamma_{3}$, 
$\gamma_{2})$.  
This sequence results in the position 
$(2(k+1)+1,-2(k+1),0,\ldots,0)$. By symmetry, we also obtain 
a divergent game sequence from the fundamental position 
$\omega_{2}$. 
The fundamental position $\omega_{3} = (0,0,1,0,\ldots,0)$ is the 
$k=0$ version of the position $(-2k,-2k,2k+1,0,\ldots,0)$.  From any such 
position with $k \geq 0$, 
the following sequence of node firings is easily seen to be legal: 
$(\gamma_{3}$, $\gamma_{4}$, $\ldots$ , $\gamma_{n-1}$, $\gamma_{n}$, 
$\gamma_{n-1}$, $\ldots$ , $\gamma_{3}$, 
$\gamma_{2}$, $\gamma_{1})$.  
This sequence results in the position 
$(-2(k+1),-2(k+1),2(k+1)+1,0,\ldots,0)$. 
For $4 \leq i \leq n-1$, 
any fundamental position $\omega_{i} = (0,\ldots,0,1,0,\ldots,0)$ is the 
$k=0$ version of the position $(0,\ldots,0,-2k,2k+1,0,\ldots,0)$.  
From any such 
position with $k \geq 0$, 
the following sequence of node firings is easily seen to be legal: 
$(\gamma_{i}$, $\gamma_{i+1}$, $\ldots$ , $\gamma_{n-1}$, $\gamma_{n}$, 
$\gamma_{n-1}$, 
$\ldots$ , $\gamma_{i+1}$, $\gamma_{i}$, 
$\gamma_{i-1}$, $\ldots$ , $\gamma_{3}$, 
$\gamma_{2}$, $\gamma_{1}$, $\gamma_{3}$, $\ldots$ , $\gamma_{i-1})$.  
This sequence results in the position 
$(0,\ldots,0,-2(k+1),2(k+1)+1,0,\ldots,0)$. 
The fundamental position $\omega_{n} = (0,\ldots,0,1)$ is the 
$k=0$ version of the position $(0,\ldots,0,-k,2k+1)$.  
From any such 
position with $k \geq 0$, 
the following sequence of node firings is easily seen to be legal: 
$(\gamma_{n}$, 
$\gamma_{n-1}$, 
$\ldots$ ,  $\gamma_{3}$, 
$\gamma_{2}$, $\gamma_{1}$, $\gamma_{3}$, $\ldots$ , $\gamma_{n-1})$.  
This sequence results in the position 
$(0,\ldots,0,-(k+1),2(k+1)+1)$. 

Next, we show why  
\hspace*{0.025in}
\parbox[c]{1.1in}{
\setlength{\unitlength}{0.5in}
\begin{picture}(2.1,0.2)
            \put(0,0.1){\circle*{0.075}}
            \put(1,0.1){\circle*{0.075}}
            \put(2,0.1){\circle*{0.075}}
            \put(0,0.1){\line(1,0){2}}
            \put(0.2,0.1){\vector(1,0){0.1}}
            \put(0.8,0.1){\vector(-1,0){0.1}}
            \put(0.7,0.1){\vector(-1,0){0.1}}
            \put(1.8,0.1){\vector(-1,0){0.1}}
            \put(1.2,0.1){\vector(1,0){0.1}}
            \put(1.3,0.1){\vector(1,0){0.1}}
           \end{picture}} 
is not admissible.  
Label the nodes as $\gamma_{1}$, $\gamma_{2}$, and $\gamma_{3}$ 
from left to right.  
For each fundamental position, we exhibit a divergent game sequence 
as a short sequence of legal 
node firings which can be repeated indefinitely. 
The fundamental position $\omega_{1} = (1,0,0)$ is the 
$k=0$ version of the position $(2k+1,-k,0)$.  From any such 
position with $k \geq 0$, 
the following sequence of node firings is easily seen to be legal: 
$(\gamma_{1}$, $\gamma_{2}$, $\gamma_{3}$, $\gamma_{2})$.  
This sequence results in the position 
$(2(k+1)+1,-(k+1),0)$. 
By symmetry, we also obtain  
a divergent game sequence from the fundamental position 
$\omega_{3}$. 
The fundamental position $\omega_{2} = (0,1,0)$ is the 
$k=0$ version of the position $(-4k,2k+1,0)$.  From any such 
position with $k \geq 0$, 
the following sequence of node firings is easily seen to be legal: 
$(\gamma_{2}$, $\gamma_{3}$, $\gamma_{2}$, $\gamma_{1})$.  
This sequence results in the position 
$(-4(k+1),2(k+1)+1,0)$. 

We finish the $\widetilde{\myB}$ family by showing why  
\hspace*{0.025in}
\parbox[c]{3.1in}{
\setlength{\unitlength}{0.5in}
\begin{picture}(4.1,0.2)
            \put(0,0.1){\circle*{0.075}}
            \put(1,0.1){\circle*{0.075}}
            \put(2,0.1){\circle*{0.075}}
            \put(4,0.1){\circle*{0.075}}
            \put(5,0.1){\circle*{0.075}}
            \put(6,0.1){\circle*{0.075}}
            \put(0,0.1){\line(1,0){2}}
            \multiput(2,0.1)(0.4,0){5}{\line(1,0){0.2}}
            \put(4,0.1){\line(1,0){2}}
            \put(0.2,0.1){\vector(1,0){0.1}}
            \put(0.8,0.1){\vector(-1,0){0.1}}
            \put(0.7,0.1){\vector(-1,0){0.1}}
            \put(5.8,0.1){\vector(-1,0){0.1}}
            \put(5.2,0.1){\vector(1,0){0.1}}
            \put(5.3,0.1){\vector(1,0){0.1}}
           \end{picture}} 
is not admissible when the graph has $n \geq 4$ nodes.  
Label the nodes as $\gamma_{1}$, $\gamma_{2}$, 
$\ldots$ , $\gamma_{n-1}$, and $\gamma_{n}$ from left to right.  
For each fundamental position, we exhibit a divergent game sequence 
as a short sequence of legal 
node firings which can be repeated indefinitely. 
The fundamental position $\omega_{1} = (1,0,\ldots,0)$ is the 
$k=0$ version of the position $(2k+1,-k,0,\ldots,0)$.  From any such 
position with $k \geq 0$, 
the following sequence of node firings is easily seen to be legal: 
$(\gamma_{1}$, $\gamma_{2}$, $\ldots$ , $\gamma_{n-1}$, $\gamma_{n}$, 
$\gamma_{n-1}$, $\ldots$ , $\gamma_{3}$, 
$\gamma_{2})$.  
This sequence results in the position 
$(2(k+1)+1,-2(k+1),0,\ldots,0)$. 
By symmetry, we also obtain  
a divergent game sequence from the fundamental position 
$\omega_{n}$. 
For $2 \leq i \leq n-1$, 
any fundamental position $\omega_{i} = (0,\ldots,0,1,0,\ldots,0)$ is the 
$k=0$ version of the position $(0,\ldots,0,2k+1,-2k,0,\ldots,0)$.  From any such 
position with $k \geq 0$, 
the following sequence of node firings is easily seen to be legal: 
$(\gamma_{i}$, $\gamma_{i-1}$, $\ldots$ , $\gamma_{2}$, $\gamma_{1}$, 
$\gamma_{2}$, 
$\ldots$ , $\gamma_{n-1}$, $\gamma_{n}$, 
$\gamma_{n-1}$, $\ldots$ , $\gamma_{i+2}$, 
$\gamma_{i+1})$.  
This sequence results in the position 
$(0,\ldots,0,2(k+1)+1,-2(k+1),0,\ldots,0)$.

\noindent 
\fbox{The $\widetilde{\myC}$ family}\ \ 
First, we show why 
\hspace*{0.025in}
\parbox[c]{3.1in}{
\setlength{\unitlength}{0.5in}
\begin{picture}(4.1,0.2)
            \put(0,0.1){\circle*{0.075}}
            \put(1,0.1){\circle*{0.075}}
            \put(2,0.1){\circle*{0.075}}
            \put(4,0.1){\circle*{0.075}}
            \put(5,0.1){\circle*{0.075}}
            \put(6,0.1){\circle*{0.075}}
            \put(0,0.1){\line(1,0){2}}
            \multiput(2,0.1)(0.4,0){5}{\line(1,0){0.2}}
            \put(4,0.1){\line(1,0){2}}
            \put(0.2,0.1){\vector(1,0){0.1}}
            \put(0.3,0.1){\vector(1,0){0.1}}
            \put(0.8,0.1){\vector(-1,0){0.1}}
            \put(5.8,0.1){\vector(-1,0){0.1}}
            \put(5.7,0.1){\vector(-1,0){0.1}}
            \put(5.2,0.1){\vector(1,0){0.1}}
           \end{picture}}  
is not admissible when the graph has $n \geq 3$ nodes.  
(Since firing the middle node in the $n=3$ case 
is comparable to firing either $\gamma_{2}$ or $\gamma_{n-1}$ in the 
$n \geq 4$ cases, then the $n=3$ case 
does not need to be considered separately here.) 
Label the nodes as $\gamma_{1}$, $\gamma_{2}$, 
$\ldots$ , $\gamma_{n-1}$, and $\gamma_{n}$ from left to right.  
For each fundamental position, we exhibit a divergent game sequence 
as a short sequence of legal 
node firings which can be repeated indefinitely. 
The fundamental position $\omega_{1} = (1,0,\ldots,0)$ is the 
$k=0$ version of the position $(2k+1,-2k,0,\ldots,0)$.  From any such 
position with $k \geq 0$, 
the following sequence of node firings is easily seen to be legal: 
$(\gamma_{1}$, $\gamma_{2}$, $\ldots$ , $\gamma_{n-1}$, $\gamma_{n}$, 
$\gamma_{n-1}$, $\ldots$ , $\gamma_{3}$, 
$\gamma_{2})$.  
This sequence results in the position 
$(2(k+1)+1,-2(k+1),0,\ldots,0)$. 
By symmetry, we also obtain 
a divergent game sequence from the fundamental position 
$\omega_{n}$. 
For $2 \leq i \leq n-1$, 
any fundamental position $\omega_{i} = (0,\ldots,0,1,0,\ldots,0)$ is the 
$k=0$ version of the position $(0,\ldots,0,2k+1,-2k,0,\ldots,0)$.  From any such 
position with $k \geq 0$, 
the following sequence of node firings is easily seen to be legal: 
$(\gamma_{i}$, $\gamma_{i-1}$, $\ldots$ , $\gamma_{2}$, 
$\gamma_{1}$, $\gamma_{2}$, 
$\ldots$ , $\gamma_{n-1}$, $\gamma_{n}$, 
$\gamma_{n-1}$, $\ldots$ , $\gamma_{i+2}$, 
$\gamma_{i+1})$.  
This sequence results in the position 
$(0,\ldots,0,2(k+1)+1,-2(k+1),0,\ldots,0)$. 


Next, we show why 
\hspace*{0.025in}
\parbox[c]{3.1in}{
\setlength{\unitlength}{0.5in}
\begin{picture}(4.1,0.2)
            \put(0,0.1){\circle*{0.075}}
            \put(1,0.1){\circle*{0.075}}
            \put(2,0.1){\circle*{0.075}}
            \put(4,0.1){\circle*{0.075}}
            \put(5,0.1){\circle*{0.075}}
            \put(6,0.1){\circle*{0.075}}
            \put(0,0.1){\line(1,0){2}}
            \multiput(2,0.1)(0.4,0){5}{\line(1,0){0.2}}
            \put(4,0.1){\line(1,0){2}}
            \put(0.2,0.1){\vector(1,0){0.1}}
            \put(0.8,0.1){\vector(-1,0){0.1}}
            \put(0.7,0.1){\vector(-1,0){0.1}}
            \put(5.8,0.1){\vector(-1,0){0.1}}
            \put(5.7,0.1){\vector(-1,0){0.1}}
            \put(5.2,0.1){\vector(1,0){0.1}}
           \end{picture}}  
is not admissible when the graph has $n \geq 3$ nodes.  
(Since firing the middle node in the $n=3$ case 
is comparable to firing $\gamma_{2}$ in the 
$n \geq 4$ cases, then the $n=3$ case 
does not need to be considered separately here.) 
Label the nodes as $\gamma_{1}$, $\gamma_{2}$, 
$\ldots$ , $\gamma_{n-1}$, and $\gamma_{n}$ from left to right.  
For each fundamental position, we exhibit a divergent game sequence 
as a short sequence of legal 
node firings which can be repeated indefinitely. 
The fundamental position $\omega_{1} = (1,0,\ldots,0)$ is the 
$k=0$ version of the position $(2k+1,-k,0,\ldots,0)$.  From any such 
position with $k \geq 0$, 
the following sequence of node firings is easily seen to be legal: 
$(\gamma_{1}$, $\gamma_{2}$, $\ldots$ , $\gamma_{n-1}$, $\gamma_{n}$, 
$\gamma_{n-1}$, $\ldots$ , $\gamma_{3}$, 
$\gamma_{2})$.  
This sequence results in the position 
$(2(k+1)+1,-(k+1),0,\ldots,0)$. 
For $2 \leq i \leq n-1$, 
any fundamental position $\omega_{i} = (0,\ldots,0,1,0,\ldots,0)$ is the 
$k=0$ version of the position $(0,\ldots,0,2k+1,-2k,0,\ldots,0)$.  From any such 
position with $k \geq 0$, 
the following sequence of node firings is easily seen to be legal: 
$(\gamma_{i}$, $\gamma_{i-1}$, $\ldots$ , $\gamma_{2}$, 
$\gamma_{1}$, $\gamma_{2}$, 
$\ldots$ , $\gamma_{n-1}$, $\gamma_{n}$, 
$\gamma_{n-1}$, $\ldots$ , $\gamma_{i+2}$, 
$\gamma_{i+1})$.  
This sequence results in the position 
$(0,\ldots,0,2(k+1)+1,-2(k+1),0,\ldots,0)$. 
The fundamental position $\omega_{n} = (0,\ldots,0,1)$ is the 
$k=0$ version of the position $(0,\ldots,0,-2k,2k+1)$.  From any such 
position with $k \geq 0$, 
the following sequence of node firings is easily seen to be legal: 
$(\gamma_{n}$, $\gamma_{n-1}$, $\ldots$ , $\gamma_{2}$, $\gamma_{1}$, 
$\gamma_{2}$, $\ldots$ , $\gamma_{n-2}$, 
$\gamma_{n-1})$.  
This sequence results in the position 
$(0,\ldots,0,-2(k+1),2(k+1)+1)$. 


We finish the $\widetilde{\myC}$ family by showing why  
\hspace*{0.025in}
\parbox[c]{3.1in}{
\setlength{\unitlength}{0.5in}
\begin{picture}(4.1,0.7)
            \put(0,0.6){\circle*{0.075}}
            \put(0,0.1){\circle*{0.075}}
            \put(1,0.35){\circle*{0.075}}
            \put(2,0.35){\circle*{0.075}}
            \put(4,0.35){\circle*{0.075}}
            \put(5,0.35){\circle*{0.075}}
            \put(6,0.35){\circle*{0.075}}
            \put(0,0.1){\line(4,1){1}}
            \put(0,0.6){\line(4,-1){1}}
            \put(1,0.35){\line(1,0){1}}
            \multiput(2,0.35)(0.4,0){5}{\line(1,0){0.2}}
            \put(4,0.35){\line(1,0){2}}
            \put(5.2,0.35){\vector(1,0){0.1}}
            \put(5.8,0.35){\vector(-1,0){0.1}}
            \put(5.7,0.35){\vector(-1,0){0.1}}
           \end{picture}}
is not admissible when the graph has $n \geq 4$ nodes.  
(Since firing the middle node in the $n=4$ case 
is comparable to firing $\gamma_{3}$ in the 
$n \geq 5$ cases, then the $n=4$ case 
does not need to be considered separately here.) 
Label the leftmost nodes as $\gamma_{1}$ and $\gamma_{2}$, and label 
the remaining nodes in succession from left to right as $\gamma_{3}$, 
$\ldots$ , $\gamma_{n-1}$, $\gamma_{n}$.  
For each fundamental position, we exhibit a divergent game sequence 
as a short sequence of legal 
node firings which can be repeated indefinitely. 
The fundamental position $\omega_{1} = (1,0,\ldots,0)$ is the 
$k=0$ version of the position $(2k+1,-2k,0,\ldots,0)$.  From any such 
position with $k \geq 0$, 
the following sequence of node firings is easily seen to be legal: 
$(\gamma_{1}$, $\gamma_{3}$, $\gamma_{4}$, $\ldots$ , 
$\gamma_{n-1}$, $\gamma_{n}$, 
$\gamma_{n-1}$, $\ldots$ , $\gamma_{3}$, 
$\gamma_{2}$, 
$\gamma_{1}$, $\gamma_{3}$, $\gamma_{4}$, $\ldots$ , $\gamma_{n-1}$, $\gamma_{n}$, 
$\gamma_{n-1}$, $\ldots$ , $\gamma_{3}$, 
$\gamma_{2})$.  
This sequence results in the position 
$(2(k+1)+1,-2(k+1),0,\ldots,0)$. By symmetry, we also obtain 
a divergent game sequence from the fundamental position 
$\omega_{2}$. 
The fundamental position $\omega_{3} = (0,0,1,0,\ldots,0)$ is the 
$k=0$ version of the position $(-2k,-2k,2k+1,0,\ldots,0)$.  From any such 
position with $k \geq 0$, 
the following sequence of node firings is easily seen to be legal: 
$(\gamma_{3}$, $\gamma_{4}$, $\ldots$ , $\gamma_{n-1}$, $\gamma_{n}$, 
$\gamma_{n-1}$, $\ldots$ , $\gamma_{3}$, 
$\gamma_{2}$, $\gamma_{1})$.  
This sequence results in the position 
$(-2(k+1),-2(k+1),2(k+1)+1,0,\ldots,0)$. 
For $4 \leq i \leq n$, 
any fundamental position $\omega_{i} = (0,\ldots,0,1,0,\ldots,0)$ is the 
$k=0$ version of the position $(0,\ldots,0,-2k,2k+1,0,\ldots,0)$.  From any such 
position with $k \geq 0$, 
the following sequence of node firings is easily seen to be legal: 
$(\gamma_{i}$, $\gamma_{i+1}$, $\ldots$ , $\gamma_{n-1}$, $\gamma_{n}$, 
$\gamma_{n-1}$, 
$\ldots$ , $\gamma_{i+1}$, $\gamma_{i}$, 
$\gamma_{i-1}$, $\ldots$ , $\gamma_{3}$, 
$\gamma_{2}$, $\gamma_{1}$, $\gamma_{3}$, $\ldots$ , $\gamma_{i-1})$.  
This sequence results in the position 
$(0,\ldots,0,-2(k+1),2(k+1)+1,0,\ldots,0)$.

\noindent 
\fbox{The $\widetilde{\myD}$ family}\ \ 
First, we show why the five-node graph 
\hspace*{0.025in}
\parbox[c]{1.1in}{
\setlength{\unitlength}{0.5in}
\begin{picture}(2.1,0.7)
            \put(0,0.6){\circle*{0.075}}
            \put(0,0.1){\circle*{0.075}}
            \put(1,0.35){\circle*{0.075}}
            \put(2,0.1){\circle*{0.075}}
            \put(2,0.6){\circle*{0.075}}
            \put(0,0.1){\line(4,1){1}}
            \put(0,0.6){\line(4,-1){1}}
            \put(1,0.35){\line(4,1){1}}
            \put(1,0.35){\line(4,-1){1}}
           \end{picture}}
is not admissible.  Label the leftmost nodes as $\gamma_{1}$ and $\gamma_{2}$, label 
the middle node as $\gamma_{3}$, and label the rightmost nodes as 
$\gamma_{4}$ and $\gamma_{5}$.   
For each fundamental position, we exhibit a divergent game sequence 
as a short sequence of legal 
node firings which can be repeated indefinitely. 
From the fundamental position $\omega_{1} = (1,0,0,0,0)$, 
play the (legal) sequence $(\gamma_{1}$, $\gamma_{3}$, $\gamma_{2}$, 
$\gamma_{4}$, $\gamma_{5}$, $\gamma_{3})$ to obtain the position 
$(2,1,-2,1,1)$.  This is the 
$k=0$ version of the position $(k+2,k+1,-2(k+1),k+1,k+1)$.  From any such 
position with $k \geq 0$, 
the following sequence of node firings is easily seen to be legal: 
$(\gamma_{1}$, $\gamma_{2}$, $\gamma_{4}$, $\gamma_{5}$, $\gamma_{3}$, 
$\gamma_{1}$, $\gamma_{2}$, $\gamma_{4}$, $\gamma_{5}$, $\gamma_{3})$.  
This sequence results in the position 
$((k+1)+2,(k+1)+1,-2[(k+1)+1],(k+1)+1,(k+1)+1)$. By symmetry, we also obtain 
divergent game sequences from the fundamental positions 
$\omega_{2}$, $\omega_{4}$, and $\omega_{5}$. 
The fundamental position $\omega_{3} = (0,0,1,0,0)$ is the 
$k=0$ version of the position $(-k,-k,2k+1,-k,-k)$.  From any such 
position with $k \geq 0$, 
the following sequence of node firings is easily seen to be legal: 
$(\gamma_{3}$, $\gamma_{1}$, $\gamma_{2}$, 
$\gamma_{4}$, $\gamma_{5})$.  
This sequence results in the position 
$(-(k+1),-(k+1),2(k+1)+1,-(k+1),-(k+1))$. 

We finish the $\widetilde{\myD}$ family by showing why  
\hspace*{0.025in}
\parbox[c]{3.1in}{
\setlength{\unitlength}{0.5in}
\begin{picture}(4.1,0.7)
            \put(0,0.6){\circle*{0.075}}
            \put(0,0.1){\circle*{0.075}}
            \put(1,0.35){\circle*{0.075}}
            \put(2,0.35){\circle*{0.075}}
            \put(4,0.35){\circle*{0.075}}
            \put(5,0.35){\circle*{0.075}}
            \put(6,0.1){\circle*{0.075}}
            \put(6,0.6){\circle*{0.075}}
            \put(0,0.1){\line(4,1){1}}
            \put(0,0.6){\line(4,-1){1}}
            \put(1,0.35){\line(1,0){1}}
            \multiput(2,0.35)(0.4,0){5}{\line(1,0){0.2}}
            \put(4,0.35){\line(1,0){1}}
            \put(5,0.35){\line(4,1){1}}
            \put(5,0.35){\line(4,-1){1}}
           \end{picture}} 
%
is not admissible when the graph has $n \geq 6$ nodes.  
Label the leftmost nodes as $\gamma_{1}$ and $\gamma_{2}$, label 
the ``isthmus'' nodes in succession from left to right as $\gamma_{3}$, 
$\ldots$ , $\gamma_{n-2}$, and label the rightmost nodes as 
$\gamma_{n-1}$ and $\gamma_{n}$.  
For each fundamental position, we exhibit a divergent game sequence 
as a short sequence of legal 
node firings which can be repeated indefinitely. 
The fundamental position $\omega_{1} = (1,0,\ldots,0)$ is the 
$k=0$ version of the position $(2k+1,-2k,0,\ldots,0)$.  From any such 
position with $k \geq 0$, 
the following sequence of node firings is easily seen to be legal: 
$(\gamma_{1}$, $\gamma_{3}$, $\gamma_{4}$, $\ldots$ , $\gamma_{n-2}$, 
$\gamma_{n-1}$, $\gamma_{n}$, 
$\gamma_{n-2}$, $\ldots$ , $\gamma_{4}$, $\gamma_{3}$, 
$\gamma_{2}$, 
$\gamma_{1}$, $\gamma_{3}$, $\gamma_{4}$, $\ldots$ , $\gamma_{n-2}$, 
$\gamma_{n-1}$, $\gamma_{n}$, 
$\gamma_{n-2}$, $\ldots$ , $\gamma_{4}$, $\gamma_{3}$, 
$\gamma_{2})$.  
This sequence results in the position 
$(2(k+1)+1,-2(k+1),0,\ldots,0)$. By symmetry, we also obtain 
divergent game sequences from the fundamental positions 
$\omega_{2}$, $\omega_{n-1}$, and $\omega_{n}$. 
The fundamental position $\omega_{3} = (0,0,1,0\ldots,0)$ is the 
$k=0$ version of the position $(-2k,-2k,2k+1,0,\ldots,0)$.  From any such 
position with $k \geq 0$, 
the following sequence of node firings is easily seen to be legal: 
$(\gamma_{3}$, $\gamma_{4}$, $\ldots$ , $\gamma_{n-2}$, 
$\gamma_{n-1}$, $\gamma_{n}$, 
$\gamma_{n-2}$, $\ldots$ , $\gamma_{4}$, $\gamma_{3}$, 
$\gamma_{2}$, 
$\gamma_{1})$.  
This sequence results in the position 
$(-2(k+1),-2(k+1),2(k+1)+1,0,\ldots,0)$. By symmetry, we also obtain 
a divergent game sequence from the fundamental position 
$\omega_{n-2}$. 
For $4 \leq i \leq n-2$, 
any fundamental position $\omega_{i} = (0,\ldots,0,1,0,\ldots,0)$ is the 
$k=0$ version of the position $(0,\ldots,0,-k,2k+1,-k,0,\ldots,0)$.  From any such 
position with $k \geq 0$, 
the following sequence of node firings is easily seen to be legal: 
$(\gamma_{i}$, $\gamma_{i+1}$, $\ldots$ , $\gamma_{n-2}$, $\gamma_{n-1}$, 
$\gamma_{n}$, $\gamma_{n-2}$,  
$\ldots$ , $\gamma_{i+1}$,  
$\gamma_{i-1}$, $\ldots$ , $\gamma_{3}$, 
$\gamma_{2}$, $\gamma_{1}$, $\gamma_{3}$, $\ldots$ , $\gamma_{i-1})$.  
This sequence results in the position 
$(0,\ldots,0,-(k+1),2(k+1)+1,-(k+1),0,\ldots,0)$.

\noindent 
\fbox{The $\widetilde{\myE}$ family}\ \ 
First, we show why\\ 
\hspace*{2.5in} 
\parbox[c]{3.1in}{
\setlength{\unitlength}{0.5in}
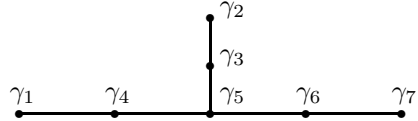
\begin{picture}(4.1,1.2)
            \put(0,0.1){\circle*{0.075}}
            \put(-0.1,0.25){\small $\gamma_{1}$}
            \put(1,0.1){\circle*{0.075}}
            \put(0.9,0.25){\small $\gamma_{4}$}
            \put(2,0.1){\circle*{0.075}}
            \put(2.1,0.25){\small $\gamma_{5}$}
            \put(2,0.6){\circle*{0.075}}
            \put(2.1,0.65){\small $\gamma_{3}$}
            \put(2,1.1){\circle*{0.075}}
            \put(2.1,1.15){\small $\gamma_{2}$}
            \put(3,0.1){\circle*{0.075}}
            \put(2.9,0.25){\small $\gamma_{6}$}
            \put(4,0.1){\circle*{0.075}}
            \put(3.9,0.25){\small $\gamma_{7}$}
            \put(0,0.1){\line(1,0){4}}
            \put(2,0.1){\line(0,1){1}}
           \end{picture}}

\vspace*{0.05in}
\noindent             
is not admissible.  
For each fundamental position, we exhibit a divergent game sequence 
as a short sequence of legal 
node firings which can be repeated indefinitely. 
The fundamental position $\omega_{1} = (1,0,0,0,0,0,0)$ is the 
$k=0$ version of the position $(2k+1,0,0,-k,0,0,0)$.  From any such 
position with $k \geq 0$, 
the following sequence of node firings is easily seen to be legal: 
$(\gamma_{1}$, $\gamma_{4}$, $\gamma_{5}$, $\gamma_{3}$, $\gamma_{2}$, 
$\gamma_{6}$, $\gamma_{5}$, $\gamma_{3}$, $\gamma_{4}$, $\gamma_{5}$, 
$\gamma_{6}$, $\gamma_{7}$, $\gamma_{6}$, $\gamma_{5}$, $\gamma_{3}$, 
$\gamma_{2}$, $\gamma_{4}$, $\gamma_{5}$, $\gamma_{3}$, $\gamma_{6}$, 
$\gamma_{5}$, $\gamma_{4})$.  
This results in the position $(2(k+1)+1,0,0,-(k+1),0,0,0)$. 
By symmetry, we also obtain 
divergent game sequences from the fundamental positions 
$\omega_{2}$ and $\omega_{7}$. 
The fundamental position $\omega_{4} = (0,0,0,1,0,0,0)$ is the 
$k=0$ version of the position $(-4k,0,0,2k+1,0,0,0)$.  From any such 
position with $k \geq 0$, 
the following sequence of node firings is easily seen to be legal: 
$(\gamma_{4}$, $\gamma_{5}$, $\gamma_{3}$, $\gamma_{2}$, $\gamma_{6}$, 
$\gamma_{5}$, $\gamma_{3}$, $\gamma_{4}$, $\gamma_{5}$, $\gamma_{6}$, 
$\gamma_{7}$, $\gamma_{6}$, $\gamma_{5}$, $\gamma_{3}$, $\gamma_{2}$, 
$\gamma_{4}$, $\gamma_{5}$, $\gamma_{3}$, $\gamma_{6}$, $\gamma_{5}$, 
$\gamma_{4}$, $\gamma_{1})$.  
This results in the position $(-4(k+1),0,0,2(k+1)+1,0,0,0)$. 
By symmetry, we also obtain 
divergent game sequences from the fundamental positions 
$\omega_{3}$ and $\omega_{6}$. 
The fundamental position $\omega_{5} = (0,0,0,0,1,0,0)$ is the 
$k=0$ version of the position $(-6k,0,0,-6k,6k+1,0,0)$.  
From any such 
position with $k \geq 0$, play the game to see that 
the following sequence of node firings is legal: 
$(\gamma_{5}$, $\gamma_{3}$, $\gamma_{2}$, $\gamma_{4}$, $\gamma_{5}$, 
$\gamma_{3}$, $\gamma_{6}$, $\gamma_{5}$, $\gamma_{3}$, $\gamma_{2}$, 
$\gamma_{4}$, $\gamma_{5}$, $\gamma_{3}$, $\gamma_{6}$, $\gamma_{5}$, 
$\gamma_{7}$, $\gamma_{6}$, $\gamma_{5}$, $\gamma_{3}$, $\gamma_{2}$, 
$\gamma_{4}$, $\gamma_{5}$, $\gamma_{3}$, $\gamma_{6}$, $\gamma_{5}$, 
$\gamma_{4}$, $\gamma_{7}$, $\gamma_{6}$, $\gamma_{5}$, $\gamma_{1}$, 
$\gamma_{4}$,
$\gamma_{5}$, $\gamma_{3}$, $\gamma_{2}$, $\gamma_{4}$, $\gamma_{5}$, 
$\gamma_{3}$, $\gamma_{6}$, $\gamma_{5}$, $\gamma_{3}$, $\gamma_{2}$, 
$\gamma_{4}$, $\gamma_{5}$, $\gamma_{3}$, $\gamma_{6}$, $\gamma_{5}$, 
$\gamma_{7}$, $\gamma_{6}$, $\gamma_{5}$, $\gamma_{3}$, $\gamma_{2}$, 
$\gamma_{4}$, $\gamma_{5}$, $\gamma_{3}$, $\gamma_{6}$, $\gamma_{5}$, 
$\gamma_{4}$, $\gamma_{7}$, $\gamma_{6}$, $\gamma_{5}$, $\gamma_{1}$, 
$\gamma_{4})$.   
This results in the position $(-6(k+1),0,0,-6(k+1),6(k+1)+1,0,0)$. 

Next, we show why\\ 
\hspace*{1.75in} 
\parbox[c]{3.1in}{
\setlength{\unitlength}{0.5in}
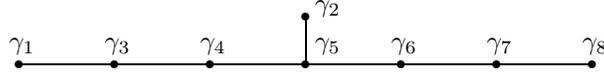
\begin{picture}(4.1,0.8)
            \put(0,0.1){\circle*{0.075}}
            \put(-0.1,0.25){\small $\gamma_{1}$}
            \put(1,0.1){\circle*{0.075}}
            \put(0.9,0.25){\small $\gamma_{3}$}
            \put(2,0.1){\circle*{0.075}}
            \put(1.9,0.25){\small $\gamma_{4}$}
            \put(3,0.1){\circle*{0.075}}
            \put(3.1,0.25){\small $\gamma_{5}$}
            \put(3,0.6){\circle*{0.075}}
            \put(3.1,0.65){\small $\gamma_{2}$}
            \put(4,0.1){\circle*{0.075}}
            \put(3.9,0.25){\small $\gamma_{6}$}
            \put(5,0.1){\circle*{0.075}}
            \put(4.9,0.25){\small $\gamma_{7}$}
            \put(6,0.1){\circle*{0.075}}
            \put(5.9,0.25){\small $\gamma_{8}$}
            \put(0,0.1){\line(1,0){6}}
            \put(3,0.1){\line(0,1){0.5}}
           \end{picture}}

\vspace*{0.05in}
\noindent             
is not admissible.  
For each fundamental position, we exhibit a divergent game sequence 
as a short sequence of legal 
node firings which can be repeated indefinitely. 
The fundamental position $\omega_{1} = (1,0,0,0,0,0,0,0)$ is the 
$k=0$ version of the position $(2k+1,0,-k,0,0,0,0,0)$.  
From any such 
position with $k \geq 0$, play the game to see that 
the following sequence of node firings is legal: 
$(\gamma_{1}$, $\gamma_{3}$, $\gamma_{4}$, $\gamma_{5}$, $\gamma_{2}$, 
$\gamma_{6}$, $\gamma_{5}$, $\gamma_{4}$, $\gamma_{3}$, $\gamma_{7}$, 
$\gamma_{6}$, $\gamma_{5}$, $\gamma_{2}$, $\gamma_{4}$, $\gamma_{5}$, 
$\gamma_{6}$, $\gamma_{7}$, $\gamma_{8}$, $\gamma_{7}$, $\gamma_{6}$, 
$\gamma_{5}$, $\gamma_{2}$, $\gamma_{4}$, $\gamma_{3}$, $\gamma_{5}$, 
$\gamma_{4}$, $\gamma_{6}$, $\gamma_{5}$, $\gamma_{2}$, $\gamma_{7}$, 
$\gamma_{6}$,
$\gamma_{5}$, $\gamma_{4}$, $\gamma_{3})$. 
This results in the position $(2(k+1)+1,0,-(k+1),0,0,0,0,0)$. 
By symmetry, we also obtain 
a divergent game sequence from the fundamental position 
$\omega_{8}$. 
From the fundamental position $\omega_{2} = (0,1,0,0,0,0,0,0)$, 
play the (legal) sequence 
$(\gamma_{2}$, $\gamma_{5}$, $\gamma_{4}$, $\gamma_{3}$, $\gamma_{6}$, 
$\gamma_{5}$, $\gamma_{2}$, $\gamma_{4}$, $\gamma_{5}$, $\gamma_{6}$, 
$\gamma_{7}$, $\gamma_{6}$, $\gamma_{5}$, $\gamma_{2}$, $\gamma_{4}$, 
$\gamma_{3}$, $\gamma_{5}$, $\gamma_{4}$, $\gamma_{6}$, $\gamma_{5}$, 
$\gamma_{2}$, $\gamma_{8}$, $\gamma_{7}$, $\gamma_{6}$, $\gamma_{5}$, 
$\gamma_{2}$, $\gamma_{4}$, $\gamma_{3}$, $\gamma_{5}$, $\gamma_{4}$, 
$\gamma_{6}$, $\gamma_{5}$, $\gamma_{2}$, $\gamma_{7}$, $\gamma_{6}$, 
$\gamma_{5}$, $\gamma_{4}$, $\gamma_{8}$, $\gamma_{7}$, $\gamma_{6}$, 
$\gamma_{5}$, $\gamma_{2}$, $\gamma_{1}$, $\gamma_{3}$, $\gamma_{4}$, 
$\gamma_{5})$  
to obtain the position $(0,3,0,0,-4,4,0,0)$.  This is the 
$k=0$ version of the position $(0,2k+3,0,0,-2k-4,2k+4,0,0)$.  
From any such 
position with $k \geq 0$, play the game to see that 
the following sequence of node firings is legal: 
$(\gamma_{2}$, $\gamma_{6}$, $\gamma_{5}$, $\gamma_{4}$, $\gamma_{3}$, 
$\gamma_{7}$, $\gamma_{6}$, $\gamma_{5}$, $\gamma_{2}$, $\gamma_{4}$, 
$\gamma_{5}$, $\gamma_{6}$, $\gamma_{7}$, $\gamma_{8}$, $\gamma_{7}$, 
$\gamma_{6}$, $\gamma_{5}$, $\gamma_{2}$, $\gamma_{4}$, $\gamma_{3}$, 
$\gamma_{5}$, $\gamma_{4}$, $\gamma_{6}$, $\gamma_{5}$, $\gamma_{2}$, 
$\gamma_{7}$, $\gamma_{6}$, $\gamma_{5}$, $\gamma_{4}$, $\gamma_{3}$, 
$\gamma_{1}$, $\gamma_{3}$, $\gamma_{4}$, $\gamma_{5})$. 
This results in the position $(0,2(k+1)+3,0,0,-2(k+1)-4,2(k+1)+4,0,0)$. 
The fundamental position $\omega_{3} = (0,0,1,0,0,0,0,0)$ is the 
$k=0$ version of the position $(-4k,0,2k+1,0,0,0,0,0)$.  
From any such 
position with $k \geq 0$, play the game to see that 
the following sequence of node firings is legal: 
$(\gamma_{3}$, $\gamma_{4}$, $\gamma_{5}$, $\gamma_{2}$, $\gamma_{6}$, 
$\gamma_{5}$, $\gamma_{4}$, $\gamma_{3}$, $\gamma_{7}$, $\gamma_{6}$, 
$\gamma_{5}$, $\gamma_{2}$, $\gamma_{4}$, $\gamma_{5}$, $\gamma_{6}$, 
$\gamma_{7}$, $\gamma_{8}$, $\gamma_{7}$, $\gamma_{6}$, $\gamma_{5}$, 
$\gamma_{2}$, $\gamma_{4}$, $\gamma_{3}$, $\gamma_{5}$, $\gamma_{4}$, 
$\gamma_{6}$, $\gamma_{5}$, $\gamma_{2}$, $\gamma_{7}$, $\gamma_{6}$, 
$\gamma_{5}$, $\gamma_{4}$, $\gamma_{3}$, $\gamma_{1})$. 
This results in the position $(-4(k+1),0,2(k+1)+1,0,0,0,0,0)$. 
By symmetry, we also obtain 
a divergent game sequence from the fundamental position 
$\omega_{7}$. 
From the fundamental position $\omega_{4} = (0,0,0,1,0,0,0,0)$, 
play the (legal) sequence 
$(\gamma_{4}$, $\gamma_{3}$, $\gamma_{5}$, $\gamma_{2}$, $\gamma_{4}$, 
$\gamma_{5}$, $\gamma_{6}$, $\gamma_{5}$, $\gamma_{2}$, $\gamma_{4}$, 
$\gamma_{3}$, $\gamma_{5}$, $\gamma_{4}$, $\gamma_{7}$, $\gamma_{6}$, 
$\gamma_{5}$, $\gamma_{2}$, $\gamma_{4}$, $\gamma_{3}$, $\gamma_{5}$, 
$\gamma_{4}$, $\gamma_{6}$, $\gamma_{5}$, $\gamma_{7}$, $\gamma_{6}$, 
$\gamma_{8}$, $\gamma_{7}$, $\gamma_{6}$, $\gamma_{5}$, $\gamma_{2}$, 
$\gamma_{4}$, $\gamma_{3}$, $\gamma_{5}$, $\gamma_{4}$, $\gamma_{6}$, 
$\gamma_{5}$, $\gamma_{2}$, $\gamma_{7}$, $\gamma_{6}$, $\gamma_{5}$, 
$\gamma_{4}$, $\gamma_{3}$, $\gamma_{8}$, $\gamma_{7}$, $\gamma_{6}$, 
$\gamma_{5}$, $\gamma_{4}$, $\gamma_{1}$, $\gamma_{3})$  
to obtain the position $(0,0,-6,5,0,0,0,0)$.  This is the 
$k=0$ version of the position $(-3k,0,-3k-6,3k+5,0,0,0,0)$.  
From any such 
position with $k \geq 0$, play the game to see that 
the following sequence of node firings is legal: 
$(\gamma_{4}$, $\gamma_{5}$, $\gamma_{2}$, $\gamma_{6}$, $\gamma_{5}$, 
$\gamma_{4}$, $\gamma_{3}$, $\gamma_{7}$, $\gamma_{6}$, $\gamma_{5}$, 
$\gamma_{2}$, $\gamma_{4}$, $\gamma_{5}$, $\gamma_{6}$, $\gamma_{7}$, 
$\gamma_{8}$, $\gamma_{7}$, $\gamma_{6}$, $\gamma_{5}$, $\gamma_{2}$, 
$\gamma_{4}$, $\gamma_{3}$, $\gamma_{5}$, $\gamma_{4}$, $\gamma_{6}$, 
$\gamma_{5}$, $\gamma_{2}$, $\gamma_{7}$, $\gamma_{6}$, $\gamma_{5}$, 
$\gamma_{4}$, $\gamma_{3}$, $\gamma_{1}$, $\gamma_{3})$. 
This results in the position $(-3(k+1),0,-3(k+1)-6,3(k+1)+5,0,0,0,0)$. 
By symmetry, we also obtain a divergent game sequence from the fundamental position 
$\omega_{6}$. 
From the fundamental position $\omega_{5} = (0,0,0,0,1,0,0,0)$, 
play the (legal) sequence 
$(\gamma_{5}$, $\gamma_{2}$, $\gamma_{4}$, $\gamma_{3}$, $\gamma_{5}$, 
$\gamma_{4}$, $\gamma_{6}$, $\gamma_{5}$, $\gamma_{2}$, $\gamma_{4}$, 
$\gamma_{3}$, $\gamma_{5}$, $\gamma_{4}$, $\gamma_{6}$, $\gamma_{5}$, 
$\gamma_{7}$, $\gamma_{6}$, $\gamma_{5}$, $\gamma_{2}$, $\gamma_{4}$, 
$\gamma_{3}$, $\gamma_{5}$, $\gamma_{4}$, $\gamma_{6}$, $\gamma_{5}$, 
$\gamma_{2}$, $\gamma_{7}$, $\gamma_{6}$, $\gamma_{5}$, $\gamma_{8}$, 
$\gamma_{7}$, $\gamma_{6}$, $\gamma_{5}$, $\gamma_{2}$, $\gamma_{4}$, 
$\gamma_{3}$, $\gamma_{5}$, $\gamma_{4}$, $\gamma_{6}$, $\gamma_{5}$, 
$\gamma_{2}$, $\gamma_{7}$, $\gamma_{6}$, $\gamma_{5}$, $\gamma_{4}$, 
$\gamma_{3}$, $\gamma_{8}$, $\gamma_{7}$, $\gamma_{6}$, $\gamma_{5}$, 
$\gamma_{2}$, $\gamma_{4}$, $\gamma_{5}$, $\gamma_{1}$, $\gamma_{3}$, 
$\gamma_{4})$  
to obtain the position $(0,0,0,-8,7,0,0,0)$.  This is the 
$k=0$ version of the position $(-4k,0,0,-4k-8,4k+7,0,0,0,0)$.  
From any such 
position with $k \geq 0$, play the game to see that 
the following sequence of node firings is legal: 
$(\gamma_{5}$, $\gamma_{2}$, $\gamma_{6}$, $\gamma_{5}$, $\gamma_{4}$, 
$\gamma_{3}$, $\gamma_{7}$, $\gamma_{6}$, $\gamma_{5}$, $\gamma_{2}$, 
$\gamma_{4}$, $\gamma_{5}$, $\gamma_{6}$, $\gamma_{7}$, $\gamma_{8}$, 
$\gamma_{7}$, $\gamma_{6}$, $\gamma_{5}$, $\gamma_{2}$, $\gamma_{4}$, 
$\gamma_{3}$, $\gamma_{5}$, $\gamma_{4}$, $\gamma_{6}$, $\gamma_{5}$, 
$\gamma_{2}$, $\gamma_{7}$, $\gamma_{6}$, $\gamma_{5}$, $\gamma_{4}$, 
$\gamma_{3}$, $\gamma_{1}$, $\gamma_{3}$, $\gamma_{4})$. 
This results in the position $(-4(k+1),0,0,-4(k+1)-8,4(k+1)+7,0,0,0,0)$. 

We finish the $\widetilde{\myE}$ family by showing why\\ 
\hspace*{2.5in} 
\parbox[c]{3.1in}{
\setlength{\unitlength}{0.5in}
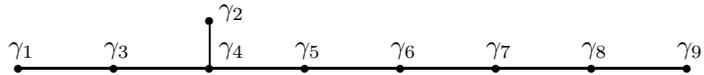
\begin{picture}(4.1,0.8)
            \put(0,0.1){\circle*{0.075}}
            \put(-0.1,0.25){\small $\gamma_{1}$}
            \put(1,0.1){\circle*{0.075}}
            \put(0.9,0.25){\small $\gamma_{3}$}
            \put(2,0.1){\circle*{0.075}}
            \put(2.1,0.25){\small $\gamma_{4}$}
            \put(2,0.6){\circle*{0.075}}
            \put(2.1,0.65){\small $\gamma_{2}$}
            \put(3,0.1){\circle*{0.075}}
            \put(2.9,0.25){\small $\gamma_{5}$}
            \put(4,0.1){\circle*{0.075}}
            \put(3.9,0.25){\small $\gamma_{6}$}
            \put(5,0.1){\circle*{0.075}}
            \put(4.9,0.25){\small $\gamma_{7}$}
            \put(6,0.1){\circle*{0.075}}
            \put(5.9,0.25){\small $\gamma_{8}$}
            \put(7,0.1){\circle*{0.075}}
            \put(6.9,0.25){\small $\gamma_{9}$}
            \put(0,0.1){\line(1,0){7}}
            \put(2,0.1){\line(0,1){0.5}}
           \end{picture}}

\vspace*{0.05in}
\noindent             
is not admissible. 
For each fundamental position, we exhibit a divergent game sequence 
as a short sequence of legal 
node firings which can be repeated indefinitely. 
The fundamental position $\omega_{1} = (1,0,0,0,0,0,0,0,0)$ is the 
$k=0$ version of the position $(2k+1,-k,0,0,0,0,0,0,0)$.  
From any such 
position with $k \geq 0$, play the game to see that 
the following sequence of node firings is legal: 
                            $(\gamma_{1}$,
                            $\gamma_{3}$,
                            $\gamma_{4}$,
                            $\gamma_{2}$,
                            $\gamma_{5}$,
                            $\gamma_{4}$,
                            $\gamma_{3}$,
                            $\gamma_{6}$,
                            $\gamma_{5}$,
                            $\gamma_{4}$,
                            $\gamma_{2}$,
                            $\gamma_{7}$,
                            $\gamma_{6}$,
                            $\gamma_{5}$,
                            $\gamma_{4}$,
                            $\gamma_{3}$,
                            $\gamma_{8}$,
                            $\gamma_{7}$,
                            $\gamma_{6}$,
                            $\gamma_{5}$,
                            $\gamma_{4}$,
                            $\gamma_{2}$,
                            $\gamma_{9}$,
                            $\gamma_{8}$,
                            $\gamma_{7}$,
                            $\gamma_{6}$,
                            $\gamma_{5}$,
                            $\gamma_{4}$,
                            $\gamma_{3})$. 
This results in the position $(2(k+1)+1,-(k+1),0,0,0,0,0,0,0)$. 
The fundamental position $\omega_{2} = (0,1,0,0,0,0,0,0,0)$ is the 
$k=0$ version of the position $(0,3k+1,0,-k,0,0,-k,0,0)$. 
From any such 
position with $k \geq 0$, play the game to see that 
the following sequence $\mathbf{s}$ 
of node firings is legal when played twice in a 
row: 
                            $\mathbf{s} := (\gamma_{2}$,
                            $\gamma_{4}$,
                            $\gamma_{3}$,
                            $\gamma_{1}$,
                            $\gamma_{5}$,
                            $\gamma_{4}$,
                            $\gamma_{3}$,
                            $\gamma_{6}$,
                            $\gamma_{5}$,
                            $\gamma_{4}$,
                            $\gamma_{7}$,
                            $\gamma_{6}$,
                            $\gamma_{5}$,
                            $\gamma_{8}$,
                            $\gamma_{7}$,
                            $\gamma_{6}$,
                            $\gamma_{9}$,
                            $\gamma_{8}$,
                            $\gamma_{7})$.
Playing $\mathbf{s}$ twice in a row results in the position 
$(0,3(k+1)+1,0,-(k+1),0,0,-(k+1),0,0)$. 
From the fundamental position $\omega_{3} = (0,0,1,0,0,0,0,0,0)$, 
play the (legal) sequence $(\gamma_{3}$,
                            $\gamma_{1}$,
                            $\gamma_{4}$,
                            $\gamma_{3}$,
                            $\gamma_{5}$,
                            $\gamma_{4}$,
                            $\gamma_{6}$,
                            $\gamma_{5}$,
                            $\gamma_{7}$,
                            $\gamma_{6}$,
                            $\gamma_{8}$,
                            $\gamma_{7}$,
                            $\gamma_{9}$,
                            $\gamma_{8})$  
to obtain the position $(0,2,0,0,0,0,0,-1,0)$.  This is the 
$k=0$ version of the position $(0,6k+2,0,-3k,0,0,-3k,-1,0)$.  
From any such 
position with $k \geq 0$, play the game to see that 
the following sequence $\mathbf{s}$ 
of node firings is legal when played three times in a 
row:  $\mathbf{s} := (\gamma_{2}$,
      $\gamma_{4}$,
      $\gamma_{3}$,
      $\gamma_{1}$,
      $\gamma_{5}$,
      $\gamma_{4}$,
      $\gamma_{3}$,
      $\gamma_{6}$,
      $\gamma_{5}$,
      $\gamma_{4}$,
      $\gamma_{7}$,
      $\gamma_{6}$,
      $\gamma_{5}$,
      $\gamma_{8}$,
      $\gamma_{7}$,
      $\gamma_{6}$,
      $\gamma_{9}$,
      $\gamma_{8}$,
      $\gamma_{7})$. 
Playing $\mathbf{s}$ three times in a row results in the position 
$(0,6(k+1)+2,0,-3(k+1),0,0,-3(k+1),-1,0)$.
The fundamental position $\omega_{4}$ is the $k=0$ version of the 
position $(0,-3k,0,2k+1,0,0,-k,0,0)$.  
From any such position with $k \geq 0$, play the game to see that the 
following sequence of node firings is legal:
                            $(\gamma_{4}$,
                            $\gamma_{3}$,
                            $\gamma_{1}$,
                            $\gamma_{5}$,
                            $\gamma_{4}$,
                            $\gamma_{3}$,
                            $\gamma_{6}$,
                            $\gamma_{5}$,
                            $\gamma_{4}$,
                            $\gamma_{7}$,
                            $\gamma_{6}$,
                            $\gamma_{5}$,
                            $\gamma_{8}$,
                            $\gamma_{7}$,
                            $\gamma_{6}$,
                            $\gamma_{9}$,
                            $\gamma_{8}$,
                            $\gamma_{7}$,
                            $\gamma_{2})$. 
This results in the position $(0,-3(k+1),0,2(k+1)+1,0,0,-(k+1),0,0)$. 
From the fundamental position $\omega_{5} = (0,0,0,0,1,0,0,0,0)$, 
play the (legal) sequence 
$(\gamma_{5}$,
$\gamma_{4}$,
$\gamma_{3}$,
$\gamma_{1}$,
$\gamma_{6}$,
$\gamma_{5}$,
$\gamma_{4}$,
$\gamma_{3}$,
$\gamma_{7}$,
$\gamma_{6}$,
$\gamma_{5}$,
$\gamma_{4}$,
$\gamma_{8}$,
$\gamma_{7}$,
$\gamma_{6}$,
$\gamma_{5}$,
$\gamma_{9}$,
$\gamma_{8}$,
$\gamma_{7}$,
$\gamma_{6})$ 
to obtain the position $(0,3,0,0,0,-1,0,0,0)$.  This is the 
$k=0$ version of the position $(0,15k+3,0,-5k,0,-1,-5k,0,0)$.  
From any such 
position with $k \geq 0$, play the game to see that 
the following sequence $\mathbf{s}$ 
of node firings is legal when played six times in a 
row:
     $\mathbf{s} := (\gamma_{2}$,     
     $\gamma_{4}$,
     $\gamma_{3}$,
     $\gamma_{1}$,
     $\gamma_{5}$,
     $\gamma_{4}$,
     $\gamma_{3}$,
     $\gamma_{6}$,
     $\gamma_{5}$,
     $\gamma_{4}$,
     $\gamma_{7}$,
     $\gamma_{6}$,
     $\gamma_{5}$,
     $\gamma_{8}$,
     $\gamma_{7}$,
     $\gamma_{6}$,
     $\gamma_{9}$,
     $\gamma_{8}$,
     $\gamma_{7})$. 
Playing $\mathbf{s}$ six times in a row results in the position 
$(0,15(k+1)+3,0,-5(k+1),0,-1,-5(k+1),0,0)$.
From the fundamental position $\omega_{6} = (0,0,0,0,0,1,0,0,0)$, 
play the (legal) sequence 
          $(\gamma_{6}$,
          $\gamma_{5}$,
          $\gamma_{4}$,
          $\gamma_{3}$,
          $\gamma_{1}$,
          $\gamma_{7}$,
          $\gamma_{6}$,
          $\gamma_{5}$,
          $\gamma_{4}$,
          $\gamma_{3}$,
          $\gamma_{8}$,
          $\gamma_{7}$,
          $\gamma_{6}$,
          $\gamma_{5}$,
          $\gamma_{4}$,
          $\gamma_{9}$,
          $\gamma_{8}$,
          $\gamma_{7}$,
          $\gamma_{6}$,
          $\gamma_{5})$ 
to obtain the position $(0,3,0,0,-1,0,0,0,0)$.  This is the 
$k=0$ version of the position $(0,12k+3,0,-4k,-1,0,-4k,0,0)$.  
From any such 
position with $k \geq 0$, play the game to see that 
the following sequence $\mathbf{s}$ 
of node firings is legal when played six times in a 
row:
     $\mathbf{s} := (\gamma_{2}$,     
      $\gamma_{4}$,
      $\gamma_{3}$,
      $\gamma_{1}$,
      $\gamma_{5}$,
      $\gamma_{4}$,
      $\gamma_{3}$,
      $\gamma_{6}$,
      $\gamma_{5}$,
      $\gamma_{4}$,
      $\gamma_{7}$,
      $\gamma_{6}$,
      $\gamma_{5}$,
      $\gamma_{8}$,
      $\gamma_{7}$,
      $\gamma_{6}$,
      $\gamma_{9}$,
      $\gamma_{8}$,
      $\gamma_{7})$. 
Playing $\mathbf{s}$ six times in a row results in the position 
$(0,12(k+1)+3,0,-4(k+1),-1,0,-4(k+1),0,0)$.
From the fundamental position $\omega_{7} = (0,0,0,0,0,0,1,0,0)$, 
play the (legal) sequence 
      $(\gamma_{7}$,
      $\gamma_{6}$,
      $\gamma_{5}$,
      $\gamma_{4}$,
      $\gamma_{3}$,
      $\gamma_{1}$,
      $\gamma_{8}$,
      $\gamma_{7}$,
      $\gamma_{6}$,
      $\gamma_{5}$,
      $\gamma_{4}$,
      $\gamma_{3}$,
      $\gamma_{9}$,
      $\gamma_{8}$,
      $\gamma_{7}$,
      $\gamma_{6}$,
      $\gamma_{5}$,
      $\gamma_{4})$
to obtain the position $(0,3,0,-1,0,0,0,0,0)$.  This is the 
$k=0$ version of the position $(0,3k+3,0,-k-1,0,0,-k,0,0)$.  
From any such 
position with $k \geq 0$, play the game to see that 
the following sequence $\mathbf{s}$ 
of node firings is legal when played twice in a 
row:
$\mathbf{s} := (\gamma_{2}$,     
      $\gamma_{4}$,
      $\gamma_{3}$,
      $\gamma_{1}$,
      $\gamma_{5}$,
      $\gamma_{4}$,
      $\gamma_{3}$,
      $\gamma_{6}$,
      $\gamma_{5}$,
      $\gamma_{4}$,
      $\gamma_{7}$,
      $\gamma_{6}$,
      $\gamma_{5}$,
      $\gamma_{8}$,
      $\gamma_{7}$,
      $\gamma_{6}$,
      $\gamma_{9}$,
      $\gamma_{8}$,
      $\gamma_{7})$. 
Playing $\mathbf{s}$ twice in a row results in the position 
$(0,3(k+1)+3,0,-(k+1)-1,0,0,-(k+1),0,0)$.
From the fundamental position $\omega_{8} = (0,0,0,0,0,0,0,1,0)$, 
play the (legal) sequence 
      $(\gamma_{8}$,
      $\gamma_{7}$,
      $\gamma_{6}$,
      $\gamma_{5}$,
      $\gamma_{4}$,
      $\gamma_{3}$,
      $\gamma_{1}$,
      $\gamma_{9}$,
      $\gamma_{8}$,
      $\gamma_{7}$,
      $\gamma_{6}$,
      $\gamma_{5}$,
      $\gamma_{4}$,
      $\gamma_{3}$,
      $\gamma_{2}$,
        $\gamma_{4}$,
        $\gamma_{3}$,
        $\gamma_{1}$,
        $\gamma_{5}$,
        $\gamma_{4}$,
        $\gamma_{3}$,
        $\gamma_{6}$,
        $\gamma_{5}$,
        $\gamma_{4}$,
        $\gamma_{7}$,
        $\gamma_{6}$,
        $\gamma_{5}$,
        $\gamma_{8}$,
        $\gamma_{7}$,
        $\gamma_{6}$,
        $\gamma_{9}$,
        $\gamma_{8}$,
        $\gamma_{7}$,        
        $\gamma_{2}$
        $\gamma_{4}$,
        $\gamma_{3}$,
        $\gamma_{1}$,
        $\gamma_{5}$,
        $\gamma_{4}$,
        $\gamma_{3}$,
        $\gamma_{6}$,
        $\gamma_{5}$,
        $\gamma_{4}$,
        $\gamma_{7}$,
        $\gamma_{6}$,
        $\gamma_{5}$,
        $\gamma_{8}$,
        $\gamma_{7}$,
        $\gamma_{6})$
to obtain the position $(0,4,0,-1,0,-1,0,0,0)$.  This is the 
$k=0$ version of the position $(0,6k+4,0,-2k-1,0,-1,-2k,0,0)$.  
From any such 
position with $k \geq 0$, play the game to see that 
the following sequence $\mathbf{s}$ 
of node firings is legal when played six times in a 
row:
$\mathbf{s} := (\gamma_{2}$,     
     $\gamma_{4}$,
     $\gamma_{3}$,
     $\gamma_{1}$,
     $\gamma_{5}$,
     $\gamma_{4}$,
     $\gamma_{3}$,
     $\gamma_{6}$,
     $\gamma_{5}$,
     $\gamma_{4}$,
     $\gamma_{7}$,
     $\gamma_{6}$,
     $\gamma_{5}$,
     $\gamma_{8}$,
     $\gamma_{7}$,
     $\gamma_{6}$,
     $\gamma_{9}$,
     $\gamma_{8}$,
     $\gamma_{7})$. 
Playing $\mathbf{s}$ six times in a row results in the position 
$(0,6(k+1)+4,0,-2(k+1)-1,0,-1,-2(k+1),0,0)$.
The fundamental position $\omega_{9} = (0,0,0,0,0,0,0,0,1)$ is the 
$k=0$ version of the position $(0,0,0,0,0,0,0,-k,2k+1)$.  
From any such 
position with $k \geq 0$, play the game to see that 
the following sequence of node firings is legal: 
                            $(\gamma_{9}$
                            $\gamma_{8}$,
                            $\gamma_{7}$,
                            $\gamma_{6}$,
                            $\gamma_{5}$,
                            $\gamma_{4}$,
                            $\gamma_{2}$,
                            $\gamma_{3}$,
                            $\gamma_{1}$,
                            $\gamma_{4}$,
                            $\gamma_{3}$,
                            $\gamma_{5}$,
                            $\gamma_{4}$,
                            $\gamma_{2}$,
                            $\gamma_{6}$,
                            $\gamma_{5}$,
                            $\gamma_{4}$,
                            $\gamma_{3}$,
                            $\gamma_{1}$,
                            $\gamma_{7}$,
                            $\gamma_{6}$,
                            $\gamma_{5}$,
                            $\gamma_{4}$,
                            $\gamma_{2}$,
                            $\gamma_{3}$,
                            $\gamma_{4}$,
                            $\gamma_{5}$,
                            $\gamma_{6}$,
                            $\gamma_{7}$,
                            $\gamma_{8}$,
                            $\gamma_{7}$,
                            $\gamma_{6}$,
                            $\gamma_{5}$,
                            $\gamma_{4}$,
                            $\gamma_{2}$,
                            $\gamma_{3}$,
                            $\gamma_{1}$,
                            $\gamma_{4}$,
                            $\gamma_{3}$,
                            $\gamma_{5}$,
                            $\gamma_{4}$,
                            $\gamma_{2}$,
                            $\gamma_{6}$,
                            $\gamma_{5}$,
                            $\gamma_{4}$,
                            $\gamma_{3}$,
                            $\gamma_{1}$,
                            $\gamma_{7}$,
                            $\gamma_{6}$,
                            $\gamma_{5}$,
                            $\gamma_{4}$,
                            $\gamma_{2}$,
                            $\gamma_{3}$,
                            $\gamma_{4}$,
                            $\gamma_{5}$,
                            $\gamma_{6}$,
                            $\gamma_{7}$,
                            $\gamma_{8})$. 
This results in the position $(0,0,0,0,0,0,0,-(k+1),2(k+1)+1)$.

\noindent 
\fbox{The $\widetilde{\myF}$ family}\ \ 
First, we show why 
\hspace*{0.05in}
\parbox[c]{2.1in}{
\setlength{\unitlength}{0.5in}
\begin{picture}(4.1,0.2)
            \put(0,0.1){\circle*{0.075}}
            \put(1,0.1){\circle*{0.075}}
            \put(2,0.1){\circle*{0.075}}
            \put(3,0.1){\circle*{0.075}}
            \put(4,0.1){\circle*{0.075}}
            \put(0,0.1){\line(1,0){4}}
            \put(1.2,0.1){\vector(1,0){0.1}}
            \put(1.3,0.1){\vector(1,0){0.1}}
            \put(1.8,0.1){\vector(-1,0){0.1}}
            \end{picture}
} is not admissible.  Label the nodes as $\gamma_{1}$, $\gamma_{2}$, 
$\gamma_{3}$, $\gamma_{4}$, and $\gamma_{5}$ from left to right.  
For each fundamental position, we exhibit a divergent game sequence 
as a short sequence of legal 
node firings which can be repeated indefinitely. 
The fundamental position $\omega_{1} = (1,0,0,0,0)$ is the 
$k=0$ version of the position $(2k+1,-k,0,0,0)$.  From any such 
position with $k \geq 0$, 
the following sequence of node firings is easily seen to be legal: 
$(\gamma_{1}$, $\gamma_{2}$, $\gamma_{3}$, $\gamma_{2}$, $\gamma_{4}$, 
$\gamma_{3}$, $\gamma_{2}$, $\gamma_{5}$, $\gamma_{4}$, $\gamma_{3}$, 
$\gamma_{2})$.  
This results in the position $(2(k+1)+1,-(k+1),0,0,0)$. 
The fundamental position $\omega_{2} = (0,1,0,0,0)$ is the 
$k=0$ version of the position $(-2k,4k+1,-2k,-2k,-2k)$.  From any such 
position with $k \geq 0$, 
the following sequence of node firings is easily seen to be legal: 
$(\gamma_{2}$, $\gamma_{3}$, $\gamma_{4}$, $\gamma_{5}$, $\gamma_{1}$, 
$\gamma_{2}$, $\gamma_{3}$, $\gamma_{4}$, $\gamma_{5}$, $\gamma_{1}$, 
$\gamma_{2}$, $\gamma_{3}$, $\gamma_{4}$, $\gamma_{5}$, $\gamma_{1})$.  
This results in the position $(-2(k+1),4(k+1)+1,-2(k+1),-2(k+1),-2(k+1))$. 
The fundamental position $\omega_{3} = (0,0,1,0,0)$ is the 
$k=0$ version of the position $(-k,-k,3k+1,-k,-k)$.  From any such 
position with $k \geq 0$, 
the following sequence of node firings is easily seen to be legal: 
$(\gamma_{3}$, $\gamma_{4}$, $\gamma_{5}$, $\gamma_{2}$, $\gamma_{1}$, 
$\gamma_{3}$, $\gamma_{4}$, $\gamma_{5}$, $\gamma_{2}$, $\gamma_{1}$, 
$\gamma_{3}$, $\gamma_{4}$, $\gamma_{5}$, $\gamma_{2}$, $\gamma_{1})$.  
This results in the position $(-(k+1),-(k+1),3(k+1)+1,-(k+1),-(k+1))$. 
The fundamental position $\omega_{4} = (0,0,0,1,0)$ is the 
$k=0$ version of the position $(-k,-k,k,2k+1,-k)$.  From any such 
position with $k \geq 0$, 
the following sequence of node firings is easily seen to be legal: 
$(\gamma_{4}$, $\gamma_{5}$, $\gamma_{3}$, $\gamma_{2}$, $\gamma_{1}$, 
$\gamma_{4}$, $\gamma_{5}$, $\gamma_{3}$, $\gamma_{2}$, $\gamma_{1}$, 
$\gamma_{4}$, $\gamma_{5}$, $\gamma_{3}$, $\gamma_{2}$, $\gamma_{1})$.  
This results in the position $(-(k+1),-(k+1),k+1,2(k+1)+1,-(k+1))$. 
The fundamental position $\omega_{5} = (0,0,0,0,1)$ is the 
$k=0$ version of the position $(0,0,0,-k,2k+1)$.  From any such 
position with $k \geq 0$, 
the following sequence of node firings is easily seen to be legal: 
$(\gamma_{5}$, $\gamma_{4}$, $\gamma_{3}$, $\gamma_{2}$, $\gamma_{3}$, 
$\gamma_{4}$, $\gamma_{1}$, $\gamma_{2}$, $\gamma_{3}$, $\gamma_{4}$, 
$\gamma_{2}$, $\gamma_{3}$, $\gamma_{1}$, $\gamma_{2}$, $\gamma_{3}$, $\gamma_{4})$.  
This results in the position $(0,0,0,-(k+1),2(k+1)+1)$. 

To finish our analysis of the $\widetilde{\myF}$ family,  we show why 
\hspace*{0.05in}
\parbox[c]{2.1in}{
\setlength{\unitlength}{0.5in}
\begin{picture}(4.1,0.2)
            \put(0,0.1){\circle*{0.075}}
            \put(1,0.1){\circle*{0.075}}
            \put(2,0.1){\circle*{0.075}}
            \put(3,0.1){\circle*{0.075}}
            \put(4,0.1){\circle*{0.075}}
            \put(0,0.1){\line(1,0){4}}
            \put(1.2,0.1){\vector(1,0){0.1}}
            \put(1.8,0.1){\vector(-1,0){0.1}}
            \put(1.7,0.1){\vector(-1,0){0.1}}
            \end{picture}
} is not admissible.  Label the nodes as $\gamma_{1}$, $\gamma_{2}$, 
$\gamma_{3}$, $\gamma_{4}$, and $\gamma_{5}$ from left to right.  
For each fundamental position, we exhibit a divergent game sequence 
as a short sequence of legal 
node firings which can be repeated indefinitely. 
The fundamental position $\omega_{1} = (1,0,0,0,0)$ is the 
$k=0$ version of the position $(2k+1,-k,0,0,0)$.  From any such 
position with $k \geq 0$, 
the following sequence of node firings is easily seen to be legal: 
$(\gamma_{1}$, $\gamma_{2}$, $\gamma_{3}$, $\gamma_{2}$, $\gamma_{4}$, 
$\gamma_{3}$, $\gamma_{2}$, $\gamma_{5}$, $\gamma_{4}$, $\gamma_{3}$, 
$\gamma_{2})$.  
This results in the position $(2(k+1)+1,-(k+1),0,0,0)$. 
The fundamental position $\omega_{2} = (0,1,0,0,0)$ is the 
$k=0$ version of the position $(-2k,4k+1,-k,-k,-k)$.  From any such 
position $k \geq 0$, 
the following sequence of node firings is easily seen to be legal: 
$(\gamma_{2}$, $\gamma_{3}$, $\gamma_{4}$, $\gamma_{5}$, $\gamma_{1}$, 
$\gamma_{2}$, $\gamma_{3}$, $\gamma_{4}$, $\gamma_{5}$, $\gamma_{1}$, 
$\gamma_{2}$, $\gamma_{3}$, $\gamma_{4}$, $\gamma_{5}$, $\gamma_{1})$.
This results in the position $(-2(k+1),4(k+1)+1,-(k+1),-(k+1),-(k+1))$.   
The fundamental position $\omega_{3} = (0,0,1,0,0)$ is the 
$k=0$ version of the position $(-2k,-2k,3k+1,-k,-k)$.  From any such 
position $k \geq 0$, 
the following sequence of node firings is easily seen to be legal: 
$(\gamma_{3}$, $\gamma_{4}$, $\gamma_{5}$, $\gamma_{2}$, $\gamma_{1}$, 
$\gamma_{3}$, $\gamma_{4}$, $\gamma_{5}$, $\gamma_{2}$, $\gamma_{1}$, 
$\gamma_{3}$, $\gamma_{4}$, $\gamma_{5}$, $\gamma_{2}$, $\gamma_{1})$.  
This results in the position $(-2(k+1),-2(k+1),3(k+1)+1,-(k+1),-(k+1))$. 
The fundamental position $\omega_{4} = (0,0,0,1,0)$ is the 
$k=0$ version of the position $(0,0,0,2k+1,-4k)$.  From any such 
position $k \geq 0$, 
the following sequence of node firings is easily seen to be legal: 
$(\gamma_{4}$, $\gamma_{3}$, $\gamma_{2}$, $\gamma_{3}$, $\gamma_{4}$, 
$\gamma_{1}$, $\gamma_{2}$, $\gamma_{3}$, $\gamma_{4}$, $\gamma_{2}$, 
$\gamma_{3}$, $\gamma_{1}$, $\gamma_{2}$, $\gamma_{3}$, $\gamma_{4}$, 
$\gamma_{5})$.  
This results in the position $(0,0,0,2(k+1)+1,-4(k+1))$. 
The fundamental position $\omega_{5} = (0,0,0,0,1)$ is the 
$k=0$ version of the position $(0,0,0,-k,2k+1)$.  From any such 
position $k \geq 0$, 
the following sequence of node firings is easily seen to be legal: 
$(\gamma_{5}$, $\gamma_{4}$, $\gamma_{3}$, $\gamma_{2}$, $\gamma_{3}$, 
$\gamma_{4}$, $\gamma_{1}$, $\gamma_{2}$, $\gamma_{3}$, $\gamma_{4}$, 
$\gamma_{2}$, $\gamma_{3}$, $\gamma_{1}$, $\gamma_{2}$, $\gamma_{3}$, $\gamma_{4})$.  
This results in the position $(0,0,0,-(k+1),2(k+1)+1)$.

\noindent 
\fbox{The $\widetilde{\myG}$ family}\ \ 
First, we show why 
\hspace*{0.05in}
\parbox[c]{1.1in}{
\setlength{\unitlength}{0.5in}
\begin{picture}(2.1,0.2)
            \put(0,0.1){\circle*{0.075}}
            \put(1,0.1){\circle*{0.075}}
            \put(2,0.1){\circle*{0.075}}
            \put(0,0.1){\line(1,0){2}}
            \put(0.2,0.1){\vector(1,0){0.1}}
            \put(0.8,0.1){\vector(-1,0){0.1}}
            \put(0.7,0.1){\vector(-1,0){0.1}}
            \put(0.6,0.1){\vector(-1,0){0.1}}
            \end{picture}
} is not admissible.  Label the nodes as $\gamma_{1}$, $\gamma_{2}$, 
and $\gamma_{3}$ from left to right.  A position $(a,b,c)$ meets 
condition ({\tt *}) if $a \leq 0$, $b \leq 0$, and  $a+2b+c > 
0$.  The following inequalities are immediate: (1) $c>0$, 
(2) $b+c>0$, (3) $a+3b+3c > 0$, (4) $a+2b+2c > 0$, (5) $2a+3b+3c > 0$, 
(6) $a+b+c > 0$, and (7) $2a+4b+3c > 0$. 
From (1) through (6) it now follows that all node firings of the sequence 
$\mathbf{s} := (\gamma_{3}, \gamma_{2}, \gamma_{1}, \gamma_{2}, \gamma_{1}, 
\gamma_{2})$ are legal: The left-hand side of each inequality is the 
number at the respective node of the sequence when that node is 
fired.  The resulting position $(a_{1},b_{1},c_{1})$ has $a_{1} = a$, 
$b_{1} = -(a+b+c)$, and $c_{1} = 2a+4b+3c$.  Clearly $a_{1} \leq 0$.  
By inequality (6), it follows that $b_{1} < 0$.  From inequality (7) 
we get $c_{1} > 0$.  Finally, $a_{1}+2b_{1}+c_{1} = a+2b+c > 0$, so 
$(a_{1},b_{1},c_{1})$ meets condition ({\tt *}).  
So from any position which meets condition ({\tt *}), the firing 
sequence $\mathbf{s}$ can be 
legally applied indefinitely, resulting in a divergent game sequence. 
Then it suffices to show that from each fundamental position we can 
reach a position which meets condition ({\tt *}) using a sequence of 
legal node firings. 
The 
fundamental position $\omega_{3} = (0,0,1)$ meets condition 
({\tt *}).  
Now take fundamental position $\omega_{2} = (0,1,0)$ 
and apply the legal firing sequence $(\gamma_{2}, \gamma_{1}, 
\gamma_{2}, \gamma_{1}, \gamma_{2})$ to get the resulting position 
$(0,-1,4)$.  The latter 
meets condition ({\tt *}).  
For the fundamental position $\omega_{1} = (1,0,0)$, 
apply the legal firing sequence $(\gamma_{1}, \gamma_{2}, 
\gamma_{1}, \gamma_{2}, \gamma_{1})$ to get the resulting position 
$(-1,0,2)$.  The latter 
meets condition ({\tt *}).  

Next, we show why 
\hspace*{0.05in}
\parbox[c]{1.1in}{
\setlength{\unitlength}{0.5in}
\begin{picture}(2.1,0.2)
            \put(0,0.1){\circle*{0.075}}
            \put(1,0.1){\circle*{0.075}}
            \put(2,0.1){\circle*{0.075}}
            \put(0,0.1){\line(1,0){2}}
            \put(0.2,0.1){\vector(1,0){0.1}}
            \put(0.8,0.1){\vector(-1,0){0.1}}
            \put(0.7,0.1){\vector(-1,0){0.1}}
            \put(0.6,0.1){\vector(-1,0){0.1}}
            \put(1.2,0.1){\vector(1,0){0.1}}
            \put(1.8,0.1){\vector(-1,0){0.1}}
            \put(1.7,0.1){\vector(-1,0){0.1}}
            \end{picture}
} is not admissible.  Our argument is similar to the previous case. 
Label the nodes as $\gamma_{1}$, $\gamma_{2}$, 
and $\gamma_{3}$ from left to right.  A position $(a,b,c)$ meets 
condition ({\tt *}) if $b \leq 0$, $c \leq 0$, $a+3b > 0$, and  $a+b+c > 
0$.  The following inequalities are easy to see: (1) $a>0$, 
(2) $a+b>0$, (3) $2a+3b > 0$, (4) $a+2b > 0$, (5) $a+3b > 0$, 
(6) $2a+3b+c > 0$, (7) $4a+7b+2c > 0$, (8) $2a+4b+c > 0$, and (9) 
$11a+18b+6c > 0$. 
From (1) through (8) it now follows that all node firings of the sequence 
$\mathbf{s} := 
(\gamma_{1}, \gamma_{2}, \gamma_{1}, \gamma_{2}, \gamma_{1}, 
\gamma_{3}, \gamma_{2}, \gamma_{3})$ 
are legal: The left-hand side of each inequality is the 
number at the respective node of the sequence when that node is 
fired.  The resulting position $(a_{1},b_{1},c_{1})$ has $a_{1} = 
11a+18b+c$, 
$b_{1} = b$, and $c_{1} = -(2a+4b+c)$.  Clearly $b_{1} \leq 0$.  
Inequality (8) gives $c_{1} < 0$. From (9) 
we get $a_{1} > 0$.  Note that $a_{1}+3b_{1} = 11a+21b+6c = 
6(a+b+c) + 5(a+3b) > 0$.  Finally, $a_{1}+b_{1}+c_{1} = 
9a+15b+5c = 5(a+b+c) + 2(a+2b) + 2(a+3b) > 0$, so 
$(a_{1},b_{1},c_{1})$ meets condition ({\tt *}).  
So from any position which meets condition ({\tt *}), the firing 
sequence $\mathbf{s}$ can be 
legally applied indefinitely, resulting in a divergent game sequence. 
Then it suffices to show that from each fundamental position we can 
reach a position which meets condition ({\tt *}) using a sequence of 
legal node firings. 
The 
fundamental position $\omega_{1} = (1,0,0)$ meets condition 
({\tt *}).  
For the fundamental position $\omega_{2} = (0,1,0)$,  
apply the legal firing sequence $(\gamma_{2}, \gamma_{3}, 
\gamma_{2})$ to get the resulting position 
$(6,-1,0)$.  The latter 
meets condition ({\tt *}).  
For the fundamental position $\omega_{3} = (0,0,1)$, 
apply the legal firing sequence $(\gamma_{3}, \gamma_{2}, 
\gamma_{3})$ to get the resulting position 
$(6,0,-1)$, which  
meets condition ({\tt *}).  

Next, we show why 
\hspace*{0.05in}
\parbox[c]{1.1in}{
\setlength{\unitlength}{0.5in}
\begin{picture}(2.1,0.2)
            \put(0,0.1){\circle*{0.075}}
            \put(1,0.1){\circle*{0.075}}
            \put(2,0.1){\circle*{0.075}}
            \put(0,0.1){\line(1,0){2}}
            \put(0.2,0.1){\vector(1,0){0.1}}
            \put(0.8,0.1){\vector(-1,0){0.1}}
            \put(0.7,0.1){\vector(-1,0){0.1}}
            \put(0.6,0.1){\vector(-1,0){0.1}}
            \put(1.2,0.1){\vector(1,0){0.1}}
            \put(1.8,0.1){\vector(-1,0){0.1}}
            \put(1.7,0.1){\vector(-1,0){0.1}}
            \put(1.6,0.1){\vector(-1,0){0.1}}
            \end{picture}
} is not admissible.  Our argument is entirely similar to the previous case. 
Label the nodes as $\gamma_{1}$, $\gamma_{2}$, 
and $\gamma_{3}$ from left to right.  A position $(a,b,c)$ meets 
condition ({\tt *}) if $b \leq 0$, $c \leq 0$, $a+3b > 0$, and  $a+b+c > 
0$.  The following inequalities are easy to see: (1) $a>0$, 
(2) $a+b>0$, (3) $2a+3b > 0$, (4) $a+2b > 0$, (5) $a+3b > 0$, 
(6) $2a+3b+c > 0$, (7) $6a+10b+3c > 0$, (8) $4a+7b+2c > 0$, (9) 
$6a+11b+3c > 0$, (10) $2a+4b+c > 0$, and (11) $35a+60b+18c > 0$. 
From (1) through (10) it now follows that all node firings of the sequence 
$\mathbf{s} := 
(\gamma_{1}, \gamma_{2}, \gamma_{1}, \gamma_{2}, \gamma_{1}, 
\gamma_{3}, \gamma_{2}, \gamma_{3}, \gamma_{2}, \gamma_{3})$ 
are legal: The left-hand side of each inequality is the 
number at the respective node of the sequence when that node is 
fired.  The resulting position $(a_{1},b_{1},c_{1})$ has $a_{1} = 
35a+60b+18c$, 
$b_{1} = b$, and $c_{1} = -(2a+4b+c)$.  Clearly $b_{1} \leq 0$.  
Inequality (10) gives $c_{1} < 0$. From (11) 
we get $a_{1} > 0$.  Note that $a_{1}+3b_{1} = 35a+63b+18c = 
18(a+b+c) + 15(a+3b) + 2a > 0$.  Finally, $a_{1}+b_{1}+c_{1} = 
33a+57b+17c = 17(a+b+c) + 8(a+2b) + 8(a+3b) > 0$, so 
$(a_{1},b_{1},c_{1})$ meets condition ({\tt *}).  
So from any position which meets condition ({\tt *}), the firing 
sequence $\mathbf{s}$ can be 
legally applied indefinitely, resulting in a divergent game sequence. 
Then it suffices to show that from each fundamental position we can 
reach a position which meets condition ({\tt *}) using a sequence of 
legal node firings. 
The 
fundamental position $\omega_{1} = (1,0,0)$ meets condition 
({\tt *}).  
For the fundamental position $\omega_{2} = (0,1,0)$,  
apply the legal firing sequence $(\gamma_{2}, \gamma_{3}, 
\gamma_{2}, \gamma_{3}, 
\gamma_{2})$ to get the resulting position 
$(12,-1,0)$.  The latter 
meets condition ({\tt *}).  
For the fundamental position $\omega_{3} = (0,0,1)$, 
apply the legal firing sequence $(\gamma_{3}, \gamma_{2}, 
\gamma_{3},\gamma_{2}, 
\gamma_{3})$ to get the resulting position 
$(18,0,-1)$, which  
meets condition ({\tt *}).

Next, we show why 
\hspace*{0.05in}
\parbox[c]{1.1in}{
\setlength{\unitlength}{0.5in}
\begin{picture}(2.1,0.2)
            \put(0,0.1){\circle*{0.075}}
            \put(1,0.1){\circle*{0.075}}
            \put(2,0.1){\circle*{0.075}}
            \put(0,0.1){\line(1,0){2}}
            \put(0.2,0.1){\vector(1,0){0.1}}
            \put(0.8,0.1){\vector(-1,0){0.1}}
            \put(0.7,0.1){\vector(-1,0){0.1}}
            \put(0.6,0.1){\vector(-1,0){0.1}}
            \put(1.8,0.1){\vector(-1,0){0.1}}
            \put(1.2,0.1){\vector(1,0){0.1}}
            \put(1.3,0.1){\vector(1,0){0.1}}
            \put(1.4,0.1){\vector(1,0){0.1}}
            \end{picture}
} is not admissible.  Our argument is entirely similar to the previous case. 
Label the nodes as $\gamma_{1}$, $\gamma_{2}$, 
and $\gamma_{3}$ from left to right.  A position $(a,b,c)$ meets 
condition ({\tt *}) if $b \leq 0$, $c \leq 0$, $a+3b > 0$, and  $3a+6b+c > 
0$.  The following inequalities are easy to see: (1) $a>0$, 
(2) $a+b>0$, (3) $2a+3b > 0$, (4) $a+2b > 0$, (5) $a+3b > 0$, 
(6) $6a+9b+c > 0$, (7) $6a+10b+c > 0$, (8) $12a+21b+2c > 0$, (9) 
$6a+11b+c > 0$, (10) $6a+12b+c > 0$, and (11) $35a+60b+6c > 0$. 
From (1) through (10) it now follows that all node firings of the sequence 
$\mathbf{s} := 
(\gamma_{1}, \gamma_{2}, \gamma_{1}, \gamma_{2}, \gamma_{1}, 
\gamma_{3}, \gamma_{2}, \gamma_{3}, \gamma_{2}, \gamma_{3})$ 
are legal: The left-hand side of each inequality is the 
number at the respective node of the sequence when that node is 
fired.  The resulting position $(a_{1},b_{1},c_{1})$ has $a_{1} = 
35a+60b+6c$, 
$b_{1} = b$, and $c_{1} = -(6a+12b+c)$.  Clearly $b_{1} \leq 0$.  
Inequality (10) gives $c_{1} < 0$. From (11) 
we get $a_{1} > 0$.  Note that $a_{1}+3b_{1} = 35a+63b+6c = 
6(3a+6b+c) + 9(a+3b) + 8a > 0$.  Finally, $3a_{1}+6b_{1}+c_{1} = 
99a+174b+17c = 17(3a+6b+c) + 24(2a+3b) > 0$, so 
$(a_{1},b_{1},c_{1})$ meets condition ({\tt *}).  
So from any position which meets condition ({\tt *}), the firing 
sequence $\mathbf{s}$ can be 
legally applied indefinitely, resulting in a divergent game sequence. 
Then it suffices to show that from each fundamental position we can 
reach a position which meets condition ({\tt *}) using a sequence of 
legal node firings. 
The 
fundamental position $\omega_{1} = (1,0,0)$ meets condition 
({\tt *}).  
For the fundamental position $\omega_{2} = (0,1,0)$,  
apply the legal firing sequence $(\gamma_{2}, \gamma_{3}, 
\gamma_{2}, \gamma_{3}, 
\gamma_{2})$ to get the resulting position 
$(12,-1,0)$, which  
meets condition ({\tt *}).  
For the fundamental position $\omega_{3} = (0,0,1)$, 
apply the legal firing sequence $(\gamma_{3}, \gamma_{2}, 
\gamma_{3},\gamma_{2}, 
\gamma_{3})$ to get the resulting position 
$(6,0,-1)$, which  
meets condition ({\tt *}).  

Next, we show why 
\hspace*{0.05in}
\parbox[c]{1.1in}{
\setlength{\unitlength}{0.5in}
\begin{picture}(2.1,0.2)
            \put(0,0.1){\circle*{0.075}}
            \put(1,0.1){\circle*{0.075}}
            \put(2,0.1){\circle*{0.075}}
            \put(0,0.1){\line(1,0){2}}
            \put(0.2,0.1){\vector(1,0){0.1}}
            \put(0.8,0.1){\vector(-1,0){0.1}}
            \put(0.7,0.1){\vector(-1,0){0.1}}
            \put(0.6,0.1){\vector(-1,0){0.1}}
            \put(1.8,0.1){\vector(-1,0){0.1}}
            \put(1.2,0.1){\vector(1,0){0.1}}
            \put(1.3,0.1){\vector(1,0){0.1}}
            \end{picture}
} is not admissible.  Our argument is entirely similar to the previous case. 
Label the nodes as $\gamma_{1}$, $\gamma_{2}$, 
and $\gamma_{3}$ from left to right.  A position $(a,b,c)$ meets 
condition ({\tt *}) if $b \leq 0$, $c \leq 0$, $a+3b > 0$, and  $2a+4b+c > 
0$.  The following inequalities are easy to see: (1) $a>0$, 
(2) $a+b>0$, (3) $2a+3b > 0$, (4) $a+2b > 0$, (5) $a+3b > 0$, 
(6) $4a+6b+c > 0$, (7) $4a+7b+c > 0$, (8) $4a+8b+c > 0$, 
and (9) $11a+18b+3c > 0$. 
From (1) through (8) it now follows that all node firings of the sequence 
$\mathbf{s} := 
(\gamma_{1}, \gamma_{2}, \gamma_{1}, \gamma_{2}, \gamma_{1}, 
\gamma_{3}, \gamma_{2}, \gamma_{3})$ 
are legal: The left-hand side of each inequality is the 
number at the respective node of the sequence when that node is 
fired.  The resulting position $(a_{1},b_{1},c_{1})$ has $a_{1} = 
11a+18b+3c$, 
$b_{1} = b$, and $c_{1} = -(4a+8b+c)$.  Clearly $b_{1} \leq 0$.  
Inequality (8) gives $c_{1} < 0$. From (9) 
we get $a_{1} > 0$.  Note that $a_{1}+3b_{1} = 11a+21b+3c = 
3(2a+4b+c) + 3(a+3b) + 2a > 0$.  Finally, $2a_{1}+4b_{1}+c_{1} = 
18a+32b+5c = 5(2a+4b+c) + 4(2a+3b) > 0$, so 
$(a_{1},b_{1},c_{1})$ meets condition ({\tt *}).  
So from any position which meets condition ({\tt *}), the firing 
sequence $\mathbf{s}$ can be 
legally applied indefinitely, resulting in a divergent game sequence. 
Then it suffices to show that from each fundamental position we can 
reach a position which meets condition ({\tt *}) using a sequence of 
legal node firings. 
The 
fundamental position $\omega_{1} = (1,0,0)$ meets condition 
({\tt *}).  
For the fundamental position $\omega_{2} = (0,1,0)$,  
apply the legal firing sequence $(\gamma_{2}, \gamma_{3}, 
\gamma_{2})$ to get the resulting position 
$(6,-1,0)$.  The latter 
meets condition ({\tt *}).  
For the fundamental position $\omega_{3} = (0,0,1)$, 
apply the legal firing sequence $(\gamma_{3}, \gamma_{2}, 
\gamma_{3})$ to get the resulting position 
$(3,0,-1)$, which  
meets condition ({\tt *}).  

To finish our analysis of the $\widetilde{\myG}$ family, 
we show why 
\hspace*{0.05in}
\parbox[c]{1.1in}{
\setlength{\unitlength}{0.5in}
\begin{picture}(2.1,0.2)
            \put(0,0.1){\circle*{0.075}}
            \put(1,0.1){\circle*{0.075}}
            \put(2,0.1){\circle*{0.075}}
            \put(0,0.1){\line(1,0){2}}
            \put(0.2,0.1){\vector(1,0){0.1}}
            \put(0.3,0.1){\vector(1,0){0.1}}
            \put(0.4,0.1){\vector(1,0){0.1}}
            \put(0.8,0.1){\vector(-1,0){0.1}}
            \end{picture}
} is not admissible.  Our argument is similar to the previous case. 
Label the nodes as $\gamma_{1}$, $\gamma_{2}$, 
and $\gamma_{3}$ from left to right.  A position $(a,b,c)$ meets 
condition ({\tt *}) if $a \leq 0$, $b \leq 0$, and  $3a+2b+c > 
0$.  The following inequalities are easy to see: (1) $c>0$, 
(2) $b+c>0$, (3) $a+b+c > 0$, (4) $3a+2b+2c > 0$, (5) $2a+b+c > 0$, 
(6) $3a+b+c > 0$, and (7) $6a+4b+3c > 0$. 
From (1) through (6) it now follows that all node firings of the sequence 
$\mathbf{s} := 
(\gamma_{3}, \gamma_{2}, \gamma_{1}, \gamma_{2}, \gamma_{1}, 
\gamma_{2})$ 
are legal: The left-hand side of each inequality is the 
number at the respective node of the sequence when that node is 
fired.  The resulting position $(a_{1},b_{1},c_{1})$ has $a_{1} = a$, 
$b_{1} = -(3a+b+c)$, and $c_{1} = 6a+4b+3c$.  Clearly $a_{1} \leq 0$.  
Inequality (6) gives $b_{1} < 0$. Note that $3a_{1}+2b_{1}+c_{1} = 
3a+2b+c > 0$, so 
$(a_{1},b_{1},c_{1})$ meets condition ({\tt *}).  
So from any position which meets condition ({\tt *}), the firing 
sequence $\mathbf{s}$ can be 
legally applied indefinitely, resulting in a divergent game sequence. 
Then it suffices to show that from each fundamental position we can 
reach a position which meets condition ({\tt *}) using a sequence of 
legal node firings. 
The 
fundamental position $\omega_{3} = (0,0,1)$ meets condition 
({\tt *}).  
For the fundamental position $\omega_{2} = (0,1,0)$,  
apply the legal firing sequence $(\gamma_{2}, \gamma_{1}, 
\gamma_{2}, \gamma_{1}, \gamma_{2})$ to get the resulting position 
$(0,-1,4)$, which  
meets condition ({\tt *}).  
For the fundamental position $\omega_{1} = (1,0,0)$, 
apply the legal firing sequence $(\gamma_{1}, \gamma_{2}, 
\gamma_{1}, \gamma_{2}, \gamma_{1})$ to get the resulting position 
$(-1,0,6)$, which  
meets condition ({\tt *}).

\noindent 
\fbox{Families of small cycles}\ \ 
First, we show why GCM graphs of the form \hspace*{0.1in}
\parbox[c]{0.5in}{
\setlength{\unitlength}{0.75in}
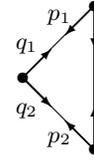
\begin{picture}(0.6,1.2)
            \put(0,0.6){\circle*{0.075}}
            \put(0.5,0.1){\circle*{0.075}}
            \put(0.5,1.1){\circle*{0.075}}
            \put(0,0.6){\line(1,1){0.5}}
            \put(0,0.6){\line(1,-1){0.5}}
            \put(0.5,0.1){\line(0,1){1}}
            \put(0.5,0.3){\vector(0,1){0.1}}
            \put(0.5,0.9){\vector(0,-1){0.1}}
            \put(0,0.6){\vector(1,1){0.2}}
            \put(-0.05,0.8){\footnotesize $q_{1}$}
            \put(0.175,1){\footnotesize $p_{1}$}
            \put(0.5,1.1){\vector(-1,-1){0.2}}
            \put(0,0.6){\vector(1,-1){0.2}}
            \put(0.5,0.1){\vector(-1,1){0.2}}
            \put(-0.05,0.35){\footnotesize $q_{2}$}
            \put(0.175,0.15){\footnotesize $p_{2}$} 
\end{picture}
}
are not admissible.  
Assign numbers $a$, $b$, 
and $c$ as follows:  
\hspace*{0.1in}
\parbox[c]{0.5in}{
\setlength{\unitlength}{0.75in}
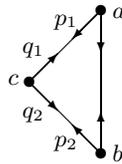
\begin{picture}(0.6,1.2)
            \put(0,0.6){\circle*{0.075}}
            \put(0.5,0.1){\circle*{0.075}}
            \put(0.5,1.1){\circle*{0.075}}
            \put(0,0.6){\line(1,1){0.5}}
            \put(0,0.6){\line(1,-1){0.5}}
            \put(0.5,0.1){\line(0,1){1}}
            \put(0.5,0.3){\vector(0,1){0.1}}
            \put(0.5,0.9){\vector(0,-1){0.1}}
            \put(0,0.6){\vector(1,1){0.2}}
            \put(-0.05,0.8){\footnotesize $q_{1}$}
            \put(0.175,1){\footnotesize $p_{1}$}
            \put(0.5,1.1){\vector(-1,-1){0.2}}
            \put(0,0.6){\vector(1,-1){0.2}}
            \put(0.5,0.1){\vector(-1,1){0.2}}
            \put(-0.05,0.35){\footnotesize $q_{2}$}
            \put(0.175,0.15){\footnotesize $p_{2}$} 
            \put(-0.15,0.58){\footnotesize $c$}
            \put(0.58,0.05){\footnotesize $b$}
            \put(0.58,1.05){\footnotesize $a$}
\end{picture}
}
Set $\kappa := (p_{1}+p_{2}-\frac{1}{q_{2}})a + 
(p_{1}+p_{2}-\frac{1}{q_{1}})b + c$.  
Assume for now that $a \geq 0$, $b \geq 0$, $c \leq 0$, and $\kappa > 0$; 
when these inequalities hold we will say the position $(a,b,c)$ meets 
condition ({\tt *}). Under 
condition ({\tt *}) notice that $a$ and $b$ cannot both be zero. 
Begin by firing only at the two 
rightmost nodes.  When this is no longer possible, fire at the 
leftmost node.  
The resulting corresponding 
numbers are $a_{1} = 
q_{1}(\kappa + \frac{1}{q_{2}}a)$, $b_{1} = q_{2}(\kappa + 
\frac{1}{q_{1}}b)$, and $c_{1} = -\kappa-\frac{1}{q_{2}}a-\frac{1}{q_{1}}b$. In 
particular, $a_{1} > 0$, 
$b_{1} > 0$, and $c_{1} < 0$.  
Next we 
check that $\kappa_{1} := (p_{1}+p_{2}-\frac{1}{q_{2}})a_{1} + 
(p_{1}+p_{2}-\frac{1}{q_{1}})b_{1} + c_{1}$ is also positive.  
Now 
\[\kappa_{1} = Q\kappa + Q_{1}a + Q_{2}b,\]
\noindent 
where $Q = q_{1}(p_{2} - \frac{1}{q_{2}}) + q_{2}(p_{1} - 
\frac{1}{q_{1}}) + (p_{1}q_{1} + p_{2}q_{2} - 1)$, 
$Q_{1} = \frac{1}{q_{2}}[q_{1}(p_{2}-\frac{1}{q_{2}}) + 
(p_{1}q_{1} - 1)]$, and 
$Q_{2} = \frac{1}{q_{1}}[q_{2}(p_{1}-\frac{1}{q_{1}}) + 
(p_{2}q_{2} - 1)]$.  
Since each 
parenthesized quantity in our expression for $Q$ is nonnegative and 
the last of these is  
positive, then $Q > 0$.  
Similar reasoning shows that each bracketed quantity in our 
expressions for $Q_{1}$ and $Q_{2}$ is nonnegative, hence $Q_{1} \geq 0$ 
and $Q_{2} \geq 0$. 
Since $\kappa > 0$ by hypothesis, it now follows that 
$\kappa_{1} > 0$. 
Then $(a_{1},b_{1},c_{1})$ meets condition ({\tt *}), so we 
can legally 
repeat the above firing sequence from position $(a_{1},b_{1},c_{1})$ 
to obtain another position 
$(a_{2},b_{2},c_{2})$ that meets condition 
({\tt *}), etc.  
Since the fundamental 
positions $(a,b,c) = (1,0,0)$ and $(a,b,c) = (0,1,0)$ meet  
condition ({\tt *}), then we see that the indicated 
legal firing sequence can be repeated 
indefinitely from these positions.  
For the fundamental 
position $(a,b,c) = (0,0,1)$, begin by firing 
at the leftmost node to 
obtain the position $(q_{1},q_{2},-1)$. This latter position meets  
condition ({\tt *}) with $\kappa = Q$, 
and so the legal firing sequence indicated above can 
be repeated indefinitely from this position.

Next, we show why GCM graphs of the form \hspace*{0.1in}
\parbox[c]{0.5in}{
\setlength{\unitlength}{0.75in}
\begin{picture}(0.6,1.2)
            \put(0,0.6){\circle*{0.075}}
            \put(0.5,0.1){\circle*{0.075}}
            \put(0.5,1.1){\circle*{0.075}}
            \put(0,0.6){\line(1,1){0.5}}
            \put(0,0.6){\line(1,-1){0.5}}
            \put(0.5,0.1){\line(0,1){1}}
            \put(0.5,0.3){\vector(0,1){0.1}}
            \put(0.5,0.9){\vector(0,-1){0.1}}
            \put(0.5,0.8){\vector(0,-1){0.1}}
            \put(0,0.6){\vector(1,1){0.2}}
            \put(-0.05,0.8){\footnotesize $q_{1}$}
            \put(0.175,1){\footnotesize $p_{1}$}
            \put(0.5,1.1){\vector(-1,-1){0.2}}
            \put(0,0.6){\vector(1,-1){0.2}}
            \put(0.5,0.1){\vector(-1,1){0.2}}
            \put(-0.05,0.35){\footnotesize $q_{2}$}
            \put(0.175,0.15){\footnotesize $p_{2}$} 
\end{picture}
}
are not admissible.  We assume that the amplitude products 
$p_{1}q_{1}$ and $p_{2}q_{2}$ are at least two. 
The argument is entirely similar to the previous case.  
Assign numbers $a$, $b$, 
and $c$ as follows:  
\hspace*{0.1in}
\parbox[c]{0.5in}{
\setlength{\unitlength}{0.75in}
\begin{picture}(0.6,1.2)
            \put(0,0.6){\circle*{0.075}}
            \put(0.5,0.1){\circle*{0.075}}
            \put(0.5,1.1){\circle*{0.075}}
            \put(0,0.6){\line(1,1){0.5}}
            \put(0,0.6){\line(1,-1){0.5}}
            \put(0.5,0.1){\line(0,1){1}}
            \put(0.5,0.3){\vector(0,1){0.1}}
            \put(0.5,0.9){\vector(0,-1){0.1}}
            \put(0.5,0.8){\vector(0,-1){0.1}}
            \put(0,0.6){\vector(1,1){0.2}}
            \put(-0.05,0.8){\footnotesize $q_{1}$}
            \put(0.175,1){\footnotesize $p_{1}$}
            \put(0.5,1.1){\vector(-1,-1){0.2}}
            \put(0,0.6){\vector(1,-1){0.2}}
            \put(0.5,0.1){\vector(-1,1){0.2}}
            \put(-0.05,0.35){\footnotesize $q_{2}$}
            \put(0.175,0.15){\footnotesize $p_{2}$} 
            \put(-0.15,0.58){\footnotesize $c$}
            \put(0.58,0.05){\footnotesize $b$}
            \put(0.58,1.05){\footnotesize $a$}
\end{picture}
}
Set $\kappa := (2p_{1}+2p_{2}-\frac{1}{q_{1}})a + 
(p_{1}+2p_{2}-\frac{1}{q_{2}})b + c$.  
Assume for now that $a \geq 0$, $b \geq 0$, $c \leq 0$, and $\kappa > 0$; 
when these inequalities hold we will say the position $(a,b,c)$ meets 
condition ({\tt *}). 
Using the same firing sequence as before, the resulting corresponding 
numbers are $a_{1} = 
q_{1}(\kappa + \frac{1}{q_{2}}b)$, $b_{1} = q_{2}(\kappa + 
\frac{1}{q_{1}}a)$, and $c_{1} = -\kappa-\frac{1}{q_{1}}a-\frac{1}{q_{2}}b$. In 
particular, $a_{1} > 0$, 
$b_{1} > 0$, and $c_{1} < 0$.  Next we 
check that $\kappa_{1} := (2p_{1}+2p_{2}-\frac{1}{q_{1}})a_{1} + 
(p_{1}+2p_{2}-\frac{1}{q_{2}})b_{1} + c_{1}$ is also positive.  
Now 
\[\kappa_{1} = Q\kappa + Q_{1}a + Q_{2}b,\]
\noindent 
where $Q = q_{1}(2p_{2}-\frac{1}{q_{1}}) + 
q_{2}(p_{1}-\frac{1}{q_{2}}) + (2p_{1}q_{1}+2p_{2}q_{2}-1)$, 
$Q_{1} = 
\frac{1}{q_{1}}[q_{2}(p_{1}-\frac{1}{q_{2}}) + 
(2p_{2}q_{2}-1)]$, 
and 
$Q_{2} = 
\frac{1}{q_{2}}[q_{1}(2p_{2}-\frac{1}{q_{1}}) + 
(2p_{1}q_{1}-1)]$. 
Since each 
parenthesized quantity in our expression for $Q$ is nonnegative and 
the last of these is  
positive, then $Q > 0$.  
Similar reasoning shows that each bracketed quantity in our 
expressions for $Q_{1}$ and $Q_{2}$ is nonnegative, hence $Q_{1} \geq 0$ 
and $Q_{2} \geq 0$. 
Since $\kappa > 0$ by hypothesis, it now follows that 
$\kappa_{1} > 0$. Conclude as in the previous case.

Next, we show why GCM graphs of the form \hspace*{0.1in}
\parbox[c]{0.5in}{
\setlength{\unitlength}{0.75in}
\begin{picture}(0.6,1.2)
            \put(0,0.6){\circle*{0.075}}
            \put(0.5,0.1){\circle*{0.075}}
            \put(0.5,1.1){\circle*{0.075}}
            \put(0,0.6){\line(1,1){0.5}}
            \put(0,0.6){\line(1,-1){0.5}}
            \put(0.5,0.1){\line(0,1){1}}
            \put(0.5,0.3){\vector(0,1){0.1}}
            \put(0.5,0.9){\vector(0,-1){0.1}}
            \put(0.5,0.8){\vector(0,-1){0.1}}
            \put(0.5,0.7){\vector(0,-1){0.1}}
            \put(0,0.6){\vector(1,1){0.2}}
            \put(0.5,1.1){\vector(-1,-1){0.2}}
            \put(0,0.6){\vector(1,-1){0.2}}
            \put(0.5,0.1){\vector(-1,1){0.2}}
            \put(-0.05,0.8){\footnotesize $q_{1}$}
            \put(0.175,1){\footnotesize $p_{1}$}
            \put(-0.05,0.35){\footnotesize $q_{2}$}
            \put(0.175,0.15){\footnotesize $p_{2}$} 
           \end{picture}
}
are not admissible. 
We assume that the amplitude products 
$p_{1}q_{1}$ and $p_{2}q_{2}$ are at least three. 
The argument is entirely similar to the previous two cases.  
Assign numbers $a$, $b$, 
and $c$ as follows:  
\hspace*{0.1in}
\parbox[c]{0.5in}{
\setlength{\unitlength}{0.75in}
\begin{picture}(0.6,1.2)
            \put(0,0.6){\circle*{0.075}}
            \put(0.5,0.1){\circle*{0.075}}
            \put(0.5,1.1){\circle*{0.075}}
            \put(0,0.6){\line(1,1){0.5}}
            \put(0,0.6){\line(1,-1){0.5}}
            \put(0.5,0.1){\line(0,1){1}}
            \put(0.5,0.3){\vector(0,1){0.1}}
            \put(0.5,0.9){\vector(0,-1){0.1}}
            \put(0.5,0.8){\vector(0,-1){0.1}}
            \put(0.5,0.7){\vector(0,-1){0.1}}
            \put(0,0.6){\vector(1,1){0.2}}
            \put(0.5,1.1){\vector(-1,-1){0.2}}
            \put(0,0.6){\vector(1,-1){0.2}}
            \put(0.5,0.1){\vector(-1,1){0.2}}
            \put(-0.05,0.8){\footnotesize $q_{1}$}
            \put(0.175,1){\footnotesize $p_{1}$}
            \put(-0.05,0.35){\footnotesize $q_{2}$}
            \put(0.175,0.15){\footnotesize $p_{2}$} 
            \put(-0.15,0.58){\footnotesize $c$}
            \put(0.58,0.05){\footnotesize $b$}
            \put(0.58,1.05){\footnotesize $a$}
           \end{picture}
}
Set $\kappa := (4p_{1}+6p_{2}-\frac{1}{q_{1}})a + 
(2p_{1}+4p_{2}-\frac{1}{q_{2}})b + c$.  
Assume for now that $a \geq 0$, $b \geq 0$, $c \leq 0$, and $\kappa > 0$; 
when these inequalities hold we will say the position $(a,b,c)$ meets 
condition ({\tt *}). 
Using the same firing sequence as in the previous two cases, 
the resulting corresponding 
numbers are $a_{1} = 
q_{1}(\kappa + \frac{1}{q_{2}}b)$, $b_{1} = q_{2}(\kappa + 
\frac{1}{q_{1}}a)$, and $c_{1} = -\kappa-\frac{1}{q_{1}}a-\frac{1}{q_{2}}b$. In 
particular, $a_{1} > 0$, 
$b_{1} > 0$, and $c_{1} < 0$.  Next we 
check that $\kappa_{1} := (4p_{1}+6p_{2}-\frac{1}{q_{1}})a_{1} + 
(2p_{1}+4p_{2}-\frac{1}{q_{2}})b_{1} + c_{1}$ is also positive.  
Now 
\[\kappa_{1} = Q\kappa + Q_{1}a + Q_{2}b,\]
\noindent 
where $Q = q_{1}(6p_{2}-\frac{1}{q_{1}}) + 
q_{2}(2p_{1}-\frac{1}{q_{2}}) + (4p_{1}q_{1}+4p_{2}q_{2}-1)$, 
$Q_{1} = 
\frac{1}{q_{1}}[q_{2}(2p_{1}-\frac{1}{q_{2}}) + 
(4p_{2}q_{2}-1)]$, 
and 
$Q_{2} = 
\frac{1}{q_{2}}[q_{1}(6p_{2}-\frac{1}{q_{1}}) + 
(4p_{1}q_{1}-1)]$. 
Since each 
parenthesized quantity in our expression for $Q$ is nonnegative and 
the last of these is  
positive, then $Q > 0$.  
Similar reasoning shows that each bracketed quantity in our 
expressions for $Q_{1}$ and $Q_{2}$ is nonnegative, hence $Q_{1} \geq 0$ 
and $Q_{2} \geq 0$. 
Since $\kappa > 0$ by hypothesis, it now follows that 
$\kappa_{1} > 0$. Conclude as in the previous two cases.

Next, we show why 
\hspace*{0.05in}
\parbox[c]{0.82in}{
\setlength{\unitlength}{0.7in}
\begin{picture}(1.1,1.2)
            \put(0,0.6){\circle*{0.075}}
            \put(0.5,0.1){\circle*{0.075}}
            \put(0.5,1.1){\circle*{0.075}}
            \put(1,0.6){\circle*{0.075}}
            \put(0,0.6){\line(1,1){0.5}}
            \put(0,0.6){\line(1,-1){0.5}}
            \put(1,0.6){\line(-1,1){0.5}}
            \put(1,0.6){\line(-1,-1){0.5}}
            \put(0.95,0.65){\vector(-1,1){0.1}}
            \put(0.55,1.05){\vector(1,-1){0.1}}
            \put(0.65,0.95){\vector(1,-1){0.1}}
           \end{picture}
} 
is not admissible.  The argument is similar to the previous three 
cases, but simpler since the amplitudes are all known.  Number the 
nodes with $\gamma_{1}$ as the North vertex, $\gamma_{2}$ as the East 
vertex, $\gamma_{3}$ as the South vertex, and $\gamma_{4}$ as the West 
vertex. 
We say an initial position $(a,b,c,d)$ meets 
condition ({\tt *}) if the following 
inequalities are satisfied: $b \geq 0$, $c \geq 0$, $d \leq 0$, and 
$a+d > 0$.  The firing sequence 
$(\gamma_{1},\gamma_{2},\gamma_{3},\gamma_{4})$ is easily seen to be 
legal from any such position.  The resulting position is 
$(a_{1},b_{1},c_{1},d_{1})$ with $a_{1} = 4a + 2b + c + d$, $b_{1} = 
c$, $c_{1} = a+d$, and $d_{1} = -(3a+b+c+d)$.  It is easy now to check 
that $(a_{1},b_{1},c_{1},d_{1})$ also meets condition ({\tt *}).  (In fact, 
the inequalities $c_{1} > 0$ and $d_{1} < 0$ are now strict.)  
So from any position which meets condition ({\tt *}), the firing 
sequence  $(\gamma_{1},\gamma_{2},\gamma_{3},\gamma_{4})$ can be 
legally applied indefinitely, resulting in a divergent game sequence. 
Then it suffices to show that from each fundamental position we can 
reach a position which meets condition ({\tt *}) using a sequence of 
legal node firings. 
The 
fundamental position $\omega_{1} = (1,0,0,0)$ meets condition 
({\tt *}).  
Now take fundamental position $\omega_{2} = (0,1,0,0)$ 
and apply the legal firing sequence $(\gamma_{2}, \gamma_{3}, 
\gamma_{4})$ to get the resulting position $(2,0,0,-1)$.  The latter 
meets condition ({\tt *}).  
Next take fundamental position $\omega_{3} = (0,0,1,0)$ 
and apply the legal firing sequence $(\gamma_{3}, \gamma_{4}, 
\gamma_{1}, \gamma_{2}, \gamma_{3}, 
\gamma_{4})$ to get the resulting position $(5,0,0,-3)$.  The latter 
meets condition ({\tt *}).  
For fundamental position $\omega_{4} = (0,0,0,1)$, 
apply the legal firing sequence $(\gamma_{4}, 
\gamma_{1}, \gamma_{2}, \gamma_{3}, 
\gamma_{4})$ to get the resulting position $(4,1,0,-3)$.  The latter 
meets condition ({\tt *}).

Next, we show why  
\hspace*{0.05in}
\parbox[c]{0.82in}{
\setlength{\unitlength}{0.7in}
\begin{picture}(1.1,1.2)
            \put(0,0.6){\circle*{0.075}}
            \put(0.5,0.1){\circle*{0.075}}
            \put(0.5,1.1){\circle*{0.075}}
            \put(1,0.6){\circle*{0.075}}
            \put(0,0.6){\line(1,1){0.5}}
            \put(0,0.6){\line(1,-1){0.5}}
            \put(1,0.6){\line(-1,1){0.5}}
            \put(1,0.6){\line(-1,-1){0.5}}
            \put(0.95,0.65){\vector(-1,1){0.1}}
            \put(0.55,1.05){\vector(1,-1){0.1}}
            \put(0.65,0.95){\vector(1,-1){0.1}}
            \put(0.45,0.15){\vector(-1,1){0.1}}
            \put(0.05,0.55){\vector(1,-1){0.1}}
            \put(0.15,0.45){\vector(1,-1){0.1}}
           \end{picture}
} 
is not admissible.  The argument is entirely similar to the previous 
case.  Number the 
nodes with $\gamma_{1}$ as the North vertex, $\gamma_{2}$ as the East 
vertex, $\gamma_{3}$ as the South vertex, and $\gamma_{4}$ as the West 
vertex. 
We say an initial position $(a,b,c,d)$ meets 
condition ({\tt *}) if the following 
inequalities are satisfied: $a > 0$, $b \geq 0$, $c \geq 0$, $d \leq 0$, and 
$3a+b+c+2d > 0$.  The firing sequence 
$(\gamma_{1},\gamma_{2},\gamma_{3},\gamma_{4})$ is easily seen to be 
legal from any such position.  The resulting position is 
$(a_{1},b_{1},c_{1},d_{1})$ with $a_{1} = 4a + 2b + c + d$, $b_{1} = 
c$, $c_{1} = 4a+b+c+2d$, and $d_{1} = -(3a+b+c+d)$.  It is easy now to check 
that $(a_{1},b_{1},c_{1},d_{1})$ also meets condition ({\tt *}).  (In fact, 
the inequalities $c_{1} > 0$ and $d_{1} < 0$ are now strict.)  
So from any position which meets condition ({\tt *}), the firing 
sequence  $(\gamma_{1},\gamma_{2},\gamma_{3},\gamma_{4})$ can be 
legally applied indefinitely, resulting in a divergent game sequence. 
Then it suffices to show that from each fundamental position we can 
reach a position which meets condition ({\tt *}) using a sequence of 
legal node firings. 
The 
fundamental position $\omega_{1} = (1,0,0,0)$ meets condition 
({\tt *}).  
Now take fundamental position $\omega_{2} = (0,1,0,0)$ 
and apply the legal firing sequence $(\gamma_{2}, \gamma_{3}, 
\gamma_{4})$ to get the resulting position $(2,0,1,-1)$.  The latter 
meets condition ({\tt *}).  
Next take fundamental position $\omega_{3} = (0,0,1,0)$ 
and apply the legal firing sequence $(\gamma_{3}, 
\gamma_{4})$ to get the resulting position $(1,1,1,-1)$.  The latter 
meets condition ({\tt *}).  
For fundamental position $\omega_{4} = (0,0,0,1)$, 
apply the legal firing sequence $(\gamma_{4})$ to get the resulting 
position $(1,0,2,-1)$.  The latter 
meets condition ({\tt *}).

Next, we show why  
\hspace*{0.05in}
\parbox[c]{0.82in}{
\setlength{\unitlength}{0.7in}
\begin{picture}(1.1,1.2)
            \put(0,0.6){\circle*{0.075}}
            \put(0.5,0.1){\circle*{0.075}}
            \put(0.5,1.1){\circle*{0.075}}
            \put(1,0.6){\circle*{0.075}}
            \put(0,0.6){\line(1,1){0.5}}
            \put(0,0.6){\line(1,-1){0.5}}
            \put(1,0.6){\line(-1,1){0.5}}
            \put(1,0.6){\line(-1,-1){0.5}}
            \put(0.95,0.65){\vector(-1,1){0.1}}
            \put(0.55,1.05){\vector(1,-1){0.1}}
            \put(0.65,0.95){\vector(1,-1){0.1}}
            \put(0.45,0.15){\vector(-1,1){0.1}}
            \put(0.05,0.55){\vector(1,-1){0.1}}
            \put(0.35,0.25){\vector(-1,1){0.1}}
           \end{picture}
} 
is not admissible.  The argument is entirely similar to the previous 
case.  Number the 
nodes with $\gamma_{1}$ as the North vertex, $\gamma_{2}$ as the East 
vertex, $\gamma_{3}$ as the South vertex, and $\gamma_{4}$ as the West 
vertex. 
We say an initial position $(a,b,c,d)$ meets 
condition ({\tt *}) if the following 
inequalities are satisfied: $a > 0$, $b \geq 0$, $c \geq 0$, $d \leq 0$, and 
$3a+b+c+d > 0$.  The firing sequence 
$(\gamma_{1},\gamma_{2},\gamma_{3},\gamma_{4})$ is easily seen to be 
legal from any such position.  The resulting position is 
$(a_{1},b_{1},c_{1},d_{1})$ with $a_{1} = 6a + 3b + 2c + d$, $b_{1} = 
c$, $c_{1} = 3a+b+c+d$, and $d_{1} = -(5a+2b+2c+d)$.  It is easy now to check 
that $(a_{1},b_{1},c_{1},d_{1})$ also meets condition ({\tt *}).  (In fact, 
the inequalities $c_{1} > 0$ and $d_{1} < 0$ are now strict.)  
So from any position which meets condition ({\tt *}), the firing 
sequence  $(\gamma_{1},\gamma_{2},\gamma_{3},\gamma_{4})$ can be 
legally applied indefinitely, resulting in a divergent game sequence. 
Then it suffices to show that from each fundamental position we can 
reach a position which meets condition ({\tt *}) using a sequence of 
legal node firings. 
The 
fundamental position $\omega_{1} = (1,0,0,0)$ meets condition 
({\tt *}).  
Now take fundamental position $\omega_{2} = (0,1,0,0)$ 
and apply the legal firing sequence $(\gamma_{2}, \gamma_{3}, 
\gamma_{4})$ to get the resulting position $(3,0,0,-2)$.  The latter 
meets condition ({\tt *}).  
Next take fundamental position $\omega_{3} = (0,0,1,0)$ 
and apply the legal firing sequence $(\gamma_{3}, 
\gamma_{4})$ to get the resulting position $(2,1,1,-2)$.  The latter 
meets condition ({\tt *}).  
For fundamental position $\omega_{4} = (0,0,0,1)$, 
apply the legal firing sequence $(\gamma_{4})$ to get the resulting 
position $(1,0,1,-1)$.  The latter 
meets condition ({\tt *}).

To finish our analysis of families of small cycles, 
we show why  
\hspace*{0.05in}
\parbox[c]{0.82in}{
\setlength{\unitlength}{0.7in}
\begin{picture}(1.1,1.2)
            \put(0,0.6){\circle*{0.075}}
            \put(0,0.1){\circle*{0.075}}
            \put(0.5,0.1){\circle*{0.075}}
            \put(0.5,1.1){\circle*{0.075}}
            \put(1,0.6){\circle*{0.075}}
            \put(0,0.6){\line(1,1){0.5}}
            \put(0,0.6){\line(0,-1){0.5}}
            \put(0,0.1){\line(1,0){0.5}}
            \put(1,0.6){\line(-1,1){0.5}}
            \put(1,0.6){\line(-1,-1){0.5}}
            \put(0.95,0.65){\vector(-1,1){0.1}}
            \put(0.55,1.05){\vector(1,-1){0.1}}
            \put(0.65,0.95){\vector(1,-1){0.1}}
           \end{picture}
} 
is not admissible.  The argument is entirely similar to the previous 
case.  Number the 
nodes with $\gamma_{1}$ as the North vertex and $\gamma_{2}$, $\gamma_{3}$, 
$\gamma_{4}$, and $\gamma_{5}$ in succession in the clockwise order 
around the cycle. 
We say an initial position $(a,b,c,d,e)$ meets 
condition ({\tt *}) if the following 
inequalities are satisfied: 
$b \geq 0$, $c \geq 0$, $d \geq 
0$, $e \leq 0$, and 
$a+e > 0$.  The firing sequence 
$(\gamma_{1},\gamma_{2},\gamma_{3},\gamma_{4},\gamma_{5})$ 
is easily seen to be 
legal from any such position.  The resulting position is 
$(a_{1},b_{1},c_{1},d_{1},e_{1})$ with $a_{1} = 4a + 2b + c + d + e$, $b_{1} = 
c$, $c_{1} = d$, $d_{1} = a+e$, and $e_{1} = -(3a+b+c+d+e)$.  
It is easy now to check 
that $(a_{1},b_{1},c_{1},d_{1},e_{1})$ also meets condition ({\tt *}).  (In fact, 
the inequalities $d_{1} > 0$ and $e_{1} < 0$ are now strict.)  
So from any position which meets condition ({\tt *}), the firing 
sequence  $(\gamma_{1},\gamma_{2},\gamma_{3},\gamma_{4},\gamma_{5})$ can be 
legally applied indefinitely, resulting in a divergent game sequence. 
Then it suffices to show that from each fundamental position we can 
reach a position which meets condition ({\tt *}) using a sequence of 
legal node firings. 
The 
fundamental position $\omega_{1} = (1,0,0,0,0)$ meets condition 
({\tt *}).  
Now take fundamental position $\omega_{2} = (0,1,0,0,0)$ 
and apply the legal firing sequence $(\gamma_{2}, \gamma_{3}, 
\gamma_{4},\gamma_{5})$ to get the resulting position $(2,0,0,0,-1)$.  The latter 
meets condition ({\tt *}).  
Next take fundamental position $\omega_{3} = (0,0,1,0,0)$ 
and apply the legal firing sequence $(\gamma_{3}, 
\gamma_{4}, \gamma_{5}, \gamma_{1}, \gamma_{2}, \gamma_{3}, \gamma_{4}, 
\gamma_{5})$ 
to get the resulting position $(5,0,0,0,-3)$.  The latter 
meets condition ({\tt *}).  
Next take fundamental position $\omega_{4} = (0,0,0,1,0)$ 
and apply the legal firing sequence $(\gamma_{4}, \gamma_{5}, \gamma_{1}, 
\gamma_{2}, \gamma_{3}, \gamma_{4}, 
\gamma_{5})$ 
to get the resulting position $(4,1,0,0,-3)$.  The latter 
meets condition ({\tt *}).  
For fundamental position $\omega_{4} = (0,0,0,0,1)$, 
apply the legal firing sequence $(\gamma_{5}, \gamma_{1}, 
\gamma_{2}, \gamma_{3}, \gamma_{4}, 
\gamma_{5})$  to get the resulting 
position $(4,0,1,0,-3)$.  The latter 
meets condition ({\tt *}).

This completes the proof of \NotMarsFriendlyCatalog.\hfill\QED

\vspace*{-0.2in}
\renewcommand{\baselinestretch}{1}

\end{document}